\documentclass[12pt,reqno]{amsart}
\usepackage{amscd}
\usepackage{amssymb}
\usepackage[mathscr]{eucal}
\usepackage{enumerate}

\makeatletter

\def\subsection{\@startsection{subsection}{2}%
  \z@{.5\linespacing\@plus.7\linespacing}{1pt}%
      {\normalfont\bfseries}}

\def\l@section{\@tocline{1}{0pt}{1pc}{}{\bfseries}}
\def\l@subsection{\@tocline{2}{0pt}{3.1em}{5pc}{}}

\makeatother

\allowdisplaybreaks[1]

\newtheorem{thm}{Theorem}[section]
\newtheorem{lem}[thm]{Lemma}
\newtheorem{cor}[thm]{Corollary}
\newtheorem{prop}[thm]{Proposition}

\theoremstyle{definition}

\newtheorem{rem}[thm]{Remark}
\newtheorem{defn}[thm]{Definition}

\theoremstyle{remark}

\numberwithin{equation}{section}

\newcommand{\C}{{\mathbb{C}}}
\newcommand{\R}{{\mathbb{R}}}
\newcommand{\Z}{{\mathbb{Z}}}

\newcommand{\N}{{\mathbb{N}}}

\DeclareMathOperator{\Ad}{Ad}
\DeclareMathOperator{\Tr}{Tr}

\DeclareMathOperator{\id}{id}

\DeclareMathOperator{\Proj}{Proj}
\DeclareMathOperator{\Projf}{Projf}

\DeclareMathOperator{\spa}{span}
\DeclareMathOperator{\ONB}{ONB}
\DeclareMathOperator{\Irr}{Irr}

\DeclareMathOperator{\Mor}{Mor}

\DeclareMathOperator{\Aut}{Aut}
\DeclareMathOperator{\Int}{Int}
\def\oInt{\ovl{\Int}}

\DeclareMathOperator{\supp}{supp}            
\DeclareMathOperator{\co}{co}                

\def\mB{\mathcal{B}}
\def\mC{\mathcal{C}}

\def\mF{\mathcal{F}}
\def\mK{\mathcal{K}}

\def\meC{\mathscr{C}}
\def\meF{\mathscr{F}}
\def\meH{\mathscr{H}}
\def\meI{\mathscr{I}}
\def\meJ{\mathscr{J}}
\def\meK{\mathscr{K}}
\def\meL{\mathscr{L}}

\def\meN{\mathscr{N}}

\def\meR{\mathscr{R}}
\def\meS{\mathscr{S}}
\def\meT{\mathscr{T}}
\def\meU{\mathscr{U}}

\def\al{\alpha}
\def\be{\beta}
\def\ga{\gamma}
\def\de{\delta}
\def\ka{\kappa}
\def\la{\lambda}
\def\ep{\epsilon}
\def\vep{\varepsilon}
\def\ph{{\phi}}
\def\ps{{\psi}}
\def\vph{\varphi}
\def\om{\omega}
\def\si{\sigma}
\def\ta{\tau}
\def\th{\theta}

\def\ze{\zeta}

\def\Si{\Sigma}
\def\Ph{\Phi}
\def\Ps{\Psi}

\def\De{\Delta}
\def\La{\Lambda}

\def\el{\ell}

\def\ovl{\overline}                 
\def\wdh{\widehat}
\def\wdt{\widetilde}
\def\opp{{\mathrm{opp}}}

\def\hA{\hat{A}}                    

\def\hal{{\hat{\al}}}

\def\hga{{\hat{\ga}}}

\def\hta{{\hat{\ta}}}

\def\hth{{\hat{\th}}}

\def\hps{{\hat{\ps}}}

\def\tV{\widetilde{V}}

\def\opi{{\ovl{\pi}}}                
\def\orho{{\ovl{\rho}}}
\def\osi{{\ovl{\si}}}

\def\subs{\subset}                   
\def\Subs{\Subset}
\def\setm{\setminus}
\def\nin{\notin}

\def\per{\perp}

\def\oti{\otimes}                    
\def\rti{\rtimes}                    

\def\col{\colon}
\def\ra{\hspace{-0.5mm}\rightarrow\!}

\def\sqr{\sqrt}                      

\def\btr{{\boldsymbol{1}}}                   
\def\bG{\mathbb{G}}
\def\bhG{\wdh{\mathbb{G}}}
\def\bhGG{\wdh{\mathbb{G}}\times \wdh{\mathbb{G}}^\opp}

\def\lhG{L^\infty(\bhG)}                    
\def\lhGG{L^\infty(\bhG\times\bhG^\opp)}    
\def\ltG{L^2(\bG)}
\def\lthG{L^2(\bhG)}                        
\def\lG{L^\infty(\bG)}                      

\def\IG{\Irr(\bG)}

\def\II{$\mbox{II}_{1}$}    

\title{Classification of minimal actions of a compact Kac algebra 
with amenable dual}

\author{Toshihiko Masuda}
\address{Toshihiko Masuda, 
Graduate School of Mathematics, Kyushu University, 
6-10-1 Hakozaki, Fukuoka, 812-8581, JAPAN}
\email{masuda@math.kyushu-u.ac.jp}

\author{Reiji Tomatsu}
\address{Reiji Tomatsu, Department of Mathematical Sciences,
University of Tokyo, 3-8-1 Komaba, Meguro, Tokyo, 153-8914, JAPAN}
\email{tomatsu@ms.u-tokyo.ac.jp}

\subjclass[2000]{Primary 46L10; Secondary 46L40}

\begin{document}

\begin{abstract}
We show the uniqueness of minimal actions of a compact Kac algebra 
with  amenable dual on the AFD factor of type II$_1$. 
This particularly implies the uniqueness of minimal 
actions of a compact group. 
Our main tools are a Rohlin type theorem, 
the 2-cohomology vanishing theorem, 
and the Evans-Kishimoto type intertwining argument.
\end{abstract}

\maketitle

\tableofcontents

\section{Introduction}

This paper presents 
the uniqueness of minimal actions of a compact Kac algebra 
with  amenable dual. More precisely, 
any two minimal actions of a compact Kac algebra with  amenable dual 
on the approximately finite dimensional (AFD) factor of type II$_1$ 
are conjugate. 
Note that every compact group action is 
a particular example of 
an action of such a compact Kac algebra. 
Here we say that an action is minimal 
if it has the full spectrum and the relative commutant 
of its fixed point algebra is trivial \cite{ILP}. 


After the completion of classification of discrete amenable group 
actions, it is natural to focus our 
attention to actions of continuous groups. 
Although 
it is very difficult to analyze actions of continuous groups in general, 
compact group actions have been extensively studied 
among them 
since the dual of a compact group is discrete. 
Indeed as emphasized in \cite{Wa1}, 
it gives us important insight about compact group actions 
to consider actions of compact group duals, namely, 
coactions of groups \cite{Nk-Tak}, 
or Roberts actions \cite{Rob-act}. 
Actions of compact abelian groups have been completely classified 
in \cite{JT} and \cite{KwT} by combining results 
in \cite{Oc1}, \cite{ST2} and \cite{KwST} with the Takesaki duality \cite{T-duality}. 
For compact non-abelian groups, 
classification of all actions is still very far from completion, 
but a few kinds of actions, such as ergodic actions and 
minimal actions, have been studied. 
For example, ergodic actions have been studied 
in \cite{OPT} for the abelian case, 
and A. Wassermann has dealt with general ergodic actions 
in \cite{Wa2}, \cite{Wa3} and \cite{Wa4}. 
In particular, he has finished the classification of 
ergodic actions of $SU(2)$. 

It is a central theme to classify all minimal actions 
in the study of compact group actions. 
The notion of minimality corresponds 
to outerness of actions in discrete case. 
So far several attempts have been made for this classification problem. 
In an unpublished work \cite{Oc2}, 
Ocneanu has announced the uniqueness of minimal actions of a compact group 
by developing the method used in \cite{Oc1}. 
A different approach has been proposed by 
S. Popa and A. Wassermann in \cite{PW}, 
which is based on the classification of subfactors by Popa \cite{Po-amen}. 
They have applied the main theorem of \cite{Po-amen} 
to Wassermann's subfactors \cite{Wa1}, 
and concluded 
the uniqueness of minimal actions of a compact Lie group. 
Unfortunately, 
the details of the both adorable 
theories have not been available to the authors. 

We explain our approach to this problem. 
Most part of this paper is devoted to 
classification of centrally free actions of 
amenable discrete Kac algebras. 
We can utilize the framework of amenable discrete Kac algebras 
for treating the duals of compact groups. 
The uniqueness of minimal actions follows 
through the duality between compact Kac algebras 
and discrete Kac algebras. 
Throughout this paper, 
the most fundamental tools are 
the ultraproducts and central sequence technique. 
Here one notes that 
actions of Kac algebras never preserve central 
sequences. 
To overcome this difficulty, 
we mainly treat approximately inner actions. 
Then we can use the central sequence technique
as demonstrated in \cite{Ma1}.

The first half of our arguments 
is similar to those in \cite{Con-auto}, \cite{Oc1}, 
that is, we formulate a Rohlin type theorem 
for actions of discrete amenable Kac algebras, 
and show approximate 
1-cohomology vanishing and 2-cohomology vanishing by the Shapiro type 
argument. 
Combining these results with the
stability of minimal actions (cf. \cite{Wa2}), 
we show that every minimal action is dual. 

The traditional method 
for classification theory of actions on von Neumann algebras 
is the model action splitting argument, 
that is, we construct the model action as 
an infinite tensor product type action at first, 
and pull out pieces of the model 
action from a given action by means of the Rohlin type theorem and 
cohomology vanishing. 
It is not so difficult 
to construct the infinite tensor product type model action for 
cyclic groups or finite groups as in \cite{Con-auto}, \cite{Co-peri} and
\cite{Jo1}, 
but it is necessary to apply the paving 
theorem of D. Ornstein and B. Weiss \cite{Ornstein1}, \cite{Ornstein2} 
to construct the model action for a general discrete amenable group. 
Because of use of paving theorem, 
Ocneanu's  Rohlin type theorem \cite{Oc1} takes a complicated form, 
and it seems difficult to 
generalize it to a discrete amenable Kac algebra
case. 
Of course, it is easy to construct the model action for coactions of 
finite groups, and the classification given in \cite{Ma1} is based on 
the model action splitting argument.

So we do not take the model action splitting argument 
in the final stage of classification. 
Instead of that, we use the 
intertwining argument, which was initiated by 
D. E. Evans and A. Kishimoto in
\cite{EvKi}, and has been further developed in \cite{Naka} and
\cite{Iz-roh1} for group actions on $C^*$-algebras. 
It also works for von Neumann algebras as is shown in \cite{Ma-inter}, 
enables us to avoid using a paving theorem, 
and makes our arguments simpler than those in \cite{Oc1}. 

\textbf{Acknowledgements}. 
The authors are grateful to Masaki Izumi and Yasuyuki Kawahigashi 
for permanent encouragement and fruitful discussions. 
The second named author was supported in part by Research Fellowship for 
Young Scientists of the Japan Society for the Promotion of 
Science.

\section{Preliminaries}

\subsection{Notations}
We treat only separable von Neumann algebras and Kac algebras 
except for ultraproduct von Neumann algebras. 
We denote by $\meR_0$ the AFD factor of type \II. 
Let $M$ be a von Neumann algebra. 
For a subset $S\subs M$, we denote by $W^*(S)$ the von Neumann 
subalgebra of $M$ generated by $S$. 
We denote 
the sets of unitaries, projections, positive elements, 
the unit ball and the center of $M$  
by $U(M)$, $\Proj(M)$, $M_+$, $M_1$ and $Z(M)$, respectively. 
For a weight $\ph$ on $M$ and $x~\in M$, 
we set $|x|_\ph=\ph(|x|)$, $\|x\|_\ph=\ph(x^*x)^{1/2}$ and 
$\|x\|_\ph^\sharp=2^{-1/2}(\ph(x^*x)+\ph(xx^*))^{1/2}$. 
For each $a\in M$ and $\ph\in M_*$, 
$a\ph, \ph a\in M_*$ are defined by 
$a\ph(y)=\ph(ya)$ and $\ph a(y)=\ph(ay)$ for $y\in M$. 
Set $[\ph,x]=\ph x-x\ph$. 
For $x\in M_+$, $s(x)$ denotes its support projection. 
For von Neumann algebras $M$ and $N$, 
we denote by $\Mor(M,N)$ 
the set of unital faithful normal $*$-homomorphisms 
from $M$ to $N$. 
For von Neumann algebras $M$ and $N$, 
$M\oti N$ is the tensor product von Neumann algebra. 
We denote by $\Re z$ the real part of a complex number $z$. 
We write $A\Subs B$ when $A$ is a finite subset of $B$. 
Let $X$ be a linear space and consider the $n$-fold 
tensor product $X^{\otimes n}$. 
Then the symmetric group $S_n$ acts on $X^{\otimes n}$ canonically. 
The action is written by $\si$. 
Note that we also use the symbol $\si$ for an irreducible representation of a 
compact Kac algebra. 

\subsection{Quick review of theory of Kac algebras}
Our basic references on theory of Kac algebras are 
\cite{BaSk} and \cite{ES}. 
A \textit{compact Kac algebra} is a triple 
$\bG=(A,\de,h)$ where 
$A$ is a von Neumann algebra, 
$\de\in\Mor(A,A\oti A)$ is a coproduct 
and $h$ is a faithful invariant tracial state, 
that is, they satisfy 
\[
(\de\oti\id)\circ\de=(\id\oti\de)\circ\de, 
\]
\[
(\th\oti h)(\de(a))=\th(1)h(a)=(h\oti\th)(\de(a))
\quad 
\mbox{for all}\ 
\th\in A_*,\  a\in A. 
\]
Let $\si\in \Aut(A\oti A)$ be the flip automorphism. 
Then the map $\si\circ\de$ is denoted by $\de^\opp$. 
Let $\{\pi_h, \hat{1}_h, H_h\}$ be the GNS representation of $A$ 
with respect to the state $h$, 
which means $H_h$ is a Hilbert space, 
$\pi_h\in\Mor(A, B(H_h))$ and $\hat{1}_h$ is the GNS cyclic vector. 
We always regard $A$ as a subalgebra of $B(H_h)$ via the map $\pi_h$. 
The canonical tracial weight on $B(H_h)$ is denoted by $\Tr$. 
We define a unitary $V\in B(H_h\oti H_h)$ by 
\[
V^* (x\hat{1}_h\oti y\hat{1}_h)=\de(y) (x\hat{1}_h\oti \hat{1}_h)
\quad 
\mbox{for}
\ 
x,\ y\in A. 
\]
Then $V$ satisfies the following pentagonal equality, and is called 
a \textit{multiplicative unitary}.
\[
V_{12}V_{13}V_{23}=V_{23}V_{12}. 
\]
By definition, $V\in A\oti B(H_h)$. 
There exists the \textit{antipode} $I$ on $A$ 
which is an antiautomorphism on $A$ satisfying $(I\oti\id)(V)=V^*$. 
Define $\hA$ by the $\si$-weak closure of the linear space 
$\{(\th\oti\id)(V)\mid \th\in A_*\}$ 
and 
a map $\De(x)=V(x\oti1)V^*$ for $x\in \hA$. 
Then $\hA$ is actually a von Neumann algebra 
and $\De\in\Mor(\hA,\hA\oti\hA)$. 
Let $\vph$ be the Planchrel weight on $\hA$ induced by a compact 
Kac algebra $\bG$. 
In fact, $\vph=\Tr$ holds on $\hA$. 
Then the triple $\bhG=(\hA,\De,\vph)$ is a 
\textit{discrete Kac algebra}. 
The tracial weight $\vph$ is invariant for $\De$, that is, 
\[
(\th\oti\vph)(\De(x))=\th(1)\vph(x)=(\vph\oti\th)(\De(x))
\quad \mbox{for all}\ 
\th\in \hA_*,\ x\in \hA_+. 
\]

For a more convenient description of von Neumann algebras $A$ and $\hA$, 
we make use of unitary representations of $\bG$. 
Let $K$ be a Hilbert space and $v\in U(A\oti B(K))$. 
The pair $\pi=(v,K)$ is called a \textit{unitary representation} 
of $\bG$ 
if $(\de\oti\id)(v)=v_{13}v_{23}$. 
The unitary representation $\btr=(1,\C)$ is called the 
\textit{trivial representation}. 
For two unitary representations $\pi=(v,K)$ and 
$\si=(v',K')$, 
an element $S\in B(K, K')$ is called an 
\textit{intertwiner} from $\pi$ to $\si$ 
if $(1\oti S)v=v'(1\oti S)$. 
The set of intertwiners from $\pi$ to $\si$ is a linear space and 
denoted by $(\pi,\si)$. 
If $(\pi,\si)$ contains an isometry, then we write 
$\pi\prec\si$. 
If $(\pi,\si)$ contains a unitary, 
$\pi$ and $\si$ are said to be equivalent and 
we write $\pi\sim\si$. 
Of course $\pi\prec\si$ and $\si\prec\pi$ implies $\pi\sim\si$. 
We define the tensor product representation 
by $\pi\cdot\si=(v_{12}v_{13}', K\oti K')$. 

For a unitary representation $\pi=(v,K)$, 
we define the \textit{conjugate unitary representation} 
$\pi^c=(v^c,\ovl{K})$ as follows. 
Let $\ovl{K}$ be the conjugate Hilbert space of $K$ 
with the conjugation map $j\col K\ra \ovl{K}$. 
Define the transpose map $t\col B(K)\ra B(\ovl{K})$ 
by $t(x)=jx^*j^{-1}$ for all $x\in B(K)$. 
Then set $v^c=(I\oti t)(v)$. 
The relation $\de\circ I=(I\oti I)\circ\de^\opp$ implies 
that $\pi^c$ is a unitary representation. 

For 
$\pi=(v,K)$, the set $(\pi,\pi)$ is a $C^*$-subalgebra 
of $B(K)$. 
We say that $\pi$ is \textit{irreducible} 
if $(\pi,\pi)=\C$. 
The irreducibility of $\pi=(v,K)$ implies finite dimensionality of $K$. 
We write $d_\pi$ for $\dim K$. 
Let $\si$ be another unitary representation. 
Then the intertwiner space $(\pi,\si)$ is a Hilbert space 
with the inner product $(S,T)1_\pi=T^*S$ for all $S,T\in (\pi,\si)$, 
where $1_\pi$ denotes the unit of $B(K)$. 
When we fix an orthonormal basis of $(\pi,\si)$, 
we denote it by $\ONB(\pi,\si)$. 

We denote by $\IG$ 
the set of equivalence classes of all irreducible unitary representations. 
We denote by $[\pi]$ 
the equivalence class of a unitary representation $\pi$. 
The set $\IG$ has the conjugation operation defined by 
$\ovl{[\pi]}=[\pi^c]$ for all $[\pi]\in\IG$. 
We fix a representative $\pi=(v_\pi, H_\pi)$ for each $[\pi]\in\IG$ 
as follows. 
If $[\pi]\neq\ovl{[\pi]}$, 
we take representatives $\pi$ and $\pi^c$ 
for $[\pi]$ and $\ovl{[\pi]}$, respectively. 
In this case, we often write $\ovl{\pi}=(v_\opi,H_\opi)$ 
for $\pi^c=(v_\pi^c, \ovl{H_\pi})$ with the conjugation 
$j_\pi\col H_\pi\ra \ovl{H_\pi}$. 
If $[\pi]$ is self-conjugate, that is, $[\pi]=\ovl{[\pi]}$, 
then we take a representative $\pi$ for $[\pi]$. 

Take a self-conjugate $[\pi]\in\IG$. 
Then there exists a unitary $\nu_\pi\in B(\ovl{H_\pi},H_\pi)$ 
such that 
$v_\pi=(1\oti \nu_\pi)v_\pi^c (1\oti \nu_\pi^*)$. 
Taking the conjugation of the both sides, 
we have 
$v_\pi^c=(1\oti\ovl{\nu_\pi})v_\pi(1\oti\ovl{\nu_\pi}^*)$, 
where $\ovl{\nu_\pi}=j_\pi\nu_\pi j_\pi\in B(H_\pi,\ovl{H_\pi})$. 
Hence $\nu_\pi \ovl{\nu_\pi}$ is in $(\pi,\pi)$ 
and we have $\nu_\pi \ovl{\nu_\pi}=\ep_\pi 1_\pi$ or equivalently 
$\ovl{\nu_\pi}\nu_\pi=\ep_\pi \ovl{1_\pi}$ 
for some $\ep_\pi\in \C$ with $|\ep_\pi|=1$. 
Taking the conjugation of 
$\ovl{\nu_\pi}\nu_\pi=\ep_\pi \ovl{1_\pi}$, 
we have $\nu_\pi \ovl{\nu_\pi}=\ovl{\ep_\pi} 1_\pi$, 
and hence $\ep_\pi=\ovl{\ep_\pi}$, that is, $\ep_\pi=\pm1$. 
We assign $\ep_\pi=1$ for a nonself-conjugate $[\pi]\in\IG$. 

Now for each $[\pi]\in\IG$, 
fix a finite index set $I_\pi$ with $|I_\pi|=d_\pi$ 
and 
an orthonormal basis $\{\vep_{\pi_i}\}_{i\in I_\pi}$ of $H_\pi$. 
Then we introduce another orthonormal basis of $H_\opi$  
$\{\vep_{\opi_{\ovl{i}}}\}_{i\in I_\pi}$ as follows. 
If $[\pi]\neq \ovl{[\pi]}$, set 
$\vep_{\opi_{\ovl{i}}}=j_\pi \vep_{\pi_i}\in \ovl{H_\pi}$. 
If $[\pi]=\ovl{[\pi]}$, then set 
$\vep_{\opi_{\ovl{i}}}=\nu_\pi j_\pi \vep_{\pi_i}\in H_\pi$. 
For each $[\pi]\in\IG$, we define an isometric intertwiner 
$T_{\opi,\pi}\in (\btr, \opi\cdot\pi)$ by 
\[
T_{\opi,\pi}=\sum_{i\in I_\pi} \frac{1}{\sqr{d_\pi}}
\vep_{\opi_{\ovl{i}}}\oti \vep_{\pi_i}. 
\]
We claim the following equality 
\[
T_{\opi,\pi}=\ep_\pi\sum_{i\in I_\opi} \frac{1}{\sqr{d_\pi}}
\vep_{\opi_i}\oti \vep_{\pi_{\ovl{i}}}. 
\]
For a nonself-conjugate $[\pi]$, it is trivial. 
We verify the equality for a self-conjugate $[\pi]$. 
By invariance of the summation in changing bases, 
we have 
\begin{align*}
T_{\opi,\pi}
=&\,
\sum_{i\in I_\pi} \frac{1}{\sqr{d_\pi}}
\vep_{\opi_{\ovl{i}}}\oti \vep_{\pi_i}
\\
=&\,(\nu_\pi\oti1)
\sum_{i\in I_\pi} \frac{1}{\sqr{d_\pi}}
j_\pi \vep_{\pi_i}\oti\vep_{\pi_i}
\\
=&\,(\nu_\pi\oti1)
\sum_{i\in I_\pi} \frac{1}{\sqr{d_\pi}} 
j_\pi \nu_\pi j_\pi \vep_{\pi_i}\oti  
\nu_\pi j_\pi \vep_{\pi_i}
\\
=&\,(\nu_\pi \ovl{\nu_\pi}\oti1)
\sum_{i\in I_\pi} \frac{1}{\sqr{d_\pi}} 
\vep_{\pi_i}\oti 
\vep_{\pi_{\ovl{i}}}
\\
=&\,
\ep_\pi\sum_{i\in I_\pi} \frac{1}{\sqr{d_\pi}} 
\vep_{\pi_i}\oti 
\vep_{\pi_{\ovl{i}}}. 
\end{align*}
On $T_{\opi,\pi}$ we have the following equalities
\[
(1_\pi\oti T_{\opi,\pi}^*)(T_{\pi,\opi}\oti1_\pi)
=
\frac{\ep_\pi}{d_\pi}1_\pi
=
(T_{\pi,\opi}^*\oti1_\pi)(1_\pi\oti T_{\opi,\pi}). 
\]
From now on, we simply write $\pi$ for $[\pi]$. 
Hence $\pi$ means an element of $\IG$ and an irreducible 
unitary representation. 

We take the systems of matrix units $\{e_{\pi_{i,j}}\}_{i,j\in I_\pi}$ and 
$\{e_{\opi_{\ovl{i},\ovl{j}}}\}_{i,j\in I_\pi}$ 
for $B(H_\pi)$ and $B(H_\opi)$ 
coming from the bases $\{\vep_{\pi_i}\}_{i\in I_\pi}$ and 
$\{\vep_{\opi_{\ovl{i}}}\}_{i\in I_\pi}$, respectively. 
We decompose $v_{\pi}$ as 
\[
v_\pi=\sum_{i,j\in I_\pi} v_{\pi_{i,j}}\oti e_{\pi_{i,j}}.
\] 
Then the elements $\{v_{\pi_{i,j}}\mid i,j\in I_\pi,\ \pi\in\IG\}$ 
are linearly independent and the linear span of them is $\si$-weakly dense 
in $A$. 
In fact the following orthogonal relations of them hold. 
For all $\pi, \rho\in\IG$, $i, j\in I_\pi$ and $k, \el\in I_\rho$, 
\[
h(v_{\pi_{i,j}}^*v_{\rho_{k,\el}})
=d_\pi^{-1}\de_{\pi,\rho}\de_{i,k}\de_{j,\el}. 
\] 
For each $\pi\in\IG$ and $i,j\in I_\pi$, 
set $\vep_{\pi_{i,j}}= d_\pi h v_{\pi_{i,j}}^* \in A_*$ 
and put $f_{\pi_{i,j}}=(\vep_{\pi_{i,j}}\oti\id)(V)$. 
Then 
$f_{\pi_{i,j}} f_{\rho_{k,\el}}=\de_{\pi,\rho}\de_{j,k}f_{\pi_{i,\el}}$ 
and $f_{\pi_{i,j}}^*=f_{\pi_{j,i}}$ hold. 
Since $\{f_{\pi_{i,j}}\}_{i,j\in I_\pi, \pi\in\IG}$ generates $\hA$, 
we have identification by putting $f_{\pi_{i,j}}=e_{\pi_{i,j}}$, 
\[
\hA=\bigoplus_{\pi\in\IG}B(H_\pi). 
\]
By this identification, we have
\[
\vph=\bigoplus_{\pi\in\IG} d_\pi \Tr_\pi
\]
and
\[
V=\bigoplus_{\pi\in\IG}\sum_{i,j\in I_\pi} v_{\pi_{i,j}}\oti e_{\pi_{i,j}}, 
\]
where $\Tr_\pi$ is the non-normalized trace on $B(H_\pi)$. 
We also use the tracial state $\ta_\pi$ on $B(H_\pi)$. 
For $x\in\hA$, we often use the notation $x_\pi=x(1\oti1_\pi)$. 
The \textit{support} of $x\in\hA$ is the subset of $\IG$ which consists of 
$\pi\in\IG$ such that $x_\pi\neq0$. 
The support of $x\in\hA$ is denoted by $\supp(x)$. 
We write $\hA_0$ for the set of finitely supported elements of $\hA$. 
We denote by $\Projf(\hA)$ the set of projections in $\hA_0$. 

The representation $\hA\subs B(H_h)$ is standard as is seen below. 
Let $(\pi_\vph,\La_\vph,H_\vph)$ be the GNS representation of $\hA$ 
with respect to $\vph$. 
We define a unitary $H_\vph\ra H_h$ which maps 
$\La_\vph(e_{\pi_{i,j}})$ to $d_{\pi}v_{\pi_{i,j}}\hat{1}_h$ for 
all $\pi\in\IG$ and $i,j\in I_\pi$. 
Then this unitary intertwines $\pi_\vph$ 
and the identity representation on $H_h$. 
Thus we always identify $H_h$ with $H_\vph$. 
Then we obtain 
\[
V (\La_\vph(x)\oti\La_\vph(y))=\La_{\vph\oti\vph}(\De(x)(1\oti y))
\quad 
\mbox{for all}
\ x, y\in \hA_0. 
\]
We also use a unitary $W$ defined by 
\[
W^* (\La_\vph(x)\oti\La_\vph(y))=\La_{\vph\oti\vph}(\De(y)(x\oti 1))
\quad 
\mbox{for }
\ x, y\in \hA_0. 
\]
The unitary $W$ also satisfies the pentagonal equality. 
We call $V$ and $W$ the \textit{right regular representation} and the {\it left
regular representation} of $\bhG$, respectively. 
Then $V\in A\oti \hA$ and $W\in \hA\oti A'$. 
By definition of $V$ and $W$, 
we have $V(x\oti1)V^*=\De(x)=W^*(1\oti x)W$ for all $x\in\hA$. 

We give a description of $\De$ in terms of intertwiners. 
We write ${}_\pi\De_\rho(x)$ for $\De(x)(1_\pi\oti1_\rho)$ for each 
$\pi,\rho\in\IG$. 
Let $\pi,\rho, \si\in \IG$ and then
\[
{}_\pi\De_\rho(x)S=S x_\si
\quad \mbox{for all}\ 
x\in\hA, \ S\in(\si,\pi\cdot\rho), 
\]
in particular,
\[
{}_\opi\De_\pi(x) T_{\opi,\pi}=x_\btr T_{\opi,\pi} 
\quad \mbox{for all}\ 
x\in\hA.
\]
By complete decomposability of $\pi\cdot\rho$, 
\[
1_\pi\oti1_\rho
=\sum_{\si\prec\pi\cdot\rho}
\sum_{S\in\ONB(\si,\pi\cdot\rho)}SS^*, 
\]
where the index $\si$ runs in $\IG$. 
This equality implies that 
\[
{}_\pi\De_{\rho}(x)
=
\sum_{\si\prec\pi\cdot\rho}
\sum_{S\in\ONB(\si,\pi\cdot\rho)}
S x_\si S^*.
\]
Putting $x=e_\btr$, we have 
\begin{align*}
{}_\opi \De_\pi(e_\btr)
=&\,
T_{\opi,\pi}T_{\opi,\pi}^*
\\
=&\,
\sum_{i,j\in I_\pi}
\frac{1}{d_\pi}e_{\opi_{\ovl{i},\ovl{j}}}\oti e_{\pi_{i,j}}
\\
=&\,
\sum_{i,j\in I_\opi}
\frac{1}{d_\pi}e_{\opi_{i,j}}\oti e_{\pi_{\ovl{i},\ovl{j}}}. 
\end{align*}

Let $\pi,\rho$ and $\si$ be elements in $\IG$. 
Then there exists a conjugate unitary map
$\wdt{\cdot}$ from $(\osi,\orho\cdot\opi)$ to $(\si,\pi\cdot\rho)$ 
as is defined below. 
For $S\in(\osi,\orho\cdot\opi)$, 
\[
\wdt{S}^*
=\sqr{d_\si d_\pi d_\rho}
(1_\si\oti T_{\orho,\rho}^*)
(1_\si\oti1_\orho\oti T_{\opi,\pi}^*\oti1_\rho)
(1_\si \oti S\oti1_\pi\oti1_\rho)(T_{\si,\osi}\oti1_\pi\oti1_\rho).
\]

The modular conjugations $J_h$ and $J_\vph$ are defined by
\[J_h x\hat{1}_h=x^*\hat{1}_h
,\quad
J_\vph \La_\vph(y)=\La_\vph(y^*)
\quad 
\mbox{for}\ 
x\in A,\ y\in\hA. 
\]
Set a unitary $U=J_h J_\vph=J_\vph J_h$. 
For $x\in \hA$, define $\hat{I}(x)=J_h x^* J_h$. 
Then $\hat{I}$ is an antiautomorphism of $\hA$, and called the 
antipode of $\bhG$. 
The equality $(I\oti \hat{I})(V)=V$ yields 
$\hat{I}(1_\pi)=1_\opi$. 
We often write $\ovl{x}$ for $\hat{I}(x)$ for $x\in Z(\hA)$. 

Since we want to illustrate Kac algebras as function algebras 
on noncommutative spaces, we prepare the new notations 
\[
\lG=A,\quad \lhG=\hA, \quad 
\ltG=H_h=H_\vph=\lthG. 
\]

Define the coproduct $\De^\opp=\si\circ\De$, where 
$\si\in\Aut(\hA\oti\hA)$ is the flip automorphism 
and then the triple $\bhG^\opp=(\hA,\De^\opp,\vph)$ is called 
the \textit{opposite discrete Kac algebra} of $\bhG$. 
We define a discrete Kac algebra 
$\bhGG=(\lhGG, \De_{\bhGG},\vph_{\bhGG})$ as follows. 
The von Neumann algebra $\lhGG$ is $\lhG\oti \lhG$. 
The coproduct $\De_{\bhGG}$ is given by 
$\De_{\bhGG}(x\oti y)=\De(x)_{13}\De^\opp(y)_{24}$. 
The invariant tracial weight $\vph_{\bhGG}$ is equal to 
$\vph\oti\vph$. 

The map $(\De\oti\id)\circ\De=(\id\oti\De)\circ\De$ is often denoted by 
$\De^{(2)}$. 
For subsets $\mF$ and $\mK\subs\IG$, the subset $\mF\cdot \mK\subs\IG$ 
is defined as 
$\{\pi\in\IG\mid \pi\prec \rho\cdot\si 
\ \mbox{for some} \ \rho\in\mF, \si\in \mK\}$.

\subsection{Amenability of a Kac algebra} 
The amenability of Kac algebras has been studied by many authors. 
We refer readers to \cite{Ru} and the references therein. 
A discrete Kac algebra $\bhG=(\lhG,\De,\vph)$ 
is called \textit{amenable} if 
there exists a left invariant state $m\in \lhG^*$, 
here the left invariance means 
$m((\th\oti \id)(\De(x)))=\th(1)m(x)$ for all $\th\in \lhG_*$ 
and $x\in \lhG$. 
When $\bG$ comes from a compact group, 
then $\bhG$ is amenable by \cite[Theorem 4.5]{Ru}, 
which is directly proved by using 
the Kakutani-Markov fixed point theorem 
due to the cocommutativity of $\bhG$. 
The amenability is also characterized by existence of 
projections with approximate invariance. 
We follow \cite{Oc2} 
for a notion of approximate invariance of a projection. 

\begin{defn}\label{appinv}
Let $F\in\Projf(\lhG)$ and $\vep>0$. 
A projection $S\in \Projf(\lhG)$ 
is said to be $(F,\vep)$\textit{-invariant} if we have  
\[|(F\oti1)\De(S)-F\oti S|_{\vph\oti\vph} 
< \vep |F|_\vph |S|_\vph
.\] 
\end{defn}

We say that a Kac algebra $\bhG$ satisfies 
the \textit{F\o lner condition} if 
for any $F\in\Projf(Z(\lhG))$ and $\vep>0$, 
there exists 
an $(F,\vep)$-invariant $K\in\Projf(Z(\lhG))$. 
Thanks to \cite[Theorem 4.5]{Ru}, 
a Kac algebra $\bhG$ is amenable if and only if 
it satisfies the F\o lner condition. 
Moreover for any $F$, $\vep>0$ as above, 
we can take an $(F,\vep)$-invariant 
$K\in\Projf(Z(\lhG))$ such that $K\geq e_\btr$. 
We sketch the proof of this fact as follows. 
Take an $(F,\vep)$-invariant $\wdt{K}\in\Projf(Z(\lhG))$. 
Take $\pi\in\IG$ such that $\wdt{K}1_\pi\neq 0$. 
Set $x=(\id\oti\ta_\pi)(\De(\wdt{K}))\leq1$. 
Actually $x$ is in $Z(\lhG)$ and 
it has sufficient invariance for $F$. 
Then applying the Day-Namioka-Connes trick 
\cite[Theorem1.2.2]{Co}, 
we obtain a central projection $K$ which 
is sufficiently invariant for $F$. 
Moreover since $x e_\btr=e_\btr$, $e_\btr$ still remains 
after the trick, that is, $K\geq e_\btr$. 

\subsection{Actions and cocycle actions}
Let $M$ be a von Neumann algebra and 
$\bhG=(\lhG,\De,\vph)$ a discrete Kac algebra. 
Let $\al\in \Mor(M, M\oti\lhG)$ and 
$u$ a unitary in $M\oti\lhG\oti\lhG$. 
The pair $(\al,u)$ is called a \textit{cocycle action} of $\bhG$ 
on $M$ 
if we have the following three conditions 
\begin{enumerate}

\item
$(\al\oti\id)\circ\al=\Ad u\circ (\id\oti\De)\circ\al$, 

\item
$(u\oti1)(\id\oti\De\oti\id)(u)
=(\al\oti\id\oti\id)(u)(\id\oti\id\oti\De)(u)$, 

\item 
$u_{\btr,\pi}=1\oti e_\btr\oti 1_\pi$, 
$u_{\pi,\btr}=1\oti 1_\pi\oti e_\btr$ for all $\pi\in\IG$. 

\end{enumerate} 
By definition, $\al_\btr=\id$. 
The unitary $u$ is called a \textit{2-cocycle}. 
If $u=1$, we say that $\al$ is an \textit{action}. 
A \textit{perturbation} $\wdt{u}$ of $u$ by $v\in U(M\oti\lhG)$ 
is defined by 
\[
\wdt{u}=(v\oti1)(\al\oti\id)(v)u(\id\oti\De)(v^*). 
\]
Then the pair $(\Ad v\circ\al,\wdt{u})$ is a cocycle action 
perturbed by $v$. 
If $\wdt{u}=1$, we say that $u$ is a \textit{2-coboundary}. 
For an action $\al$, $v\in U(M\oti \lhG)$ 
is called a \textit{1-cocycle} or \textit{$\al$-cocycle} 
if we have 
$(v\oti1)(\al\oti\id)(v)=(\id\oti\De)(v)$. 
Note that 
if $v$ perturbs an action $(\al,1)$ to the action $(\Ad v\circ\al,1)$, 
$v$ is an $\al$-cocycle. 
A perturbation of an $\al$-cocycle $v$ by 
$w\in U(M)$ is $\wdt{v}=(w\oti1)v\al(w^*)$. 
If $\wdt{v}=1$, $v$ is a \textit{1-coboundary}. 
We say that a cocycle is small if it is close to 1. 

To simplify notations, 
we often omit the symbol $\oti\id$ after 
symbols of cocycle actions such as $\al$, $\be$, $\ga$ and so on. 
For example, we write $\al(u)$ for $(\al\oti\id\oti\id)(u)$. 

Decompose $u$ as 
\[
u=\sum_{\pi,\rho\in\IG} \sum_{i,j\in I_\pi}\sum_{k,\el\in I_\rho} 
u_{\pi_{i,j},\rho_{k,\el}}\oti e_{\pi_{i,j}}\oti e_{\rho_{k,\el}}.
\]
Then each element $u_{\pi_{i,j},\rho_{k,\el}}$ is called the entry of $u$. 
We simply say that an element $x\in M$ commutes with $u$ 
if it does with all entries of $u$. 

\begin{defn}
Let $M$ be a von Neumann algebra and $K$ a finite dimensional 
Hilbert space. 
Let $\al\in \Mor(M,M\oti B(K))$. 
A faithful normal unital completely positive map 
$\Ph\col M\oti B(K)\ra M$ is called a 
\textit{left inverse} of $\al$ 
if it satisfies $\Ph\circ\al=\id$. 
\end{defn}

We denote by $\Mor_0(M,M\oti B(K))$ the subset consisting of 
an element in $\Mor(M, M\oti B(K))$ with a left inverse. 
Note that if $\Ph$ is a left inverse of $\al$, 
$\al\circ\Ph$ is a conditional expectation from 
$M\oti B(K)$ to $\al(M)$. 
Although a left inverse is not uniquely determined in general, 
we always treat the following left inverses for cocycle actions. 

\begin{defn}
For a cocycle action $(\al,u)$ of $\bhG$ on $M$ and 
$\pi\in\IG$, we define 
the map $\Ph_\pi^\al\col M\oti B(H_\pi)\ra M$ by 
\[
\Phi_\pi^\al(x)=(1\oti T_{\opi,\pi}^*)
u_{\opi,\pi}^*(\al_{\opi}\oti\id)(x)
u_{\opi,\pi}
(1\oti T_{\opi,\pi})
\quad \mbox{for all}\  x\in M\oti B(H_\pi). 
\]
\end{defn}
We simply write $\Ph_\pi$ for $\Ph_\pi^\al$ if no confusion arises.
The next lemma shows that $\Ph_\pi$ is actually a left inverse 
of $\al_\pi$. 

\begin{lem}\label{lem: left inverse}
Let $(\al,u)$ be a cocycle action of $\bhG$ on $M$ 
and $\pi\in\IG$. 

\begin{enumerate}

\item 
The map $\Phi_\pi$ is faithful, normal, 
unital, and completely positive.

\item 
$\Phi_\pi(\al_\pi(a)b\al_\pi(c))=a\Phi_\pi(b)c$ 
for all $a,c\in M$ and $b\in M\oti B(H_\pi)$. 

\item 
For any $x\in M$, the following equality holds.
\[
(\Ph_\pi\oti\id)(u_{\pi,\opi}
(x\oti T_{\pi,\opi}T_{\pi,\opi}^*)
u_{\pi,\opi}^*)
=
d_\pi^{-2}\al_{\opi}(x)
.\]
\end{enumerate}
\end{lem}

\begin{proof}
(1). 
It is clear that $\Ph_\pi$ is a normal unital completely positive map. 
Assume 
$\Ph_\pi(x^*x)=0$, 
and  then we have 
$(\al_{\opi}\oti\id)(x)u_{\opi,\pi}
(1\oti T_{\opi,\pi})=0$. 
Applying $\al_\pi$ to the first leg, we have 
\[
\big{(}
(\al_\pi\oti\id)\circ\al_{\opi}\oti\id
\big{)}(x)
(\al_\pi\oti\id\oti\id)(u_{\opi,\pi})
(1\oti 1_\pi\oti T_{\opi,\pi})=0.
\]
The equality implies 
\begin{align*}
0=&\,
(1\oti T_{\pi,\opi}^*\oti 1_\pi)
(u_{\pi,\opi}^*\oti1_\pi)
\big{(}
(\al_\pi\oti\id)\circ\al_{\opi}\oti\id
\big{)}(x)
\\
&\quad\cdot
(\al_\pi\oti\id\oti\id)(u_{\opi,\pi})
(1\oti 1_\pi\oti T_{\opi,\pi})
\\
=&\,
(1\oti T_{\pi,\opi}^*\oti 1_\pi)
(\id\oti_{\pi}\De_{\opi}\oti\id)
\big{(}(\al\oti\id)(x)\big{)}
\\
&\quad\cdot
(u_{\pi,\opi}^*\oti1_\pi)
(\al_\pi\oti\id\oti\id)(u_{\opi,\pi})
(1\oti 1_\pi\oti T_{\opi,\pi})
\\
=&\,
x (1\oti T_{\pi,\opi}^*\oti 1_\pi)
(u_{\pi,\opi}^*\oti1_\pi)
(\al_\pi\oti\id\oti\id)(u_{\opi,\pi})
(1\oti 1_\pi\oti T_{\opi,\pi})
\\
=&\,
x (1\oti T_{\pi,\opi}^*\oti 1_\pi)
(\id\oti_{\pi}\De_{\opi}\oti\id)
(u)
(\id\oti\id\oti_{\opi}\De_{\pi})(u^*)
(1\oti 1_\pi\oti T_{\opi,\pi})
\\
=&\,
x u_{\btr,\pi}(1\oti T_{\pi,\opi}^*\oti 1_\pi)
(1\oti 1_\pi\oti T_{\opi,\pi})
u_{\pi,\btr}
\\
=&\,
d_\pi^{-1}\ep_\pi x.
\end{align*}
Hence $\Ph_\pi$ is faithful. 

(2). 
It follows from 
\begin{align*}
\Phi_\pi(\al_\pi(a)b)
=&\,
(1\oti T_{\opi,\pi}^*)
u_{\opi,\pi}^*
(\al_{\opi}\oti\id)(\al_\pi(a)b)
u_{\opi,\pi}
(1\oti T_{\opi,\pi})
\\
=&\,
(1\oti T_{\opi,\pi}^*)
(\id\oti_{\opi}\De_{\pi})(\al(a))
u_{\opi,\pi}^*
(\al_{\opi}\oti\id)(b)
u_{\opi,\pi}
(1\oti T_{\opi,\pi})
\\
=&\,
a
(1\oti T_{\opi,\pi}^*)
u_{\opi,\pi}^*
(\al_{\opi}\oti\id)(b)
u_{\opi,\pi}
(1\oti T_{\opi,\pi})
\\
=&\,
a\Ph_\pi(b).
\end{align*}

(3).
For $x=1$, we have 
\begin{align*}
&(\Ph_\pi\oti\id)(u_{\pi,\opi}
(1\oti
T_{\pi,\opi}T_{\pi,\opi}^*))
u_{\pi,\opi}^*)
\\
=&\,
(\Ph_\pi\oti\id)(u_{\pi,\opi}
(1\oti{}_{\pi}\De_{\opi}(e_\btr))
u_{\pi,\opi}^*)
\\
=&\,
(1\oti T_{\opi,\pi}^*\oti1_{\opi})
(u_{\opi,\pi}^*\oti1_{\opi})
(\al_{\opi}\oti\id\oti\id)(u_{\pi,\opi})
\\
&\quad\cdot
(1\oti1_{\opi}\oti{}_\pi\De_{\opi}(e_\btr))
(\al_{\opi}\oti\id\oti\id)(u_{\pi,\opi}^*)
(u_{\opi,\pi}\oti1_{\opi})
(1\oti T_{\opi,\pi}\oti1_{\opi})
\\
=&\,
(1\oti T_{\opi,\pi}^*\oti1_{\opi})
(\id\oti{}_\opi\De_\pi\oti\id)(u)
(\id\oti\id\oti{}_\pi\De_\opi)(u^*)
\\
&\quad\cdot
(1\oti1_{\opi}\oti{}_{\pi}\De_{\opi}(e_\btr))
(\id\oti\id\oti{}_\pi\De_\opi)(u)(\id\oti{}_\opi\De_\pi\oti\id)(u^*)
(1\oti T_{\opi,\pi}\oti1_{\opi})
\\
=&\,
(1\oti T_{\opi,\pi}^*\oti1_{\opi})
(1\oti1_{\opi}\oti{}_{\pi}\De_{\opi}(e_\btr))
(1\oti T_{\opi,\pi}\oti1_{\opi})
\\
=&\,
d_\pi^{-2}. 
\end{align*}
Then the desired equality follows from 
\begin{align*}
(\Ph_\pi\oti\id)(u_{\pi,\opi}
(x\oti
& 
T_{\pi,\opi}T_{\pi,\opi}^*))
u_{\pi,\opi}^*)
\\
=&\,
(\Ph_\pi\oti\id)(u_{\pi,\opi}
(1\oti{}_{\pi}\De_{\opi}(e_\btr))(\id\oti{}_{\pi}\De_{\opi})(\al(x))
u_{\pi,\opi}^*)
\\
=&\,
(\Ph_\pi\oti\id)(u_{\pi,\opi}
(1\oti{}_{\pi}\De_{\opi}(e_\btr))
u_{\pi,\opi}^*(\al_{\pi}\oti\id)(\al_{\opi}(x)))
\\
=&\,
(\Ph_\pi\oti\id)(u_{\pi,\opi}
(1\oti{}_{\pi}\De_{\opi}(e_\btr))
u_{\pi,\opi}^*)
\al_{\opi}(x)
\\
=&\,
d_\pi^{-2}\al_{\opi}(x). 
\end{align*}
\end{proof}

On a composition of left inverses, the next lemma holds. 
\begin{lem}\label{lem: composition}
Let $\{\Ph_\pi\}_{\pi\in\IG}$ be the left inverses of a cocycle action 
$(\al,u)$ as before. 
Then for all $x\in M\oti B(H_\pi)\oti B(H_\rho)$, one has
\[
\Ph_\rho\big{(}
(\Ph_\pi\oti\id)(u_{\pi,\rho}\,xu_{\pi,\rho}^*)\big{)}
=
\sum_{\si\prec\pi\cdot\rho \atop \si \in \IG}
\sum_{T\in \ONB(\si,\pi\cdot\rho)}
\frac{d_\si}{d_\pi d_\rho}
\Ph_\si((1\oti T^*)x(1\oti T)).
\]
\end{lem}
\begin{rem}
This summation does not depend on the choice of an orthonormal basis 
of $(\si,\pi\cdot\rho)$. 
\end{rem}
\begin{proof}
We use the leg notations indexed by irreducible representations 
to represent positions in tensor products. 
Set $T_{\pi\cdot\rho}=T_{\opi,\pi}T_{\orho,\rho}=T_{\orho,\rho}T_{\opi,\pi}
\in (\btr,\orho\cdot\opi\cdot\pi\cdot\rho)$. 
We verify the desired equality as follows. 
\begin{align*}
&\Ph_\rho\big{(}
(\Ph_\pi\oti\id)(uxu^*)\big{)}
\\
=&\,
(1\oti T_{\pi\cdot\rho}^*)
u_{\orho,\rho}^*\al_\orho(u_{\opi,\pi}^*)
\al_\orho(\al_\opi(uxu^*))
\al_\orho(u_{\opi,\pi})u_{\orho,\rho}
(1\oti T_{\pi\cdot\rho})
\\
=&\,
(1\oti T_{\pi\cdot\rho}^*)
u_{\orho,\rho}^*
\al_\orho(u_{\opi,\pi}^*)
u_{\orho,\opi}
((\id\oti{}_\orho\De_{\opi})\circ\al(uxu^*))
u_{\orho,\opi}^*
\al_\orho(u_{\opi,\pi})u_{\orho,\rho}
(1\oti T_{\pi\cdot\rho})
\\
=&\,
(1\oti T_{\pi\cdot\rho}^*)
u_{\orho,\rho}^*
(\id\oti\id_\orho\oti{}_\opi\De_\pi)(u)
(\id\oti{}_\orho\De_\opi\oti\id_\pi)(u^*)
((\id\oti{}_\orho\De_{\opi})\circ\al(uxu^*))
\\
&\quad\cdot
(\id\oti{}_\orho\De_\opi\oti\id_\pi)(u)
(\id\oti\id_\orho\oti{}_\opi\De_\pi)(u^*)
u_{\orho,\rho}
(1\oti T_{\pi\cdot\rho})
\\
=&\,
(1\oti T_{\pi\cdot\rho}^*)
u_{\orho,\rho}^*
(\id\oti{}_\orho\De_\opi\oti\id_\pi)(u^*)
((\id\oti{}_\orho\De_{\opi})\circ\al(uxu^*))
\\
&\quad\cdot
(\id\oti{}_\orho\De_\opi\oti\id_\pi)(u)
u_{\orho,\rho}
(1\oti T_{\pi\cdot\rho})
\\
=&\,
(1\oti T_{\pi\cdot\rho}^*)
u_{\orho,\rho}^*
(\id\oti{}_\orho\De_\opi\oti\id_\pi\oti\id_\rho)
\big{(}
(u^*\oti1_\rho)\al(uxu^*)(u\oti1_\rho)
\big{)}
u_{\orho,\rho}
(1\oti T_{\pi\cdot\rho})
\\
=&\,
(1\oti T_{\pi\cdot\rho}^*)
u_{\orho,\rho}^*
(\id\oti{}_\orho\De_\opi\oti\id_\pi\oti\id_\rho)
\big{(}
(\id\oti\De\oti\id)(u)(\id\oti\id\oti\De)(u^*)
\al(x)
\\
&\quad\cdot
(\id\oti\id\oti\De)(u)
(\id\oti\De\oti\id)(u^*)
\big{)}
u_{\orho,\rho}
(1\oti T_{\pi\cdot\rho})
\\
=&\,
(1\oti T_{\pi\cdot\rho}^*)
u_{\orho,\rho}^*
(\id\oti\De^{(2)}\oti\id_\rho)(u)
(\id\oti{}_{\orho}\De_\opi\oti{}_\pi\De_\rho)(u^*)
(\id\oti{}_\orho\De_\opi)(\al(x))
\\
&\quad\cdot
(\id\oti{}_{\orho}\De_\opi\oti{}_\pi\De_\rho)(u)
(\id\oti\De^{(2)}\oti\id_\rho)(u^*)
u_{\orho,\rho}
(1\oti T_{\pi\cdot\rho})
\\
=&\,
(1\oti T_{\pi\cdot\rho}^*)
u_{\orho,\rho}^*
u_{\orho,\rho}
(\id\oti{}_{\orho}\De_\opi\oti{}_\pi\De_\rho)(u^*)
(\id\oti{}_\orho\De_\opi)(\al(x))
\\
&\quad\cdot
(\id\oti{}_{\orho}\De_\opi\oti{}_\pi\De_\rho)(u)
u_{\orho,\rho}^*
u_{\orho,\rho}
(1\oti T_{\pi\cdot\rho})
\\
=&\,
(1\oti T_{\pi\cdot\rho}^*)
(\id\oti{}_{\orho}\De_\opi\oti{}_\pi\De_\rho)(u^*)
(\id\oti{}_\orho\De_\opi)(\al(x))
(\id\oti{}_{\orho}\De_\opi\oti{}_\pi\De_\rho)(u)
(1\oti T_{\pi\cdot\rho})
\\
=&\,
(1\oti T_{\pi\cdot\rho}^*)
(\id\oti{}_{\orho}\De_\opi\oti\id_\pi\oti\id_\rho)
\big{(}
(\id\oti\id\oti{}_\pi\De_{\rho})(u^*)\al(x)
(\id\oti\id\oti{}_\pi\De_{\rho})(u)
\big{)}
\\
&\quad\cdot
(1\oti T_{\pi\cdot\rho})
\\
=&\,
\sum_{\si\prec \pi\cdot\rho}
\sum_{S\in\ONB(\osi,\orho\cdot\opi)}
(1\oti T_{\pi\cdot\rho}^*)
(1\oti S\oti 1_\pi\oti1_\rho)
(\id\oti\id_\osi\oti{}_\pi\De_{\rho})(u^*)
\al_\osi(x)
\\
&\hspace{5cm}\cdot
(\id\oti\id_\osi\oti{}_\pi\De_{\rho})(u)
(1\oti S^*\oti 1_\pi\oti1_\rho)
(1\oti T_{\pi\cdot\rho})
\\
=&\,
\sum_{\si\prec \pi\cdot\rho}
\sum_{S\in\ONB(\osi,\orho\cdot\opi)}
\frac{d_\si}{d_\pi d_\rho}
(1\oti T_{\osi,\si}^*)
(1\oti 1_\osi \oti \wdt{S}^*)
(\id\oti\id_\osi\oti{}_\pi\De_{\rho})(u^*)
\al_\osi(x)
\\
&\hspace{5cm}\cdot
(\id\oti\id_\osi\oti{}_\pi\De_{\rho})(u)
(1\oti 1_\osi \oti \wdt{S})
(1\oti T_{\osi,\si})
\\
=&\,
\sum_{\si\prec \pi\cdot\rho}
\sum_{S\in\ONB(\osi,\orho\cdot\opi)}
\frac{d_\si}{d_\pi d_\rho}
(1\oti T_{\osi,\si}^*)
u_{\osi,\si}^*
(1\oti 1_\osi \oti \wdt{S}^*)
\al_\osi(x)
(1\oti 1_\osi \oti \wdt{S})
\\
&\hspace{5cm}\cdot
u_{\osi,\si}
(1\oti T_{\osi,\si})
\\
=&\,
\sum_{\si\prec \pi\cdot\rho}
\sum_{S\in\ONB(\osi,\orho\cdot\opi)}
\frac{d_\si}{d_\pi d_\rho}
\Ph_\si((1\oti\wdt{S}^*)x (1\oti \wdt{S}))
\\
=&\,
\sum_{\si\prec \pi\cdot\rho}
\sum_{T\in\ONB(\si,\pi\cdot\rho)}
\frac{d_\si}{d_\pi d_\rho}
\Ph_\si((1\oti T^*)x (1\oti T)). 
\end{align*}
\end{proof}

Next we recall the notion of freeness for a cocycle action. 

\begin{defn}
Let $(\al,u)$ be a cocycle action of $\bhG$ 
on a von Neumann algebra $M$. 
Then it is said to be \textit{free} 
if for any $\pi\in\IG\setm\{\btr\}$, 
there exists no nonzero element $a\in M\oti B(H_\pi)$ 
with $a(x\oti1_\pi)=\al_\pi(x)a$ for all $x\in M$. 
\end{defn}

Note that freeness is stable under  perturbation, 
that is, a perturbed cocycle action of a 
free cocycle action is also free. 
The following lemma is essentially proved in \cite[Lemma 5.1]{Iz1}. 

\begin{lem}\label{lem: relcom}
Let $(\al,u)$ be a free cocycle action of $\bhG$ 
on a von Neumann algebra $M$. 
Then $\al_\pi(M)'\cap (M\oti B(H_\pi))=\al_\pi(Z(M))$ holds 
for any $\pi\in\IG$. 
\end{lem}

\begin{proof}
Take an element $a$ in $\al_\pi(M)'\cap (M\oti B(H_\pi))$. 
Then for any $x\in M$, the equality 
$a\al_\pi(x)=\al_\pi(x)a$ holds. 
Applying $\al_{\opi}$ to the first leg, we have 
\[(\al_{\opi}\oti\id_\pi)(a)
(\al_{\opi}\oti\id_\pi)(\al_\pi(x))
=
(\al_{\opi}\oti\id_\pi)(\al_\pi(x))
(\al_{\opi}\oti\id_\pi)(a).
\] 
By definition of a cocycle action, 
\[
u_{\opi,\pi}^*
(\al_{\opi}\oti\id_\pi)(a)u_{\opi,\pi}
(\id\oti{}_\opi\De_{\pi})(\al(x))
=
(\id\oti{}_\opi\De_{\pi})(\al(x))
u_{\opi,\pi}^*
(\al_{\opi}\oti\id_\pi)(a)u_{\opi,\pi}.
\]
Then 
\begin{align*}
&(1\oti T_{\opi,\pi}^*)
u_{\opi,\pi}^*
(\al_{\opi}\oti\id_\pi)(a)
u_{\opi,\pi}
(\id\oti{}_\opi\De_{\pi})(\al(x))
\\
=&\,
(1\oti T_{\opi,\pi}^*)
(\id\oti{}_\opi\De_{\pi})(\al(x))
u_{\opi,\pi}^*
(\al_{\opi}\oti\id_\pi)(a)u_{\opi,\pi}
\\
=&\,
x(1\oti T_{\opi,\pi}^*)
u_{\opi,\pi}^*
(\al_{\opi}\oti\id_\pi)(a)u_{\opi,\pi}.
\end{align*}
Multiplying $S\in (\si, \opi\cdot\pi)$ 
from the right, we have
\begin{align*}
(1\oti T_{\opi,\pi}^*)
u_{\opi,\pi}^*
(\al_{\opi}\oti\id_\pi)(a)u_{\opi,\pi}
&
(1\oti S)\al_\si(x)
\\
=&\,
x(1\oti T_{\opi,\pi}^*)
u_{\opi,\pi}^*
(\al_{\opi}\oti\id_\pi)(a)u_{\opi,\pi}(1\oti S).
\end{align*}
Since 
$(1\oti T_{\opi,\pi}^*)
u_{\opi,\pi}^*
(\al_{\opi}\oti\id_\pi)(a)u_{\opi,\pi}
(1\oti S)
\in M\oti B(H_\si,H_\btr)$, 
this is equal to 0 for $\si\neq\btr$ by freeness of $\al$. 
Using the equality
\[
1_{\opi}\oti1_\pi
=\sum_{\si\prec\opi\cdot\pi}\sum_{S\in \ONB(\si,\opi\cdot\pi)}SS^*,
\] 
we obtain 
\begin{align*}
(1\oti T_{\opi,\pi}^*)u_{\opi,\pi}^*
(\al_{\opi}\oti\id_\pi)(a)u_{\opi,\pi}
=&\,
(1\oti T_{\opi,\pi}^*)
u_{\opi,\pi}^*
(\al_{\opi}\oti\id_\pi)(a)u_{\opi,\pi}
(1\oti T_{\opi,\pi}T_{\opi,\pi}^*)
\\
=&\,
\Ph_\pi(a)(1\oti T_{\opi,\pi}^*).
\end{align*}
This equality implies that 
$\Ph_\pi(ab)=\Ph_\pi(a)\Ph_\pi(b)$ for all 
$a,b\in \al_\pi(M)'\cap (M\oti B(H_\pi))$. 
Note that $\Ph_\pi(a)\in Z(M)$. 
Hence on $\al_\pi(M)'\cap (M\oti B(H_\pi))$, 
the map $\Phi_\pi$ is a faithful $*$-homomorphism
to $Z(M)$. 
Since $\Ph_\pi$ maps $\al_\pi(Z(M))$ onto $Z(M)$, 
we have $\al_\pi(M)'\cap (M\oti B(H_\pi))=\al_\pi(Z(M))$ 
by the faithfulness of $\Ph_\pi$. 
\end{proof}

As an application of the previous lemma, 
we can show that a free cocycle action of a discrete Kac algebra 
preserves a center, although 
this result is unnecessary for our study. 
Let $(\al,u)$ be a free cocycle action of $\bhG$ on $M$. 
By Lemma \ref{lem: relcom}, we have 
$\al_\pi(M)'\cap (M\oti B(H_\pi))=\al_\pi(Z(M))$. 
Clearly $Z(M)\oti\C1_\pi\subs \al_\pi(M)'\cap (M\oti B(H_\pi))$ holds. 
Hence $Z(M)\oti\C1_\pi\subs \al_\pi(Z(M))$. 
In fact, they are equal as is shown below. 
For any $z\in Z(M)$, 
there exists $\th_\opi(z)\in Z(M)$ such that 
$\al_\pi(\th_\opi(z))=z\oti1_\pi$. 
Applying  $\Ph_\pi$ to the both sides, 
we have $\th_\opi(z)=\Ph_\pi(z\oti1_\pi)$. 
Thus
\begin{align*}
\th_\opi(\th_\pi(z))
=&\,
\Ph_\pi(\th_\pi(z)\oti1_\pi)
\\
=&\,
(1\oti T_{\opi,\pi}^*)u_{\opi,\pi}^* 
(\al_\opi(\th_\pi(z))\oti1_\pi)
u_{\opi,\pi}(1\oti T_{\opi,\pi})
\\
=&\,
(1\oti T_{\opi,\pi}^*)u_{\opi,\pi}^* 
(z\oti1_\opi\oti1_\pi)
u_{\opi,\pi}(1\oti T_{\opi,\pi})
\\
=&\,
z. 
\end{align*}
Hence $\th_\pi\in\Aut(Z(M))$ and 
$\al_\pi(z)=\th_\pi(z)\oti1_\pi$ for all $z\in Z(M)$ and $\pi\in\IG$. 
Moreover by definition of a cocycle action, 
for $z\in Z(M)$ and $\pi,\si\in\IG$ we have 
\begin{align*}
\th_\pi(\th_\si(z))\oti1_\pi\oti1_\si
=&\,
\al_\pi(\th_\si(z))\oti1_\si
\\
=&\,
(\al_\pi\oti\id_\si)(\al_\si(z))
\\
=&\,
u_{\pi,\si}(\id\oti{}_\pi\De_\si)(\al(z))u_{\pi,\si}^*
\\
=&\,
\sum_{\rho\prec \pi\cdot\si}
u_{\pi,\si}
(\id\oti{}_\pi\De_\si)(\al_\rho(z))
u_{\pi,\si}^*
\\
=&\,
\sum_{\rho\prec \pi\cdot\si}
u_{\pi,\si} \th_\rho(z)\oti {}_\pi\De_\si(1_\rho) u_{\pi,\si}^*. 
\end{align*}
Applying $\Ad u_{\pi,\si}$ to the above equality, we obtain
\[
\th_\pi(\th_\si(z))\oti1_\pi\oti1_\si
=
\sum_{\rho\prec \pi\cdot\si}
\th_\rho(z)\oti{}_\pi\De_\si(1_\rho). 
\]
Hence if $\rho,\rho'\prec \pi\cdot\si$, 
then $\th_\rho=\th_{\rho'}$. 
Summarizing these arguments, we have the following lemma. 
\begin{lem}
Let $(\al,u)$ be a free cocycle action of $\bhG$ on $M$. 
Then there exists a map $\th\col \IG\ra \Aut(Z(M))$ such that 
\begin{enumerate}

\item 
$\al_\pi(z)=\th_\pi(z)\oti1_\pi$ for all $z\in Z(M)$ and $\pi\in\IG$. 

\item
$\th_\pi\circ\th_\opi=\id$ for all $\pi\in\IG$. 

\item
If $\rho\prec \pi\cdot\si$, then $\th_\pi\circ\th_\si=\th_\rho$. 

\end{enumerate}
In particular, $\al$ preserves the center $Z(M)$. 
\end{lem}

We close this subsection with a criterion on the invariance 
of a trace. 

\begin{prop}\label{prop: free-left}
Let $(\al,u)$ be a free cocycle action on a 
finite von Neumann algebra $M$ with a normalized trace $\ta$. 
\begin{enumerate}

\item
The equality 
$(\ta\oti \ta_\pi)\circ\al_\pi\circ\Ph_\pi=\ta\oti\ta_\pi$ 
holds for all $\pi\in\IG$. 

\item 
If $u=1$ and $(\ta\oti \ta_\pi)\circ\al_\pi=\ta$ for all $\pi\in\IG$, 
then the action $\al$ preserves $\ta$, i.e. 
$(\ta\oti\id)(\al(x))=\ta(x)1$ holds for all $x\in M$. 

\end{enumerate}
\end{prop}

\begin{proof}
(1). 
Since $\al_\pi(M)'\cap (M\oti B(H_\pi))=\al_\pi(Z(M))$, 
$\al_\pi\circ\Ph_\pi$ is the unique conditional expectation 
from $M\oti B(H_\pi)$ to $\al_\pi(M)$. 
Hence the map preserves the tracial state $\ta\oti\ta_\pi$, 
that is, 
$(\ta\oti \ta_\pi)\circ\al_\pi\circ\Ph_\pi=\ta\oti\ta_\pi$. 

(2). 
By (1), the equality $\ta\circ \Ph_\pi=\ta\oti\ta_\pi$ holds. 
Since 
\begin{align*}
\Ph_\pi(x\oti e_{\pi_{i,j}})
=&\,
(1\oti T_{\opi,\pi}^*)(\al_\opi(x)\oti e_{\pi_{i,j}})
(1\oti T_{\opi,\pi})
\\
=&\,
(1\oti T_{\opi,\pi}^*)(1\oti e_{\opi_{\ovl{i},\ovl{i}}}\oti1_\pi)
(\al_\opi(x)\oti 1_\pi)(1\oti e_{\opi_{\ovl{j},\ovl{j}}}\oti1_\pi)
(1\oti T_{\opi,\pi})
\\
=&\,
d_\pi^{-1}\al_{\opi_{\ovl{i},\ovl{j}}}(x),
\end{align*}
we have 
\[d_\pi^{-1}\ta(\al_{\opi_{\ovl{i},\ovl{j}}}(x))
=(\ta\oti\ta_\pi)(x\oti e_{\pi_{i,j}})
=d_\pi^{-1}\de_{i,j}\ta(x).
\] 
\end{proof}

\subsection{Crossed products and dual functionals}

We collect well-known results on crossed product von Neumann algebras. 
Let $M$ be a von Neumann algebra and $\al$ an action of $\bhG$ on $M$. 
The \textit{crossed product} is the von Neumann subalgebra $M\rti_\al \bhG$ in 
$M\oti B(\lthG)$ defined by 
\[
M\rtimes_{\al} \bhG=\al(M)\vee (\C\oti \lG).
\] 
Set an element $\la_{\pi}=1\oti V_\pi\in (M\rti_\al \bhG)\oti \lhG$. 
Then in fact the crossed product $M\rti_\al \bhG$ is the weak closure 
of the linear space spanned by 
$\{\al(x)\la_{\pi_{i,j}}\mid x\in M,\ \pi\in\IG,\ i,j\in I_\pi\}$. 
Now define a unitary $\tV$ by 
\[
\tV=\Si (1\oti U)V (1\oti U) \Si, 
\]
where $\Si$ is the flip unitary on $\lthG \oti \lthG$. 
Then $\tV$ is a multiplicative unitary in $\lhG'\oti \lG$ 
and satisfies 
\[
\tV_{13}V_{12} \tV_{13}^*=V_{12} V_{32}. 
\]
Consider the action $\Ad(1\oti \tV)(\cdot\oti1)$ on $M\oti B(\lthG)$. 
It actually preserves $M\rti_\al\bhG$, and the restriction gives an 
action of $\bG$ on $M\rti_\al\bhG$. 
It is called the \textit{dual action} of $\al$ and denoted by $\hal$. 
By definition, we have 
\[
\hal(\al(x))=\al(x)\oti1,\quad 
\hal(\la_{\pi_{i,j}})
=\sum_{k\in I_\pi} \la_{\pi_{i,k}}\oti v_{\pi_{k,j}}
\]
for all $\pi\in\IG$, $i,j\in I_\pi$. 

Recall a normal functional $\vep_{\pi_{i,j}}\in \lG_*$ 
defined in \S2.2. 
Set linear maps 
$P_{\pi_{ij}}$ and $P_\pi$ on $M\rti_{\al}\bhG$ by 
$P_{\pi_{i,j}}=(\id\oti\id\oti \vep_{\pi_{i,j}})\circ\hal$ 
and 
$P_\pi=\sum_{i\in I_\pi}P_{\pi_{i,i}}$, respectively.
In particular, $P_\btr=(\id\oti h)\circ\hal$ 
is a faithful normal conditional expectation from 
$M\rti_\al \bhG$ onto the fixed point algebra 
$(M\rti_\al \bhG)^{\hal}$. 
We write $E_\hal$ for $P_\btr$. 
The following equalities are directly verified by easy calculation. 
Let $\pi, \rho\in\IG$, $i,j\in I_\pi$ and $k,\el\in I_\rho$. 
Then we have 
\begin{enumerate}[(i)]

\item
$
P_{\pi_{i,j}}(\al(a)\la_{\rho_{k,\el}}\al(b))
=
\de_{\pi,\rho}\de_{j,\el}
\al(a) \la_{\pi_{k,i}}\al(b)$,

\item 
$P_\pi(\al(a)\la_{\rho_{k,\el}}\al(b))
=
\de_{\pi,\rho}
\al(a) \la_{\pi_{k,\el}}\al(b)$, 

\item 
$P_\pi^2=P_\pi$. 
\end{enumerate}
In particular, the map $E_\hal$ is a conditional expectation 
from $M\rti_{\al}\bhG$ onto $\al(M)$. 
The above equalities yield 
\begin{enumerate}[(i)]

\item 
$E_\hal 
\big{(}\big{(}\al(b)\la_{\rho_{k,\el}}\big{)}^*
\al(a)\la_{\pi_{i,j}}
\big{)}
=
\de_{\pi,\rho}\de_{j,\el}
\al(\Ph_\pi(b^*a\oti e_{\pi_{k,i}}))
$,

\item 
$E_\hal\big{(}
\al(a)\la_{\pi_{i,j}}
\big{(}
\al(b)\la_{\rho_{k,\el}}
\big{)}^*
\big{)}
=
\de_{\pi,\rho}\de_{i,k}\de_{j,\el}d_\pi^{-1}
\al(ab^*)
$.
\end{enumerate}
The first equality implies 
$E_{\hal}(P_\pi(x)^*y)=E_{\hal}(x^*P_\pi(y))$ 
for all $x,y\in M\rti_{\al}\bhG$. 

For a normal functional $\th\in M_*$, 
put $\hth=\th\circ \al^{-1}\circ E_{\hal}$. 
It is called a \textit{dual functional} of $\th$. 
When $\th$ is a state, $\hth$ is called a dual state. 
From the equalities on $E_\hal$, we have 
\begin{enumerate}[(i)]

\item
$\hth
\big{(}
\big{(}\al(b)\la_{\rho_{k,\el}}\big{)}^*
\al(a)\la_{\pi_{i,j}}
\big{)}
=
\de_{\pi,\rho}\de_{j,\el}
\th(\Ph_\pi(b^*a\oti e_{\pi_{k,i}})) 
$,

\item
$\hth
\big{(}
\al(a)\la_{\pi_{i,j}}
\big{(}
\al(b)\la_{\rho_{k,\el}})
\big{)}^*
\big{)}
=
\de_{\pi,\rho}\de_{i,k}\de_{j,\el}d_\pi^{-1}
\th(ab^*)$. 
\end{enumerate}

In constructing a tower base in \S 5, 
we need these equalities for a dual state of a trace. 
The next proposition characterizes 
when a dual state of a trace is also a trace. 

\begin{prop}\label{prop: dualtrace}
Let $M$ be a finite von Neumann algebra with 
a faithful normal tracial state $\ta$. 
Then the dual state $\hta$ is tracial if and only if 
all the left inverses $\{\Ph_\pi\}_{\pi\in\IG}$ preserve 
the trace $\ta$. 
\end{prop}

If $\al$ is a free action on a finite factor, 
then the left inverse $\Ph_\pi$ preserves the trace for all $\pi\in\IG$ 
by Proposition \ref{prop: free-left}. 
Hence a dual state of the unique tracial state is also tracial. 

Finally in this subsection, 
we study the relative commutant $\al(A)'\cap (M\rti_\al\bhG)$ for 
a free cocycle action $(\al,u)$ and a von Neumann subalgebra 
$A\subs M$. 
In order to avoid ambiguous arguments on formal expansion 
of elements in a crossed product, 
we introduce the map $\meF$, which picks up coefficients 
of arbitrary elements of $M\rti_{\al}\bhG$. 

\begin{lem}
Set $\lG_\pi=\spa\{v_{\pi_{i,j}}\mid i,j\in I_\pi\}$ for each 
$\pi\in\IG$.  
\begin{enumerate}

\item
For each $\pi\in\IG$, 
the elements $\{\la_{\pi_{i,j}}\}_{i,j\in I_\pi}$ 
is a basis for the $M$-left module $\al(M)(\C\oti \lG_\pi)$. 

\item 
For each $\pi\in\IG$, 
the linear space $\al(M)(\C\oti \lG_\pi)$ 
is $\si$-weakly closed. 
In particular, $P_\pi (M\rti_{\al}\bhG)=\al(M)(\C\oti \lG_\pi)$. 

\item 
For any element $x\in M\rti_{\al}\bhG$, 
consider the element $\meF(x)=(P_\pi(x))_\pi$ 
in the linear space $\prod_{\pi\in\IG} \al(M)(\C\oti \lG_\pi)$. 
Then the map 
$\meF\col M\rti_{\al}\bhG \ra \prod_{\pi\in\IG} \al(M)
(\C\oti \lG_\pi)$ 
is an injective linear map. 
\end{enumerate}
\end{lem}

\begin{proof}
(1), (2). 
We know the equality 
$E_{\hal}\big{(}
\al(a)\la_{\pi_{i,j}}\la_{\si_{k,\el}}^*
\big{)}
=
d_{\pi}^{-1}\de_{\pi,\si}\de_{i,k}\de_{j,\el}\al(a)
$. 
For $x\in M\rti_{\al}\bhG$, define 
a element $Q_{\pi_{i,j}}(x)\in M$ by 
$\al(Q_{\pi_{i,j}}(x))
=d_\pi E_{\hal}(x\la_{\pi_{i,j}}^*)$. 
Hence the map 
$Q_{\pi_{i,j}}$ catches the coefficient of 
$\la_{\pi_{i,j}}$, so that 
they give a basis over $M$. 
Indeed, 
the $\si$-weak continuity of this map shows 
the $\si$-weak closedness of $\al(M)(\C\oti \lG_\pi)$. 

(3). 
Assume $P_\pi(x)=0$ for all $\pi\in\IG$. 
Since $E_{\hal}(P_\pi(x)^*y)=E_{\hal}(x^*P_\pi(y))$ 
for all $y\in M\rti_{\al}\bhG$, 
$E_{\hal}(x^*P_\pi(y))=0$ for all $\pi\in\IG$ and 
$y\in M\rti_{\al}\bhG$. 
Since the linear space generated by 
$\{P_\pi (M\rti_{\al}\bhG)\}_{\pi\in\IG}$ 
is dense in $M\rti_{\al}\bhG$ 
and $E_\hal$ is faithful, we have $x=0$. 
\end{proof}
Let $A$ be a von Neumann subalgebra of $M$. 
Its global invariance for $\al$ is not assumed here. 
Set $R=\al(A)'\cap (M\rti_\al \bhG)$. 
Since $\bG$ trivially acts on $\al(A)$, it preserves $R$. 
Hence the von Neumann algebra $R$ is the weak closure of 
the linear space generated by $\{P_\pi(R)\}_{\pi\in\IG}$. 
Assume that for some $\pi\in\IG$, $P_\pi(R)$ is not zero. 
Take a nonzero 
$\displaystyle 
a=\sum_{i,j\in I_\pi}\al(a_{i,j}^*)\la_{\pi_{i,j}}$ in $P_\pi(R)$. 
For $x\in A$, $a\al(x)=\al(x)a$ holds. 
Using $\displaystyle
\la_{\pi_{i,j}}(\al(x))
=\sum_{k\in I_\pi}\al(\al_{\pi_{i,k}}(x))\la_{\pi_{k,j}}
$, 
we have 
\[
\sum_{i\in I_\pi} a_{i,j}^* \al_{\pi_{i,k}}(x)
=
x a_{j,k}^*. 
\]
Then 
$\displaystyle a=\sum_{i,j\in I_\pi}a_{i,j}\oti e_{\pi_{i,j}}$ 
satisfies
$
a(x\oti1_\pi)=\al_\pi(x)a 
$ for all $x\in A$. 
Summarizing these arguments, we have the following lemma. 
\begin{lem}\label{lem: free-relcom}
Let $A$ be a von Neumann subalgebra of $M$. 
Then the relative commutant $\al(A)'\cap (M\rti_\al \bhG)$ 
is not equal to $\al(A'\cap M)$ 
if and only if there exists $\pi\in\IG\setm\{\btr\}$ and a nonzero 
$a\in M\oti B(H_\pi)$ such that 
\[
a(x\oti1_\pi)=\al_\pi(x)a \quad\mbox{for all}\ x\in A. 
\]
\end{lem}
Applying it to the case that $\al$ is free and $A=M$, 
we can derive the following well-known results. 

\begin{thm}
Let $\al$ be an action of $\bhG$ on a 
von Neumann algebra $M$. 
Then it is free if and only if 
the relative commutant $\al(M)'\cap (M\rti_{\al}\bhG)$ 
is equal to $\al(Z(M))$. 
\end{thm}

\begin{cor}\label{cor: free}
Let $\al$ be a free action of $\bhG$ on a 
factor $M$ of type \II\ with the tracial state $\ta$. 
Then the relative commutant $\al(M)'\cap (M\rti_\al\bhG)$ is trivial. 
In fact, $M\rti_{\al}\bhG$ is 
a factor of type \II\ whose tracial state is given by $\hta$. 
\end{cor}

\section{Ultraproduct von Neumann algebras}

After Connes's classification of cyclic group actions on 
the AFD factor of type \II, 
we have enjoyed benefits of ultraproduct techniques in 
studying actions of various groups. 
In those circumstances, 
it is a key that 
an automorphism on a von Neumann algebra 
induces an automorphism on its central sequence algebra. 
Although general Kac algebra actions do not have 
such a property, 
the ultraproduct technique 
plays an essential and important role for approximately inner actions. 
We begin with the notion of a convergence of homomorphisms. 

\subsection{Liftable and semiliftable homomorphisms}

Fix a free ultrafilter $\om$ on $\N$. 
For a separable von Neumann algebra $M$, 
we recall the definition of 
the ultraproduct von Neumann algebras $M^\om$ and $M_\om$. 
Let $\meT_\om$ be the set of bounded sequences 
converging to 0 strongly* in the ultralimit. 
We denote by $\meN(\meT_\om)$ 
the $C^*$-subalgebra of $\el^\infty(\N,M)$ normalizing $\meT_\om$. 
An element $(x_n)_n$ in $\el^\infty(\N,M)$ is called
$\om$-\textit{centralizing} if 
\[
\lim_{n\to\om}\|[\ph, x_n]\|=0\quad 
\mbox{for all}\ \ph\in M_*. 
\]
The $C^*$-algebra of $\om$-centralizing sequences is denoted by $\meC_\om$ 
which is a $C^*$-subalgebra of $\meN(\meT_\om)$. 
Set the quotient $C^*$-algebras 
$M^\om=\meN(\meT_\om)/\meT_\om$ and $M_\om=\meC_\om/\meT_\om$. 
The quotient map is denoted by $q$. 
Then they also have the preduals and hence are von Neumann algebras. 
We say that $(x_n)_n\in \el^\infty(\N, M)$ 
is a representing sequence of $x\in M^\om$ 
if $x=q((x_n)_n)$. 
For $u\in U(M^\om)$, 
we always take a representing sequence $(u_n)_n$ of $u$ 
such that $u_n\in U(M)$ for all $n$. 
Define a map 
$\ta^\om\col M^\om \ra M$ by 
$\displaystyle\ta^\om(x)=\lim_{n\to\om}x_n$ 
for a representing sequence $(x_n)_n$ of $x\in M^\om$. 
The convergence is taken in the $\si$-weak topology of $M$. 
Then it is a faithful normal conditional expectation 
from $M^\om$ onto $M$. 
Note that $\ta^\om$ is tracial on $M_\om$ and 
$\ta^\om(M_\om)=Z(M)$. 
The restriction of $\ta^\om$ on $M_\om$ is denoted by 
$\ta_\om$. 

Now we define the notion of convergence of homomorphisms 
and their left inverses. 

\begin{defn}
Let $M$ be a von Neumann algebra and $K$ a finite dimensional 
Hilbert space. 
Let $\al_n, \be \in \Mor_0(M, M\oti B(K))$, $n\in\N$, 
with
left inverses $\Ph_n$ and $\Ph$, respectively. 
We say that the sequence of the pairs $\{(\al_n,\Ph_n)\}_{n\in\N}$ 
\textit{converges to} $(\be,\Ph)$ 
if 
\[
\lim_{n\to\infty}\|\ph\circ \Ph_n-\ph\circ\Ph\|=0 
\quad \mbox{for all}\ \ph\in M_*. 
\]
\end{defn}

For a finite dimensional Hilbert space $K$, 
we always identify $(M\oti B(K))^\om$ with $M^\om\oti B(K)$ in a natural way. 
Let
$\al_n,\be\in \Mor_0(M,M\oti B(K))$, $n\in\N$, with left inverses 
$\Ph_n$, $\Ph^\be$, $n\in \N$, respectively. 
Assume that $(\al_n,\Phi_n)$ converges to $(\be,\Phi^\be)$. 
Define the maps 
$\ovl{\al}\col \el^\infty(\N,M)\ra \el^\infty(\N,M)\oti B(K)$ 
and 
$\ovl{\Ph}\col \el^\infty(\N,M)\oti B(K)\ra \el^\infty(\N,M)$ 
by
\[
\ovl{\al}((x_n)_n)=(\al_n(x_n))_n,\quad 
\ovl{\Ph}((x_n)_n)=(\Ph_n(x_n))_n. 
\]

\begin{lem}
In the above situation, the following conditions hold. 
\begin{enumerate}

\item
$\ovl{\al}(\meT_\om)\subs \meT_\om\oti B(K)$, 
$\ovl{\Ph}(\meT_\om\oti B(K))\subs \meT_\om$. 

\item
$\ovl{\al}(\meN(\meT_\om))\subs \meN(\meT_\om)\oti B(K)$, 
$\ovl{\Ph}(\meN(\meT_\om)\oti B(K))\subs \meN(\meT_\om)$. 
\end{enumerate}
\end{lem}
\begin{proof}
(1). 
Let $(x_n)_n\in \meT_\om$ with $\sup_n\|x_n\|\leq1$ and 
$\ph\in M_*$ a faithful normal state. 
Then 
\begin{align*}
2\|\al_n(x_n)\|_{\ph\circ\Ph^\be}^{\sharp2}
=&\,
\ph\circ\Ph^\be(\al_n(x_n^* x_n+x_n x_n^*))
\\
=&\,
(\ph\circ\Ph^\be-\ph\circ\Ph_n)(\al_n(x_n^* x_n+x_n x_n^*))
+
(\ph\circ\Ph_n)(\al_n(x_n^* x_n+x_n x_n^*))
\\
\leq&\,
\|\ph\circ\Ph^\be-\ph\circ\Ph_n\| \cdot 2 \|x_n\|^2
+
\ph(x_n^* x_n+x_n x_n^*)
\\
\leq&\,
2\|\ph\circ\Ph^\be-\ph\circ\Ph_n\|
+
2\|x_n\|_\ph^{\sharp2}
\\
\to&\,0
\end{align*}
as $n\to\om$. 
It shows that $\ovl{\al}$ preserves $\meT_\om$. 

We next show that $\ovl{\Ph}$ preserves $\meT_\om$. 
Let $(x_n)_n\in \meT_\om\oti B(K)$ with $\sup_n\|x_n\|\leq1$ 
and $\ph$ be a faithful normal state 
on $M$. 
Then 
\begin{align*}
2\|\Ph_n(x_n)\|_\ph^{\sharp2}
=&\,
\ph(\Ph_n(x_n^*)\Ph_n(x_n)+\Ph_n(x_n)\Ph_n(x_n^*))
\\
\leq&\,
\ph(\Ph_n(x_n^* x_n+x_n x_n^*))
\\
=&\,
(\ph\circ\Ph_n-\ph\circ\Ph^\be)(x_n^* x_n+x_n x_n^*)
+
\ph\circ\Ph^\be(x_n^* x_n+x_n x_n^*)
\\
\leq&\,
2\|\ph\circ\Ph_n-\ph\circ\Ph^\be\|
+
2\|x_n\|_{\ph\circ\Ph^\be}^{\sharp2}
\\
\to&\,0
\end{align*}
as $n\to\om$. 
Hence $\ovl{\Ph}(\meT_\om\oti B(K))\subs\meT_\om$. 

(2) 
Let $(x_n)_n\in \meN(\meT_\om)$ with $\sup_n\|x_n\|\leq1$ 
and 
$(y_n)_n\in \meT_\om\oti B(K)$ with $\sup_n\|y_n\|\leq1$. 
Then we have 
\begin{align*}
\|y_n\al_n(x_n)\|_{\ph\circ\Ph^\be}^2
=&\,
\ph\circ\Ph^\be(\al_n(x_n^*)y_n^* y_n \al_n(x_n))
\\
=&\,
(\ph\circ\Ph^\be-\ph\circ\Ph_n)
(\al_n(x_n^*)y_n^* y_n \al_n(x_n))
\\
&\quad
+
\ph\circ\Ph_n(\al_n(x_n^*)y_n^* y_n \al_n(x_n))
\\
\leq&\,
\|\ph\circ\Ph^\be-\ph\circ\Ph_n\|
+
\ph(x_n^*\Ph_n(y_n^* y_n)x_n).
\end{align*}
By (1), $(\Ph_n(y_n^* y_n))_n\in \meT_\om$. 
Hence the right hand side converges to 0 as $n\to\om$. 
Similarly we can show that $y_n^*\al(x_n^*)$ converges to 0 
strongly as $n\to\om$. 
Hence $\ovl{\al}$ preserves $\meN(\meT_\om)$. 

Next we show that 
$\ovl{\Ph}(\meN(\meT_\om)\oti B(K))\subs\meN(\meT_\om)$. 
Let $(x_n)_n\in \meN(\meT_\om)\oti B(K)$ with $\sup_n\|x_n\|\leq1$ and 
$(y_n)_n\in\meT_\om$ with $\sup_n\|y_n\|\leq1$. 
Then 
\begin{align*}
\|y_n \Ph_n(x_n)\|_\ph^2
=&\,
\ph(\Ph_n(x_n^*)y_n^* y_n \Ph_n(x_n))
\\
=&\,
\ph(\Ph_n(x_n^*\al_n(y_n^*))\Ph_n(\al_n(y_n)x_n))
\\
\leq&\,
\ph\circ\Ph_n(x_n^*\al_n(y_n^* y_n)x_n)
\\
=&\,
(\ph\circ\Ph_n-\ph\circ\Ph^\be)(x_n^*\al_n(y_n^* y_n)x_n)
+
\ph\circ\Ph^\be(x_n^*\al_n(y_n^* y_n)x_n)
\\
\leq&\,
\|\ph\circ\Ph_n-\ph\circ\Ph^\be\|
+
\|\al_n(y_n)x_n\|_{\ph\circ\Ph^\be}^2. 
\end{align*}
Since $(\al_n(y_n))_n\in \meT_\om\oti B(K)$, 
the right hand side converges to 0 as $n\to\om$. 
Similarly we can show that 
$(y_n^*\Ph_n(x_n^*))_n$ strongly converges to 0 as $n\to\om$. 
Hence $\ovl{\Ph}(\meN(\meT_\om)\oti B(K))\subs\meN(\meT_\om)$. 
\end{proof}

By the previous lemma, we can induce the maps 
$\al\col M^\om\ra M^\om\oti B(K)$ and $\Ph\col M^\om\oti B(K)\ra M^\om$ 
defined by 
\[
\al(x)=(q\oti\id)((\al_n(x_n))_n),
\quad
\Ph(x\oti y)=q((\Ph_n(x_n\oti y))_n)
\] 
for all $x=q((x_n)_n) \in M^\om$ and $y\in B(K)$. 
Note that $\Ph$ preserves $M_\om$, that is, $\Ph(M_\om\oti\C)\subs M_\om$. 
Indeed, for any $(x_n)_n\in\meC_\om$, $\ph\in M_*$ and $y\in M$
we have 
$
[\ph,\Ph_n(x_n\oti1)](y)
=
[\ph\circ\Ph_n, x_n\oti1](\al_n(y))
$ and hence
\begin{align*}
\|[\ph,\Ph_n(x_n\oti1)]\|
\leq&\,
\|[\ph\circ\Ph_n, x_n\oti1]\|
\\
\leq&\,
\|[\ph\circ\Ph_n-\ph\circ\Ph^\be, x_n\oti1]\|
+
\|[\ph\circ\Ph^\be, x_n\oti1]\|
\\
\leq&\,
2\|\ph\circ\Ph_n-\ph\circ\Ph^\be\|\|x_n\|
+
\|[\ph\circ\Ph^\be, x_n\oti1]\|
\\
\to&\,0
\end{align*}
as $n\to\om$. 
We verify $\al\in \Mor(M^\om, M^\om\oti B(K))$ and $\Ph$ is a left inverse 
of $\al$. 
The nontrivial points are faithfulness and normality of them. 
Since $\Ph\circ \al=\id$, $\al$ is faithful. 
For $\Ph$, 
we claim that the following equality holds. 
\[
\ta^\om\circ\Ph=\Ph^\be\circ(\ta^\om\oti\id).
\]
Once we prove this equality, the faithfulness and the normality of $\Ph$ 
immediately follow. 
Moreover with the equality $\Ph\circ\al=\id$, it also 
derives the normality of $\al$. 
Now we prove the claim as follows. 
Let $\ph\in M_*$ and $x\in M^\om\oti B(K)$ with a representing sequence 
$(x_n)_n$. 
Then we have
\begin{align*}
\ph(\ta^\om\circ\Ph(x))
=&\,
\lim_{n\to\om}
\ph(\Ph_n(x_n))
\\
=&\,
\lim_{n\to\om}
(\ph\circ\Ph_n-\ph\circ\Ph^\be)(x_n)
+
\ph\circ\Ph^\be(x_n)
\\
=&\,
\ph(\Ph^\be((\ta^\om\oti\id)(x))),
\end{align*}
where we have used 
$\displaystyle\lim_{n\to\infty}\|\ph\circ\Ph_n-\ph\circ\Ph^\be\|=0$ and 
the normality of $\Ph^\be$. 
We summarize these arguments in the following lemma. 
\begin{lem}
Let $M$ be a von Neumann algebra. 
Consider $\al_n, \be\in \Mor_0(M,M\oti B(K))$, $n\in\N$, with 
left inverses $\Ph_n$, $\Ph^\be$, respectively. 
Assume that $(\al_n,\Ph_n)$ converges to $(\be,\Ph^\be)$ and 
define the maps $\al$ and $\Ph$ as before. 
Then $\al\in \Mor_0(M^\om,M^\om\oti B(K))$ and $\Ph$ is a left inverse 
of $\al$. 
Moreover we have 
\[
\ta^\om\circ\Ph=\Ph^\be\circ(\ta^\om\oti\id). 
\]
\end{lem}

\begin{defn}\label{defn: semiliftable}
Let $\al\in\Mor_0(M^\om, M^\om\oti B(K))$ with a left inverse $\Ph$. 
\begin{enumerate}
\item
We say that the pair $(\al,\Ph)$ is \textit{semiliftable} if 
there exists $(\al_n, \Ph_n)$ and $(\be,\Ph^\be)$ which induce 
$(\al,\Ph)$ as in the previous lemma. 

\item
We say that the pair $(\al,\Ph)$ is \textit{liftable} 
if the pair is semiliftable and 
we can take $\al_n=\be$ and $\Ph_n=\Ph^\be$ for all $n\in\N$. 
In this case we write $\be^\om$ for $\al$. 

\item 
For a cocycle action $(\be,w)$ of $\bhG$ on $M^\om$, 
we say that it is \textit{semiliftable} or \textit{liftable} 
if the pairs $(\be_\pi,\Ph_\pi^\be)$, $\pi\in\IG$, are 
semiliftable or liftable, respectively. 
\end{enumerate}
\end{defn}

\begin{defn}\label{defn: approx-inner}
Let $(\al,\Ph)$ be a pair of $\al\in \Mor_0(M, M\oti B(K))$ 
and a left inverse $\Ph$ of $\al$. 
\begin{enumerate}

\item
We say that the pair is \textit{approximately inner} 
if there exists a sequence of unitaris $\{u_n\}_n\subset M\oti B(K)$ such that 
the pair $(\Ad u_n(\cdot\oti1), (\id\oti\ta_K)\circ\Ad u_n^*)$ converges to 
$(\al,\Ph)$. 

\item
For a cocycle action $(\al,u)$ of $\bhG$ on $M$, 
we say that it is an \textit{approximately inner cocycle action} 
if for each $\pi\in\IG$, 
the pair $(\al_\pi,\Ph_\pi)$ is approximately inner. 
\end{enumerate}
\end{defn}

In the above situation, 
the pairs $\{(\Ad u_n, (\id\oti\ta_K)\circ\Ad u_n^*)\}_n$ 
induce the semiliftable pair 
$(\Ad U, (\id\oti\ta_K)\circ\Ad U)$. 
By the previous lemma, we have 
\[
(\ta^\om\oti\ta_K)\circ\Ad U^*=\Ph\circ (\ta^\om\oti\id).
\]
If we apply this equality to $\al^\om$, we have 
the following equality by using 
$(\ta^\om\oti\id)\circ\al^\om=\al\circ\ta^\om$,
\[
(\ta^\om\oti\ta_K)\circ\Ad U^*\circ\al^\om=\ta^\om.
\]

Next we prepare a useful lemma which characterizes 
elements in $\meC_\om\oti B(K)$. 

\begin{lem} Let $M$ be a von Neumann algebra, 
$K$ a finite dimensional Hilbert space and 
$\tau_K$ a normalized tracial state on $B(K)$. 
Then an element $(a_n)_n\in \el^\infty(\N,M)\otimes B(K)$ 
is in $\meC_\omega \otimes B(K)$ 
if and only if 
$\displaystyle\lim_{n\to\om}\|[\psi\otimes \tau_K,a_n]\|=0$, $\psi\in M_*$. 
\end{lem}
\begin{proof}
The following equality is easily verified. 
\[
[\psi\otimes \tau_K,a_n]=\sum_{i,j}[\ps,(a_n)_{i,j}]\oti (\ta_K)e_{i,j},
\]
where $\{e_{i,j}\}_{i,j}$ is a matrix unit of $B(K)$. 
Hence if $(a_n)_n\in \meC_\om\oti B(K)$, 
then $\displaystyle\lim_{n\to\infty}\|[\psi\otimes \tau_K,a_n]\|=0$. 
Suppose $\displaystyle\lim_{n\to\infty}\|[\psi\otimes \tau_K,a_n]\|=0$.
Since 
$[\ps,(a_n)_{i,j}](y)=[\psi\otimes \tau_K,a_n](y\oti e_{j,i})$, 
we have 
$\|[\ps,(a_n)_{i,j}]\|\leq\|[\psi\otimes \tau_K,a_n]\|$. 
Hence $(a_n)_n\in \meC_\om\oti B(K)$. 
\end{proof}

If $(\al,\Ph)$ is approximately inner, 
there exists a unitary $U\in M^\om\oti B(K)$ as before. 
Then consider the map $\ga=\Ad U^* \circ\al^\om$. 
By definition of the approximate innerness, 
$(\ga,\Ph\circ\Ad U)$ is semiliftable. 
Indeed, the pairs $(\Ad u_n^*\circ \al,\Ph\circ\Ad u_n)$ 
converges to $(\cdot\oti1, \id\oti\ta_K)$. 
Since $\al=\Ad U(\cdot\oti1)$ on $M$, $\ga$ fixes $M$ 
and hence preserves $M'\cap M^\om$. 
In fact, $\ga$ preserves $M_\om$ as is shown in the 
following lemma. 

\begin{lem}\label{lem: approx-inner}
Let $\alpha\in \mathrm{Mor}_0(M,M\otimes B(K))$ be 
a $*$-homomorphism with a left inverse $\Phi$. 
Assume that the pair $(\al,\Ph)$ is approximately inner. 
Take a sequence of unitaries 
$\{u_n\}_n\subset M\otimes B(K)$ with 
\[
\displaystyle
\lim_{n\to\infty}\|(\phi\otimes \tau_K)\circ \Ad u_n^*-\phi\circ \Phi\|=0 
\quad\mbox{for all}\ \ph\in M_*.
\] 
Set $U:=(u_n)_n\in \el^\infty(\N, M\otimes B(K))$ 
and then the following statements hold. 
\begin{enumerate}
\item
$U\in \meN(\meT_\om)\oti B(K)$. 

\item 
The $*$-homomorphisms $\Ad U^* \circ\ovl{\al}$ and $\Ad U^*(\cdot\oti1)$ 
from $\el^\infty(\N,M)$ to $\el^\infty(\N,M)\oti B(K)$ preserve 
$\meT_\om$, $\meN(\meT_\om)$ and $\meC_\om$. 
\end{enumerate}
\end{lem}
\begin{proof} 
(1). 
Let $(x_n)_n\in \meT_\om \otimes B(K)$. 
Fix a faithful normal state $\phi$ on $M$. 
We show
$\displaystyle\lim_{n\to\om}
\|u_n x_n\|_{\phi\circ\Ph}^\sharp=
\lim_{n\to\om}\|x_n u_n\|_{\phi\otimes \tau_K}^\sharp=0$. 
It is trivial that 
$\displaystyle\lim_{n\to\om}\|u_nx_n\|_{\phi\otimes \tau_K}=0
=\lim_{n\to\om}\|u_n^*x_n^*\|_{\phi\otimes \tau_K}$. 
On $\|x_n u_n\|_{\phi\otimes \tau_K}$, we have 
\begin{align*}
\|x_nu_n\|_{\phi\otimes \tau_K}^2
=&\,
(\phi\otimes \tau_{K})(u_n^*x_n^*x_n u_n) \\
=&\,
(\phi\otimes \tau_{K})\circ \Ad u_n^*(x_n^*x_n) \\
=&\,
((\phi\otimes \tau_{K})\circ \Ad u_n^*-\phi\circ\Phi)(x_n^*x_n)
+ \phi\circ\Phi(x_n^*x_n) \\
\leq&\,
\|(\phi\otimes \tau_{K})\circ \Ad u_n^*-\phi\circ\Phi\|\|x_n\|^2 +
\|x_n\|_{\phi\circ\Phi}^2 \\
\to&\, 0
\end{align*}
as $n\to\om$. 
Next we show that $\{x_n^*u_n^*\}$ converges to 0 $\sigma$-strongly 
as follows, 
\begin{align*}
\|x_n^*u_n^*\|_{\ph\circ \Ph}^2
=&\,
\ph(\Ph(u_n x_n x_n^* u_n^*))
\\
=&\,
(\ph\circ\Ph-(\ph\oti\ta_K)\circ\Ad u_n^*)(u_n x_n x_n^* u_n^*)
+(\ph\oti\ta_K)(x_n x_n^*)
\\
\leq&\,
\|\ph\circ\Ph-(\ph\oti\ta_K)\circ\Ad u_n^*\|\|x_n\|^2
+
\|x_n^*\|_{\ph\oti\ta_K}^2
\\
\to&\, 0 
\end{align*}
as $n\to\om$. 
Hence $U=(u_n)_n\in \meN(\meT_\om)\oti B(K)$. 

(2). 
Set $\ga_n^1:=\Ad u_n^*\circ\alpha$ and 
$\Phi_n(x):=\Phi(u_n xu_n^*)$. 
Note that $\Phi_n$ is a left inverse of $\ga_n^1$.
By assumption, we have 
$\displaystyle\lim_{n\to\infty}\|\psi\circ \Phi_n-\psi\otimes \tau_K\|=0$.
Let $(a_n)_n\in \meC_\om$ be a centralizing sequence. 
By the previous lemma, it suffices to show
$\displaystyle
\lim_{n\to\om}\|[\ga_n^1(a_n),\psi\otimes \tau_K]\|=0$. 
We may assume $\|a_n\|\leq 1$. 
Then for $x\in M\otimes B(K)$,
\begin{align*}
| [\ga_n^1(a_n),\psi\otimes \tau_K](x) |
=&\,| [\ga_n^1(a_n),\psi\circ\Phi_n](x) |+
| [\ga_n^1(a_n),\psi\otimes \tau_K-\psi\circ\Phi_n](x) | \\
\leq&\,
\| [a_n,\psi]\|\|\Phi_n(x) \|+
2\|\psi\otimes \tau_K-\psi\circ\Phi_n\|\|x\| \\
\leq &\,
\| [a_n,\psi]\|\|x \|+
2\|\psi\otimes \tau_K-\psi\circ\Phi_n\|\|x\|. 
\end{align*}
Hence we have 
$\|[\ga_n^1(a_n),\psi\otimes \tau_K]\|
\leq \| [a_n,\psi]\|+
2\|\psi\otimes \tau_K-\psi\circ\Phi_n\|$, and it follows that 
$\displaystyle\lim_{n\to\om}\|[\ga_n^1(a_n),\psi\otimes \tau_K]\|=0$ 
because $(a_n)_n\in\meC_\om$. 

Next set $\ga_n^2=\Ad u_n^*(\cdot\oti1)\in \Mor(M,M\oti B(K))$. 
Let $(a_n)_n\in \meC_\om$ with $\sup_n \|a_n\|\leq1$.  
Then 
\begin{align*}
\|[\ph\oti\ta_K,\ga_n^2(a_n)]\|
\leq&\,
\|[\ph\oti\ta_K-\ph\circ\Ph_n,\ga_n^2(a_n)]\|
+
\|[\ph\circ\Ph_n, \ga_n^2(a_n)]\|
\\
\leq&\,
2\|\ph\oti\ta_K-\ph\circ\Ph_n\|
+
\|[\ph\circ\Ph, a_n\oti1]\circ\Ad u_n\|
\\
=&\,
2\|\ph\oti\ta_K-\ph\circ\Ph_n\|
+
\|[\ph\circ\Ph, a_n\oti1]\|
\\
\to&\,0
\end{align*}
as $n\to\om$. 
\end{proof}

In the end of this subsection, 
we state a simple criterion of approximate innerness of 
homomorphisms of the AFD factor $\meR_0$ of type \II. 
For the sake of this, we prove the following lemma. 

\begin{lem}
Let $\al$, $\al_n \in\Mor_0(M, M\oti B(K))$, $n\in\N$ 
with left inverses $\Ph$ and $\Ph_n$, $n\in\N$, respectively. 
Fix a faithful normal state $\phi\in M_*$. 
The following conditions are equivalent:
\begin{enumerate}
\item
$\displaystyle\lim_{n\to\infty}
\|\psi\circ \Phi_n-\psi\circ \Phi\|=0$ 
for all $\psi\in M_*$. 

\item 
$\displaystyle\lim_{n\to\infty}
\|\phi\circ \Phi_n-\phi\circ \Phi\|=0$ 
and 
$\displaystyle\lim_{n\to\infty}\alpha_n(a)
=\alpha(a)$ strongly for all $a\in M$.

\end{enumerate}
\end{lem}

\begin{proof}
First note that $\|\psi a\|\leq \|\psi\|\|a\|$, 
$\|a\psi\|\leq\|\psi\|\|a\|$, 
$\|\phi a\|\leq \sqrt{\|\phi\|}\|a^*\|_\phi$ 
and $\|a\|_\phi^2\leq\|a\phi\|\|a\| $ 
for $\psi\in M_*$ and a positive $\phi\in M_*$. 

$(1)\Rightarrow (2)$. 
We will show 
$\displaystyle
\lim_{n\to\infty}\|\alpha_n(a)-\alpha(a)\|_{\phi\circ \Phi}=0$ as follows.
\begin{align*}
\|\alpha_n(a)-\alpha(a)\|_{\phi\circ \Phi}^2
\leq&\,
\|(\alpha_n(a)-\alpha(a))\left(\phi\circ\Phi\right)\|
\|\alpha_n(a)-\alpha(a)\| \\ 
\leq&\,
2\|a\| 
\|\alpha_n(a)\cdot
\left(\phi\circ\Phi-\phi\circ\Phi_n\right)\|
\\
&\quad+
2\|a\|\|\alpha_n(a)\cdot\left(\phi\circ\Phi_n\right)
-\alpha(a)\cdot\left(\phi\circ\Phi\right)\|
\\
\leq&\, 
2\|a\|^2 \|\phi\circ\Phi-\phi\circ\Phi_n\|
+2\|a\|\|(a\phi)\circ\Phi_n-(a\phi)\circ\Phi\| 
\\
\to &\, 0
\end{align*}
as $n\to\infty$. 

$(2)\Rightarrow (1)$. 
At first we verify 
$\displaystyle 
\lim_{n\to\infty}
\|(\phi a)\circ \Phi-(\phi a)\circ \Phi_n\|=0$. 
This is shown as follows.
\begin{align*}
\|(\phi a)\circ \Phi-(\phi a)\circ \Phi_n\|
=&\,
\|\left(\phi\circ \Phi\right)\cdot \alpha(a)
-\left(\phi\circ \Phi_n\right)\cdot\alpha_n(a)\| 
\\
\leq &\,
\|\left(\phi\circ \Phi\right)\cdot \alpha(a)
-\left(\phi\circ \Phi\right)\cdot \alpha_n(a)\| 
\\
&\quad+
\|\left(\phi\circ \Phi\right)\cdot\alpha_n(a)
-\left(\phi\circ \Phi_n\right)\cdot \alpha_n(a)\| 
\\
\leq &\,
\|\alpha(a^*)-\alpha_n(a^*)\|_{\phi\circ \Phi} 
+
\|\phi\circ \Phi-\phi\circ \Phi_n\|\|a\| \\
\to &\, 0
\end{align*}
as $n\to\infty$. 
Since $\{\phi a\}_{a\in M}$ is dense in $M_*$, we are done.
\end{proof}

\begin{lem}
Let $\meR_0$ be the AFD factor of type \II\ with the tracial state $\ta$. 
Let $\al\in \Mor(\meR_0,\meR_0\oti B(K))$ with a left inverse $\Ph$, 
where $K$ is a finite dimensional Hilbert space. 
Assume that $\Ph$ preserves the trace, i.e. $\ta\circ\Ph=\ta\oti\ta_K$. 
Then the pair $(\al,\Ph)$ is approximately inner. 
\end{lem}
\begin{proof}
Let $M_1\subs M_2\subs\dots$ be an ascending sequence of finite dimensional 
subfactors of $\meR_0$ whose union is dense in $\meR_0$. 
Let $\{e_{i,j}^n\}_{i,j}$ be a matrix unit for $M_n$. 
By the uniqueness of the trace on $\meR_0$, 
$(\ta\oti\ta_K)\circ\al=\ta$ holds. 
Hence the projections $\al(e_{i,i}^n)$ and $e_{i,i}^n\oti1$ are 
equivalent in $\meR_0\oti B(K)$ for all $i$. 
Take an element $i_0$ and a partial isometry 
$\ovl{u}_n$ in $\meR_0\oti B(K)$ such that 
$\ovl{u}_n\ovl{u}_n^*=\al(e_{i_0,i_0}^n)$ 
and 
$\ovl{u}_n^*\ovl{u}_n=e_{i_0,i_0}^n\oti1$. 
Set 
\[
u_n=\sum_{i}\al(e_{i,i_0}^n)\ovl{u}_n (e_{i_0,i}^n\oti1).
\]
Then it is a unitary in $\meR_0\oti B(K)$ 
and 
$\al(x)=\Ad u_n(x\oti1)$ holds for all $x\in M_n$. 
Set $\Ph_n=(\id\oti\ta_K)\circ\Ad u_n^*$. 
We show that the pair $(\Ad u_n(\cdot\oti1), \Ph_n)$ 
converges to $(\al,\Ph)$. 
It is easy to see that $\Ad u_n(x\oti1)$ strongly converges 
to $\al(x)$ for all $x\in \meR_0$. 
The condition 
$
\displaystyle
\lim_{n\to\infty}\|\ta\circ\Ph_n-\ta\circ\Ph\|=0
$ 
is trivial 
because the maps $\Ph_n$ and $\Ph$ preserve the trace $\ta$. 
Hence by the previous lemma, $(\al,\Ph)$ is approximately inner. 
\end{proof}

\subsection{Strongly free cocycle actions}
In this paper, we frequently make use of 
the next two results. 
Since they are proved in a similar way to proofs of 
\cite[Lemma 5.3, Lemma 5.5]{Oc1} with a little modification, 
we omit proofs. 

\begin{lem}[Fast Reindexation Trick]\label{lem: fast}
Let $M$ be a von Neumann algebra. 
Let $N$ and $S$ be countably generated von Neumann subalgebras of 
$M^\om$, 
and $\mB$ a countable family of liftable homomorphisms whose 
elements are of the form $\be^\om$ with 
$\be\in \Mor_0(M, M\oti B(K_\be))$, where 
$K_\be$ is a finite dimensional Hilbert space. 
Take a countably generated von Neumann subalgebra 
$\wdt{N}\subset M^\om$ 
satisfying $\be^\om(N)\subs \wdt{N}\oti B(K_\be)$ 
for all $\be^\om\in \mB$ and $N\subs \wdt{N}$. 
Then there exists a map 
$\Ps\in \Mor(\wdt{N},M^\om)$ such that 
\begin{enumerate}

\item $\Ps$ is identity on $\wdt{N}\cap M$,

\item $\Ps(\wdt{N}\cap M_\om)\subs S' \cap M_\om$,

\item $\ta^\om(a\Ps(x))=\ta^\om(a)\ta^\om(x)$ 
for all $x\in \wdt{N}, a\in S$, 

\item $\be^\om(\Ps(x))=(\Ps\oti \id)(\be^\om(x))$ for 
all $x\in N$, $\be^\om\in \mB$.

\end{enumerate}
\end{lem}

\begin{lem}[Index Selection Trick]\label{lem: select}
Let $M$ be a von Neumann algebra. 
Let $\mC$ be a separable $C^*$-subalgebra of $\el^\infty(\N,M^\om)$, 
and $\mB$ a countable family of semiliftable homomorphisms, 
which is
of the form 
$\be\in\Mor(M^\om,M^\om\oti B(K_\be))$ with a finite dimensional 
Hilbert space $K_\be$ 
and also acts term by term on 
$\el^\infty(\N, M^\om)$. 
Take a $C^*$-algebra $\wdt{\mC}$ containing $\mC$ and 
preserved by each $\be\in\mB$. 
Then there exists a $*$-homomorphism 
$\Psi\col \wdt{\mC}\ra M^\om$ such that 
for any $\wdt{x}=(x_n)_n\in \mC$
\begin{enumerate}

\item 
$\displaystyle
\ta^\om(\Psi(\wdt{x}))=\lim_{n\to\om}\ta^\om(x_n)$ weakly,

\item
$\Psi(\wdt{x})=x$ if $x_n=x$ for all $n$,

\item
$\Psi(\wdt{x})\in M_\om$ if $x_n\in M_\om$ for all $n$, 

\item 
$(\Psi\oti\id)(\wdt{y})=\be(\Psi(\wdt{x}))$ 
for $\be\in\mB$ and 
$\wdt{y}=(\be(x_n))_n\in \el^\infty(\N, M^\om\oti B(K_\be))$. 
\end{enumerate}
\end{lem}

We define the notion of strong outerness of 
homomorphisms, which plays a central role 
in making Rohlin projections. 

\begin{defn}
Let $M$ be a von Neumann algebra and 
$\ga\in \Mor(M^\om, M^\om\oti B(K))$. 
Then we say that 
$\ga$ is \textit{strongly outer on $M_\om$} 
(or simply \textit{strongly outer}) 
if for any countably generated von Neumann subalgebra 
$S\subs M^\om$, 
there exists no nonzero element $a\in M^\om \oti B(K)$ 
with 
$\ga(y)a=a(y\oti1)$ for all $y\in S'\cap M_\om$. 
\end{defn}
For 
$\be\in\Mor_0(M, M\oti B(K))$, 
$\be$ is said to be strongly outer 
if $\be^\om$ is strongly outer on $M_\om$. 
Actually the strong outerness does not depend on $\om\in\beta\N\setm \N$. 
\begin{defn}
Let $M$ be a von Neumann algebra and 
$(\al,u)$ a cocycle action of $\bhG$ on $M$. 
We say that $(\al,u)$ is \textit{strongly free} 
if $\al_\pi$ is strongly outer for any nontrivial $\pi\in \IG$.
\end{defn}

It is easy to see that strong freeness implies freeness, 
but the converse does not hold in general. 
Note that if $M$ is the AFD factor of type \II, then 
the strong freeness and freeness of a cocycle action 
are equivalent, 
which is obtained by a similar argument to 
\cite[Lemma 3.4]{Co-peri} 
(see Corollary \ref{cor: cent-strong-free} in Appendix).

\section{Cohomology vanishing I}

Let $M$ be a von Neumann algebra and
$\bhG=(\lhG,\De,\vph)$ an amenable discrete Kac algebra.
Consider a cocycle action $(\al,u)$ of $\bhG$ on $M$. 
On $M\oti B(\lthG)$, 
we set the map $\ga=\si_{23}\circ(\al\oti\id)$ and 
the unitary $u_{134}$. 
Then $(\ga,u_{134})$ is a cocycle action on $M\oti B(\lthG)$. 
It is well-known that $(\ga,u_{134})$ is stabilized to an action. 
However if $\IG$ is infinite, 
the stabilization is not appropriate for our work 
on finite von Neumann algebras. 
The amenability gives us a prescription for the problem. 
Indeed, by the amenability of $\bhG$, 
we can take a sufficiently large 
finitely supported projection in $\lhG$. 
We cut $B(\lthG)$ by the projection 
and stabilize a 2-cocycle approximately. 

\subsection{2-cohomology vanishing in ultraproduct von Neumann algebras}
By making use of the relation $\vph=\Tr$ on $\lhG$, 
for an $(F,\de)$-invariant projection $K$ 
we can conclude the approximate commutativity of 
$F\oti K$ and the multiplicative unitary $W$. 

\begin{lem}\label{lem: FK}
Let $F, K\in \Projf(Z(\lhG))$. 
If $K$ is $(F,\de)$-invariant, then 
\begin{enumerate}

\item
$\|W(F\oti K)-(F\oti K)W\|_{\vph\oti\Tr}<\de^{1/2} \|F\|_\vph \|K\|_\vph$, 

\item
$\|(F\oti K)W(F\oti K)-(F\oti K)W\|_{\vph\oti\Tr}
<\de^{1/2} \|F\|_\vph \|K\|_\vph$,

\item
$\|(F\oti K)W(F\oti K)-W(F\oti K)\|_{\vph\oti\Tr}
<\de^{1/2} \|F\|_\vph\|K\|_\vph$,

\item
$\|(F\oti K)W^*(F\oti K)\cdot(F\oti K)W(F\oti K)
-F\oti K\|_{\vph\oti\Tr}
<\de^{1/2} \|F\|_\vph \|K\|_\vph$.
\end{enumerate}
\end{lem}

\begin{proof}
(1). Use Powers-St\o rmer inequality (\cite{PS}) to the left hand side of 
$|W(F\oti K)W^*-(F\oti K)|_{\vph\oti\Tr}<\de |F|_\vph |K|_\vph$. 
This inequality immediately implies the conditions (2) and (3). 

(4). This is shown as follows. 
\begin{align*}
&
\|(F\oti K)W^*(F\oti K)\cdot(F\oti K)W(F\oti K)-F\oti K\|_{\vph\oti\Tr}
\\
=&\,
\|
(F\oti K)\De(K)(F\oti K)-F\oti K
\|_{\vph\oti\vph}
\\
\leq&\,
\|F\oti K\| \|(F\oti1)\De(K)-F\oti K\|_{\vph\oti\vph}
\\
<&\,
\de^{1/2}\|F\|_\vph\|K\|_\vph. 
\end{align*}
\end{proof}
The next lemma shows that we can perturb a given 2-cocycle to 
a smaller 2-cocycle. 
The outline is as follows. 
Let $F$, $K$ be projections as in the previous lemma. 
By \cite[Lemma 3.2.1]{Jo1}, 
we can take a unitary $v\in \lhG F\oti KB(\lthG)K$ 
with
\[
\|v-(F\oti K)W(F\oti K)\|_{\vph\oti\Tr}<4\de^{1/2}\|F\|_\vph \|K\|_\vph. 
\]
If $F\geq e_\btr$, $v$ is taken as $v(e_\btr\oti K)=e_\btr\oti K$. 
The unitary $v$ plays a role of 
finite dimensional cut of the left regular representation $W$. 
The error coming from the cut is controlled by the trace norm, 
that is, strong operator topology. 
By making use of the unitary $v$, 
we can approximately stabilize a 2-cocycle 
for a cocycle action on a von Neumann algebra of type \II. 

\begin{lem}\label{lem: perturb}
Let $M$ be a von Neumann algebra of type \II\ with a faithful tracial state 
$\ta$. 
Let $(\al,u)$ be a cocycle action of $\bhG$ on $M$. 
Then for any $\vep>0$ and $F\in \Projf(Z(\lhG))$, 
there exists a unitary $w\in M\oti \lhG$ satisfying 
$w_{\btr}=1\oti e_\btr$ and 
\[
\left\|\big{(}(w\oti1)\al(w)u(\id\oti\De)(w^*)-1\big{)}
(1\oti F\oti F)
\right\|_{\ta\oti\vph\oti\vph}
<
\vep. 
\]
\end{lem}

\begin{proof}
\textbf{Step A}. 
We construct a unitary close to a cut of $W$ which acts on 
sufficiently large finite dimensional subspace of $\lthG$. 

Let $\mF$ be the support of $F$ and 
$\wdt{F}$ the central projection whose support is $\mF\cdot \mF$. 
We may assume $\btr\in \mF$, and $F\leq\wdt{F}$ holds. 
Let $\de>0$ with 
$21\de^{1/2}\|F\|_\vph\|\wdt{F}\|_\vph<\vep$ 
and $K$ a finitely supported central 
$(\wdt{F},\de)$-invariant projection. 
Set $\meH=K\lthG$.
Let $\ta_\meH$ be the normalized trace on $B(\meH)$. 
From now in this proof, 
we use $\ta_\meH$ on $B(\meH)$ for measurements of norms. 
Then by \cite[Lemma 3.2.1]{Jo1}, 
we can take a unitary $v$ in $\lhG \wdt{F}\oti B(\meH)$ satisfying 
\[\
\big{\|}
v-(\wdt{F}\oti K)W(\wdt{F}\oti K)\big{\|}_{\vph\oti\ta_\meH}
<4\de^{1/2}\big{\|}\wdt{F}\big{\|}_\vph. 
\]
By (2) and (3) in Lemma \ref{lem: FK}, we also have the inequalities
\begin{align*}
&
\|v-(\wdt{F}\oti K)W\|_{\vph\oti\ta_\meH}<5\de^{1/2}\|\wdt{F}\|_\vph, 
\\
&
\|v-W(\wdt{F}\oti K)\|_{\vph\oti\ta_\meH}<5\de^{1/2}\|\wdt{F}\|_\vph. 
\end{align*}

\textbf{Step B}. 
We regard $B(\meH)\subs M$ and 
perturb $(\al,u)$ to $(\wdt{\al},\wdt{u})$ which fixes $B(\meH)\subs M$. 

Since $M$ is of type \II, 
we can take a unital embedding $B(\meH)$ into $M$. 
Let $\{e_{i,j}\}_{i,j}$ be a system of matrix units 
generating $B(\meH)$. 
For all $\pi\in\IG$, 
$\{\al_\pi(e_{i,j})\}_{i,j}$ is also a system of 
matrix units in $M\oti B(H_\pi)$. 
Hence projections $\al_\pi(e_{i,i})$ and $e_{j,j}\oti1_\pi$ 
are equivalent. 
Then there exists a unitary $\ovl{w}_\pi$ in $M\oti B(H_\pi)$ 
such that 
$\al_\pi=\Ad \ovl{w}_\pi^*(\cdot\oti1_\pi)$ on $B(\meH)$. 
Set $\ovl{w}=(\ovl{w}_\pi)_{\pi\in\IG}$ and we have 
$\al=\Ad \ovl{w}^*(\cdot\oti1)$ on $B(\meH)$. 
Then perturb $(\al,u)$ by the unitary $\ovl{w}$ and we obtain 
a cocycle action $(\wdt{\al},\wdt{u})$ on $M$ which fixes $B(\meH)$. 

\textbf{Step C}. 
We utilize $B(\meH)$ like $B(\lthG)$ 
in order to approximately stabilize the 2-cocycle $\wdt{u}$. 

Let $M=N\oti B(\meH)$ be the tensor product decomposition. 
Since $\wdt{\al}$ fixes $B(\meH)$, 
$\wdt{u}\in N\oti\C\oti\lhG\oti\lhG$. 
Then we have a cocycle action $(\ga,\ovl{u})$ of $\bhG$ on $N$ 
such that 
$\wdt{\al}=\si_{23}\circ (\ga\oti\id)$ and 
$\wdt{u}=\si_{234}(\ovl{u}\oti1_\meH)$. 
Set $\ovl{v}=v+\wdt{F}^\per \oti K\in \lhG\oti B(\meH)$. 
Consider a unitary $w$ in $M\oti\lhG$ defined by
$w=\si_{23}((1\oti \ovl{v})\ovl{u}^{\,*})$. 
Then 
\begin{align}
&(w_F\oti F)\wdt{\al}_F(w_F) \wdt{u} (\id\oti_F\De_F)(w^*)
\notag\\
=&\,
(1\oti1\oti F\oti F)
\si_{23}((1\oti \ovl{v})\ovl{u}^{\,*}\oti1) 
\cdot 
\si_{23} ((\ga\oti\id\oti\id)(\si_{23}((1\oti \ovl{v})\ovl{u}^{\,*})))
\notag\\
&\quad
\cdot 
\si_{234}(\ovl{u}\oti1)
\cdot
\si_{234}((\id\oti\De\oti\id)(\ovl{u}(1\oti \ovl{v}^{\,*})))
\notag\\
=&\,
(1\oti1\oti F\oti F)
\si_{23}((1\oti \ovl{v})\ovl{u}^{\,*}\oti1) 
\cdot 
\si_{23}\si_{34}((1\oti1\oti v)\ga(\ovl{u}^{\,*}))
\notag\\
&\quad
\cdot 
\si_{234}(\ovl{u}\oti1)
\cdot
\si_{234}((\id\oti\De\oti\id)(\ovl{u}(1\oti \ovl{v}^{\,*})))
\notag\\
=&\,
\si_{234}
\Big{(}
(1\oti F \oti F\oti 1)
((1\oti v)\ovl{u}^{\,*}\oti1)_{1243}
\cdot (1\oti1\oti v)
\notag\\
&\quad
\cdot\ga(\ovl{u}^{\,*})(\ovl{u}\oti1) 
\cdot (\id\oti \De \oti \id)(\ovl{u}( 1\oti \ovl{v}^{\,*}))
\Big{)}
\notag\\
=&\,
\si_{234}
\Big{(}
(1\oti F \oti F\oti 1)
((1\oti v)\ovl{u}^{\,*}\oti1)_{1243}
\cdot (1\oti1\oti v)
\notag\\
&\quad
\cdot(\id\oti\id\oti\De)(\ovl{u})
(1\oti(\De\oti\id)(\ovl{v}^{\,*}))
\Big{)}
\notag\\
=&\,
\si_{234}
\Big{(}
(1\oti F \oti F\oti 1)
v_{24}v_{34}\cdot v_{34}^*\ovl{u}_{124}^{\,*}
v_{34}
\notag\\
&\quad
\cdot(\id\oti\id\oti\De)(\ovl{u})
(1\oti(\De\oti\id)(\ovl{v}^*))
\Big{)}
. \label{pert1}
\end{align}
We estimate the size of the difference of 
the right hand side and $1\oti K\oti F\oti F$ as follows. 
Since 
\begin{align*}
&
\left\|(F\oti F\oti K)(W_{13}W_{23}
-(\De\oti\id)(\ovl{v}))(F\oti F\oti K)
\right\|_{\vph\oti\vph\oti\ta_\meH}
\\
=&\,
\left\|(F\oti F\oti K)((\De\oti\id)(W)
-(\De\oti\id)(\ovl{v}))(F\oti F\oti K)
\right\|_{\vph\oti\vph\oti\ta_\meH}
\\
=&\,
\left\|({}_F\De_F\oti \id)(\ovl{v}-(1\oti K)W(1\oti K))
\right\|_{\vph\oti\vph\oti\ta_\meH}
\\
=&\,
\big{\|}({}_F\De_F\oti \id)
(v-(\wdt{F}\oti K)W(\wdt{F}\oti K))\big{\|}_{\vph\oti\vph\oti\ta_\meH}
\\
\leq&\,
\big{\|}({}_F\De\oti \id)
(v-(\wdt{F}\oti K)W(\wdt{F}\oti K))\big{\|}_{\vph\oti\vph\oti\ta_\meH}
\\
=&\,
\|F\|_\vph 
\big{\|}v-(\wdt{F}\oti K)W(\wdt{F}\oti K)
\big{\|}_{\vph\oti\ta_\meH}
\\
<&\,
4\de^{1/2} \|F\|_\vph \|\wdt{F}\|_\vph, 
\end{align*} 
we have 
\begin{align}
&\|(F\oti F\oti K)(\ovl{v}_{13}\ovl{v}_{23}
-(\De\oti\id)(\ovl{v}))
\|_{\vph\oti\vph\oti\ta_\meH}
\notag\\
\leq&\,
\|(F\oti F\oti K)(v_{13}-W_{13})v_{23}
\|_{\vph\oti\vph\oti\ta_\meH}
\notag\\
&\quad+
\|(F\oti K\oti K)W_{13}(v_{23}
-W_{23}(1\oti F\oti K))
\|_{\vph\oti\vph\oti\ta_\meH}
\notag\\
&\quad
+
\|(F\oti F\oti K)(W_{13}W_{23}
-(\De\oti\id)(\ovl{v}))(F\oti F\oti K)
\|_{\vph\oti\vph\oti\ta_\meH}
\notag\\
<&\,
\|F\|_\vph
\|v-(F\oti K)W(F\oti K)
\|_{\vph\oti\ta_\meH}
+
\|F\|_\vph
\|v-W(F\oti K)
\|_{\vph\oti\ta_\meH}
\notag\\
&\quad
+
4\de^{1/2}\|F\|_\vph\|\wdt{F}\|_\vph
\notag\\
<&\,
4\de^{1/2}\|F\|_\vph^2+5\de^{1/2}\|F\|_\vph^2
+4\de^{1/2}\|F\|_\vph\|\wdt{F}\|_\vph
\notag\\
\leq&\,
13\de^{1/2}\|F\|_\vph\|\wdt{F}\|_\vph. 
\label{pert2}
\end{align}
Set $\ovl{u}_K=\ovl{u}(1\oti1\oti K)$ and then 
\begin{align}
&
\|(1\oti F\oti F\oti 1)
(v_{34}^* \ovl{u}_{124}^{\,*} v_{34}-(\id\oti\id\oti\De)(\ovl{u}_K^{\,*}))
\|_{\ta\oti\vph\oti\vph\oti\vph}
\notag\\
\leq&\,
\|(1\oti F\oti F\oti K)
(v_{34}^*- ((F\oti K)W^*(F\oti K))_{34} 
\ovl{u}_{124}^{\,*} v_{34}
\|_{\ta\oti\vph\oti\vph\oti\ta_\meH}
\notag\\
&\quad+
\big{\|}
(1\oti F\oti F\oti K)
W_{34}^* \ovl{u}_{124}^{\,*} 
\notag\\
&\qquad\cdot
\big{(}
((F\oti K)v(F\oti K))_{34}-((F\oti K)W(F\oti K))_{34}
\big{)}
\big{\|}_{\ta\oti\vph\oti\vph\oti\ta_\meH}
\notag\\
&\quad+
\|(1\oti F\oti F\oti K)
(W_{34}^*(\ovl{u}_K^{\,*})_{124}W_{34}-
(\id\oti\id\oti\De)(\ovl{u}_K^{\,*}))
\|_{\ta\oti\vph\oti\vph\oti\ta_\meH}
\notag\\
\leq&\,
\|F\|_\vph 
\|(F\oti K)v(F\oti K)-(F\oti K)W(F\oti K)
\|_{\vph\oti\ta_\meH}
\notag\\
&\quad+
\|F\|_\vph 
\|(F\oti K)v(F\oti K)-(F\oti K)W(F\oti K)\|_{\vph\oti\ta_\meH}
\notag
\\
<&\,
4\de^{1/2}\|F\|_\vph^2
+
4\de^{1/2}\|F\|_\vph^2
\notag\\
=&\,
8\de^{1/2}\|F\|_\vph \|\wdt{F}\|_\vph.
\label{pert3}
\end{align}
By using (\ref{pert1}), (\ref{pert2}) and (\ref{pert3}), 
we have 
\begin{align*}
&
\|(w_F\oti F)\wdt{\al}_F(w_F) \wdt{u} (\id\oti_F\De_F)(w^*)
-1\oti K \oti F\oti F
\|_{\ta\oti\vph\oti\vph}
\\
=&\,
\big{\|}
\si_{234}
\big{(}
(1\oti F \oti F\oti 1)
v_{24}v_{34}\cdot v_{34}^* \ovl{u}_{124}^{\,*}
v_{34}
\cdot(\id\oti\id\oti\De)(\ovl{u})
(1\oti(\De\oti\id)(\ovl{v}^*))
\big{)}
\\
&\qquad
-1\oti K\oti  F\oti F
\big{\|}_{\ta\oti\ta_\meH\oti\vph\oti\vph}
\\
=&\,
\|
v_{24}v_{34}\cdot v_{34}^*\ovl{u}_{124}^{\,*}
v_{34}
\cdot(\id\oti\id\oti\De)(\ovl{u})
(1\oti(\De\oti\id)(\ovl{v}^*))(1\oti F\oti F\oti K)
\\
&\qquad
-1\oti F\oti  F\oti K
\|_{\ta\oti\vph\oti\vph\oti\ta_\meH}
\\
<&\,
13\de^{1/2}\|F\|_\vph\|\wdt{F}\|_\vph
\\
&\quad
+
\|
(1\oti(\De\oti\id)(\ovl{v}))
\cdot v_{34}^*\ovl{u}_{124}^{\,*}
v_{34}
\cdot(\id\oti\id\oti\De)(\ovl{u})
(1\oti(\De\oti\id)(\ovl{v}^*))
\\
&\hspace{3cm}\cdot(1\oti F\oti F\oti K)
-1\oti F\oti  F\oti K
\|_{\ta\oti\vph\oti\vph\oti\ta_\meH}
\\
<&\,
13\de^{1/2}\|F\|_\vph\|\wdt{F}\|_\vph
+
8\de^{1/2}\|F\|_\vph\|\wdt{F}\|_\vph
\\
=&\,
21\de^{1/2}\|F\|_\vph\|\wdt{F}\|_\vph
\\
<&\,
\vep. 
\end{align*}
\end{proof}

We have shown that any 2-cocycle can be approximately stabilized 
in Lemma \ref{lem: perturb}. 
When we consider the stabilization problem in an ultraproduct 
von Neumann algebra, 
approximate stabilization yields exact stabilization 
by the Index Selection Trick. 

\begin{lem}\label{lem: vanish-ultra}
Let $M$ be a von Neumann algebra such that $M_\om$ is of type \II\ 
and 
$(\ga,w)$ a cocycle action of $\bhG$ on $M^\om$ 
preserving $M_\om$ and $w\in M_\om\oti\lhG\oti\lhG$. 
Assume that $\ga$ is of the form $\ga=\Ad U\circ \be$ 
where $U\in U(M^\om\oti \lhG)$ and $\be\in\Mor(M^\om, M^\om\oti \lhG)$ 
with a semiliftable $\be_\pi$ for all $\pi$. 
Then the 2-cocycle $w$ is a coboundary in $M_\om$. 
\end{lem}

\begin{proof}
Take an increasing sequence of projections $\{F_n\}_{n=1}^\infty$ 
in $\Projf(Z(\lhG))$ with $F_n\to1$ strongly, 
and a decreasing positive numbers $\{\vep_n\}_{n=1}^\infty$ 
with $\vep_n\to0$. 
Let $\ph$ be a faithful normal state on $M$ and 
set $\ps=\ph\circ\ta^\om$, which is a trace on $M_\om$. 
By using the previous lemma, for each $n\in\N$ 
we can find a unitary $v_n$ in $M_\om \oti \lhG$ satisfying 
\[
\big{\|}
\big{(}(v_n\oti1)\ga(v_n)w(\id\oti\De)(v_n^*)-1\big{)}
(1\oti F_n\oti F_n)
\big{\|}_{\ps\oti\vph\oti\vph}
<
\vep_n. 
\]
Then set a unitary $\wdt{v}=(v_n)_n$ and 
$\wdt{U}=(U)_n$ 
in $\el^\infty(\N,M^\om\oti \lhG)$ 
and $\wdt{w}=(w)_n$ in $\el^\infty(\N,M^\om\oti \lhG\oti\lhG)$. 
Let $\mC$ be a $C^*$-subalgebra generated by 
$\wdt{v}_{\pi_{i,j}}$, $\wdt{U}_{\pi_{i,j}}$ 
$\wdt{w}_{\pi_{i,j},\rho_{k,\el}}$ 
for all $\pi,\rho\in\IG$, 
$i,j\in I_\pi$ and $k,\el\in I_\rho$. 
Let $\mB=\{\be\}$. 
Then applying the Index Selection Trick, we have a 
$*$-homomorphism $\Ps\col \wdt{\mC}\ra M^\om$ satisfying 
the conditions in Lemma \ref{lem: select} 
for $\mC$ and $\mB$. 
Set $v=(\Ps\oti\id)(\wdt{v})\in M_\om\oti\lhG$ and 
\[
\wdt{x}
=
(\wdt{v}\oti1)\ga(\wdt{v})
\wdt{w}(\id\oti\De)(\wdt{v}^*)
-1\oti1\oti1.
\]
By definition of $\Ps$, 
\[
(\Ps\oti\id\oti\id)(\wdt{x})
=
(v\oti1)\ga(v)
w(\id\oti\De)(v^*)
-1\oti1\oti1. 
\]
The right hand side is equal to 0. 
Indeed for any $\pi,\rho\in\IG$,  
\begin{align*}
&
\big{\|}
(v_\pi\oti1_\rho)\ga_\pi(v_\rho)
w_{\pi,\rho}(\id\oti{}_\pi\De_\rho)(v^*)-1\oti1_\pi\oti1_\rho
\big{\|}_{\ps\oti\vph\oti\vph}^2
\\
=&\,
(\ps\oti\vph\oti\vph)
(|(\Ps\oti\id\oti\id)(\wdt{x}_{\pi,\rho})|^2)
\\
=&\,
(\ph\circ\ta^\om\circ\Ps\oti\vph\oti\vph)
(|\wdt{x}_{\pi,\rho}|^2)
\\
=&\,
\lim_{n\to\om}
(\ph\circ\ta^\om\oti\vph\oti\vph)(|(\wdt{x}_n)_{\pi,\rho}|^2)
\\
\leq&\,
\lim_{n\to\om}\vep_n^2
\\
=&\,
0.
\end{align*}
\end{proof}

The previous 2-cohomology vanishing result 
yields two results 
about approximately inner (cocycle) actions, 
which play crucial roles in our study. 
We separately discuss them in the following subsections. 
We prepare the equivalence relation $\sim$ for sequences 
of normal functionals. 
Let $(\ph_n)_n$ and $(\ps_n)_n$ be sequences of normal functionals 
on a von Neumann algebra. 
We write $(\ph_n)_n \sim (\ps_n)_n$ or simply $\ph_n\sim \ps_n$ 
when 
$\displaystyle\lim_{n\to\om}\|\ph_n-\ps_n\|=0$. 

\subsection{Cocycle actions on central sequence algebras}

Let $M$ be a von Neumann algebra such that $M_\om$ 
is of type \II\  and 
$(\al,u)$  an approximately inner cocycle action of 
$\bhG$ on $M$. 
Let $\Ph_\pi^\al$ be the left inverse of $(\al,u)$ and 
$\pi\in\IG$. 
Then there exists a unitary $U\in M^\om\oti \lhG$ satisfying 
the conditions in Lemma \ref{lem: approx-inner}. 
Then $\al=\Ad U(\cdot\oti1)$ on $M$. 
Set $\ga=\Ad U^*\circ\al^\om$ and a unitary 
$w=(U^*\oti1)\al^\om(U^*)u(\id\oti\De)(U)$. 
Then $(\ga,w)$ is a cocycle action on $M^\om$ fixing $M$. 
Note that each $\ga_\pi$ is semiliftable. 
Indeed, let $(u_n)_n$ be a representing sequence of $U$. 
By the proof of Lemma \ref{lem: approx-inner}, 
we know that 
$(\Ad u_n^*\circ\al_\pi,\Ph_\pi^\al\circ\Ad u_n)$ converges to 
$(\cdot\oti1_\pi, \id\oti \ta_\pi)$. 
This result and the equality 
$\Ph_\pi^\ga=\Ph_\pi^{\al^\om}\circ\Ad U$ 
yield the semiliftability of $(\ga_\pi,\Ph_\pi^\ga)$. 

The map $\ga_\pi$ preserves $M_\om$ by Lemma \ref{lem: approx-inner}. 
Since $\ga$ fixes $M$, $w\in M'\cap M^\om\oti\lhG\oti\lhG$, 
but in fact $w$ is a 2-cocycle whose entries are evaluated in $M_\om$. 

\begin{lem}
The unitary $w$ is in $M_\om\oti\lhG\oti\lhG$. 
\end{lem}
\begin{proof}
Let $\pi,\rho\in\IG$ and $\ph\in M_*$. 
Let $(U^n)_n$ be a representing sequence of $U$ and 
set $w^n=(U^{n*}\oti1)\al^\om(U^{n*})u(\id\oti\De)(U^{n})$. 
Then $(w^n)_n$ is a representing sequence of $w$. 
We show that 
$\displaystyle\lim_{n\to\om}
\|[\ph\oti\ta_\pi\oti\ta_\rho, w_{\pi,\rho}^n]\|=0$. 
This is verified as follows. 
\begin{align*}
&w_{\pi,\rho}^{n*} 
(\ph\oti\ta_\pi\oti\ta_\rho) w_{\pi,\rho}^n
\\
=&\,
(\id\oti\De)(U^{n*})u^*\al^\om(U^n)(U^n\oti1)
(\ph\oti\ta_\pi\oti\ta_\rho)(U^{n*}\oti1)\al^\om(U^{n*})u(\id\oti\De)(U^{n})
\\ 
\sim&\,
(\id\oti\De)(U^{n*})u^*\al^\om(U^n)
(\ph\circ\Ph_\pi^{\al^\om} \oti\ta_\rho)\al^\om(U^{n*})u(\id\oti\De)(U^{n})
\\
=&\,
(\id\oti\De)(U^{n*})u^*
\big{(}
(U(\ph\oti\ta_\rho)U^*)\circ(\Ph_\pi^{\al^\om}\oti\id_\rho)
\big{)}
u(\id\oti\De)(U^{n})
\\
\sim&\,
(\id\oti\De)(U^{n*})u^*
\big{(}
\ph\circ\Ph_\rho^\al\circ(\Ph_\pi^{\al^\om}\oti\id_\rho)
\big{)}
u(\id\oti\De)(U^{n})
\\
=&\,
\sum_{\si\prec\pi\cdot\rho}
\sum_{T \in\ONB(\si, \pi\cdot\rho)}
\frac{d_\si}{d_\pi d_\rho}
(\id\oti\De)(U^{n*})
\big{(}
(1\oti T)
\ph\circ \Ph_\si^{\al^\om}
(1\oti T^*)
\big{)}
(\id\oti\De)(U^{n})
\\
=&\,
\sum_{\si\prec\pi\cdot\rho}
\sum_{T \in\ONB(\si, \pi\cdot\rho)}
\frac{d_\si}{d_\pi d_\rho}
(1\oti T) 
\ph\circ(U_\si^{n*} \Ph_\si^{\al^\om} U_\si^n)
(1\oti T^*)
\\
\sim&\,
\sum_{\si\prec\pi\cdot\rho}
\sum_{T \in\ONB(\si, \pi\cdot\rho)}
\frac{d_\si}{d_\pi d_\rho}
(1\oti T) 
(\ph\oti\ta_\si)
(1\oti T^*)
\\
=&\,
\ph\oti\ta_\pi\oti\ta_\rho,
\end{align*}
where we have used the composition rule of the left inverses 
in Lemma \ref{lem: composition} 
and the property of $U$ in Lemma \ref{lem: approx-inner}. 
\end{proof}

Hence the restriction of $(\ga,w)$ on $M_\om$ is a cocycle action. 
Since $\ga$ is semiliftable, 
we can apply Lemma \ref{lem: vanish-ultra}. 
Hence there exists a unitary $\ovl{v}\in M_\om\oti\lhG$ such 
that 
\[
(\ovl{v}\oti1)\ga(\ovl{v})w(\id\oti\De)(\ovl{v}^*)=1. 
\]
Since $w=(U^*\oti1)\al^\om(U^*)u(\id\oti\De)(U)$, we have 
\[
(\ovl{v}U^*\oti1)\al^\om(\ovl{v}U^*)u(\id\oti\De)(U\ovl{v}^*)=1. 
\]
Setting $V=U\ovl{v}$, we have the following lemma. 

\begin{lem}\label{lem: reduction-to-action}
Let $(\al,u)$ be an approximately inner cocycle action of $\bhG$ 
on a von Neumann algebra $M$ such that $M_\om$ is of type \II. 
Then there exists a unitary $V\in M^\om\oti\lhG$ such that 
\begin{enumerate}

\item For a representing sequence $(v_n)_n$ of $V$, we have 
\[
\lim_{n\to\om}
\|(\ph\oti\ta_\pi)\circ\Ad v_n^*-\ph\circ\Ph_\pi^\al\|=0
\quad
\mbox{for all}\  \ph\in M_*,\ \pi\in\IG,
\]
where $\Ph_\pi^\al$ is the left inverse of $(\al,u)$. 

\item 
$\ga=\Ad V^*\circ \al^\om$ is an action of $\bhG$ on $M^\om$ 
fixing $M$ and preserving $M_\om$. 

\item $u=(V\oti1)\ga(V)(\id\oti\De)(V^*)$

\item 
$\ta^\om\circ\Ph_\pi^\ga=(\ta^\om\oti\ta_\pi)$ 
for all $\pi\in\IG$, 
where $\Ph_\pi^\ga$ is the left inverse of $\ga$. 

\item $(\ta^\om\oti\ta_\pi)\circ\ga_\pi(x)=\ta^\om(x)$
for all $x\in M^\om$ and $\pi\in\IG$. 

\end{enumerate}
\end{lem}
\begin{proof}
Let $V=U\ovl{v}$ as before. 
Let $(U^n)_n$ and $(\ovl{v}^n)_n$ be 
representing sequences of $U$ and $\ovl{v}$, repectively. 
Set $(v^n)_n=(U^n \ovl{v}^n)_n$ which represents $V$. 

(1). It is verified as 
\begin{align*}
(\ph\oti\ta_\pi)\circ\Ad v^{n*}
=&\,
(\ph\oti\ta_\pi)\circ\Ad \ovl{v}^{n*} U^{n*}
\\
\sim&\,
(\ph\oti\ta_\pi)\circ\Ad U^{n*}
\\
\sim&\,
\ph\circ\Ph_\pi^\al. 
\end{align*}

The conditions (2) and (3) have been already shown. 

(4). 
The left inverse of $\ga$ is given by 
$\Ph_\pi^\ga=\Ph_\pi^{\al^\om}\circ\Ad V$. 
Then for $\ph\in M_*$ and $x\in M^\om\oti B(H_\pi)$, 
\begin{align*}
\ph(\ta^\om(\Ph_\pi^\ga(x)))
=&\,
\ph(\ta^\om(\Ph_\pi^\al(VxV^*)))
\\
=&\,
\lim_{n\to\om}
\ph(\Ph_\pi^\al(v_nx_n {v_n}^*))
\\
=&\,
\lim_{n\to\om}(\ph\oti\ta_\pi)(x_n)
\\
=&\,
\ph((\ta^\om\oti\ta_\pi)(x)).
\end{align*}

(5). It is a direct consequence of (4). 

\end{proof}

Let $(\al,u)$ be an approximately inner cocycle action on $M$. 
By the previous lemma, there exists a unitary $V\in M^\om\oti\lhG$ 
such that 
\[
(V^*\oti1)\al^\om(V^*)u(\id\oti\De)(V)=1. 
\]
Therefore we can perturb the cocycle action to an action on $M^\om$. 
By taking a representing sequence of $V$, 
we can make $u$ close to $1$ with an arbitrarily small error in $M$. 
A problem is that we have no estimates of perturbation unitaries. 
For the sake of solving that, we will use the Rohlin type theorem 
presented in Theorem \ref{thm: jRohlin}. 

\subsection{Intertwining cocycles}
In this subsection, we study two approximately inner actions. 
By the 2-cohomology vanishing result, 
we can take a 1-cocycle intertwining them 
in an ultraproduct von Neumann algebra. 

\begin{lem}\label{lem: coc-ultra}
Let $M$ be a von Neumann algebra such that $M_\om$ is of type \II. 
Let $\al$ and $\be$ be approximately inner actions of $\bhG$ on $M$. 
Then there exists an $\al^\om$-cocycle $W$ in $M^\om\oti \lhG$ with 
$\be=\Ad W\circ \al$ on $M\subs M^\om$. 
\end{lem}

We prove this result after proving Lemma \ref{lem: wmom}. 
We denote  
the left inverses of $\al$ and $\be$ by 
$\{\Ph_\pi^\al\}_{\pi\in\IG}$ and $\{\Ph_\pi^\be\}_{\pi\in\IG}$, 
respectively.
Take unitaries $U$ and $V$ in $M^\om\oti\lhG$ such that 
they satisfy the conditions in Lemma \ref{lem: approx-inner} 
for $\al$ and $\be$, respectively. 
Then we have $\al=\Ad U(\cdot\oti1)$, 
$\be=\Ad V(\cdot\oti1)$ on $M$. 
Define a map $\ga\in \Mor(M^\om,M^\om\oti\lhGG)$ by
\[
\ga(x)=U_{12}^{*}\al^\om(V^{*})(\al^\om(x)\oti1)
\al^\om(V)U_{12}
\quad \mbox{for}\ x\in M^\om.
\]
Set $\wdt{\al}(x)=\al(x)\oti1$ 
and then $\wdt{\al}$ is an action of $\bhGG$ on $M$. 
Since $\ga$ is the perturbation of the action $\wdt{\al}$ by 
the unitary $U_{12}^{*}\al^\om(V^{*})$, 
a 2-cocycle $w\in M^\om\oti\lhGG^{\oti2}$ is given by 
\[
w=U_{12}^{*}\al^\om(V^{*})
\wdt{\al}^\om\left(U_{12}^{*}\al^\om(V^{*})\right)
(\id\oti\De_{\bhGG})(\al^\om(V)U_{12}). 
\]
Then $(\ga,w)$ is a cocycle action of $\bhGG$ on $M^\om$. 
The map $\ga$ is a composition of the maps 
$\Ad U_{12}^{*}\circ\wdt{\al}^\om$ and $\Ad V^{*}(\cdot\oti1)$. 
Since they preserve $M_\om$, so does $\ga$. 
We prove $w\in M_\om\oti\lhGG^{\oti2}$ as follows. 
\begin{lem}\label{lem: wmom}
The unitary $w$ is in $M_\om\oti\lhGG^{\oti2}$. 
\end{lem}
\begin{proof}
Let $(u^n)_n$ and $(v^n)_n$ be representing sequences of 
$U$ and $V$, respectively. 
Set 
\[
w_n=u_{12}^{n*}\al(v^{n*})
\wdt{\al}\left(u_{12}^{n*}\al(v^{n*})\right)
(\id\oti\De_{\bhGG})(\al(v^n)u_{12}^n).
\]
Then $(w^n)_n$ is a representing sequence of $w$. 
Let $\ph\in M_*$ and $\pi, \rho, \si, \ze\in\IG$. 
We show 
$\displaystyle
\lim_{n\to\om}\|[\ph\oti\ta_\pi\oti\ta_\rho\oti\ta_\si\oti\ta_\ze,w^n]\|
=0$. 
In order to do, we estimate 
\[
\ps_n=w^{n*}(\ph\oti\ta_\pi\oti\ta_\rho\oti\ta_\si\oti\ta_\ze)w^n.
\]
Use 
$u_\pi^{n}(\th\oti\ta_\pi)u_\pi^{n*}
\sim \th\circ\Ph_\pi^\al
$, 
$v_\pi^{n}(\th\oti\ta_\pi)v_\pi^{n*}
\sim \th\circ\Ph_\pi^\be
$ for all $\th\in M_*$ 
and then 
\begin{align*}
\ps_n\sim&\,
(\id\oti\De_{\bhGG})(u_{12}^{n*}\al(v^{n*}))
\wdt{\al}\left(\al(v^{n})u_{12}^{n}\right)
\al(v^{n})
\\
&\quad\cdot(\ph\circ\Ph_\pi^\al\oti\ta_{\rho}\oti\ta_\si\oti\ta_\ze)
\\
&\quad\cdot
\al(v^{n*})
\wdt{\al}\left(u_{12}^{n*}\al(v^{n*})\right)
(\id\oti\De_{\bhGG})(\al(v^n)u_{12}^n)
\\
=&\,
(\id\oti\De_{\bhGG})(u_{12}^{n*}\al(v^{n*}))
\al_\pi\big{(}\al(v^{n})u_{12}^{n}\big{)}_{1245}
\al_\pi(v^{n})
\\
&\quad\cdot(\ph\circ\Ph_\pi^\al\oti\ta_{\rho}\oti\ta_\si\oti\ta_\ze)
\\
&\quad\cdot
\al_\pi(v^{n*})
\al_\pi\big{(}u_{12}^{n*}\al(v^{n*})\big{)}_{1245}
(\id\oti\De_{\bhGG})(\al(v^n)u_{12}^n)
\\
=&\,
(\id\oti\De_{\bhGG})(u_{12}^{n*}\al(v^{n*}))
\al_\pi\big{(}\al(v^{n})u_{12}^{n}\big{)}_{1245}
\\
&\quad\cdot
(v^{n}(\ph\oti\ta_{\rho})v^{n*}\oti\ta_\si\oti\ta_\ze)
\circ(\Ph_\pi^\al\oti\id_\rho\oti\id_\si\oti\id_\ze)
\\
&\quad\cdot
\al_\pi\big{(}u_{12}^{n*}\al(v^{n*})\big{)}_{1245}
(\id\oti\De_{\bhGG})(\al(v^n)u_{12}^n)
\\
\sim&\,
(\id\oti\De_{\bhGG})(u_{12}^{n*}\al(v^{n*}))
\al_\pi\big{(}\al(v^{n})u_{12}^{n}\big{)}_{1245}
\\
&\quad\cdot
(\ph\circ\Ph_\rho^\be\oti\ta_\si\oti\ta_\ze)
\circ(\Ph_\pi^\al\oti\id_\rho\oti\id_\si\oti\id_\ze)
\\
&\quad\cdot
\al_\pi\big{(}u_{12}^{n*}\al(v^{n*})\big{)}_{1245}
(\id\oti\De_{\bhGG})(\al(v^n)u_{12}^n)
\\
=&\,
(\id\oti\De_{\bhGG})(u_{12}^{n*}\al(v^{n*}))
\al_\pi\big{(}\al(v^{n})\big{)}_{1245}
\\
&\quad\cdot
\big{(}
u_{13}^{n}
(\ph\circ\Ph_\rho^\be\oti\ta_\si\oti\ta_\ze)
u_{13}^{n*}
\big{)}
\circ(\Ph_\pi^\al\oti\id_\rho\oti\id_\si\oti\id_\ze)
\\
&\quad\cdot
\al_\pi\big{(}\al(v^{n*})\big{)}_{1245}
(\id\oti\De_{\bhGG})(\al(v^n)u_{12}^n)
\\
\sim&\,
(\id\oti\De_{\bhGG})(u_{12}^{n*}\al(v^{n*}))
\al_\pi\big{(}\al_\si(v_\ze^{n})\big{)}_{1245}
\\
&\quad\cdot
(\ph\oti\ta_\ze)
\circ(\Ph_\rho^\be\oti\id_\ze)\circ\Ph_\si^\al
\circ(\Ph_\pi^\al\oti\id_\rho\oti\id_\si\oti\id_\ze)
\\
&\quad\cdot
\al_\pi\big{(}\al_\si(v_\ze^{n*})\big{)}_{1245}
(\id\oti\De_{\bhGG})(\al(v^n)u_{12}^n)
\\
=&\,
(\id\oti\De_{\bhGG})(u_{12}^{n*}\al(v^{n*}))
\\
&\quad\cdot
\big{(}
v_\ze\cdot
(\ph\oti\ta_\ze)
\circ(\Ph_\rho^\be\oti\id_\ze)
\cdot v_\ze^*
\big{)}
\circ\Ph_\si^\al
\circ(\Ph_\pi^\al\oti\id_\rho\oti\id_\si\oti\id_\ze)
\\
&\quad\cdot
(\id\oti\De_{\bhGG})(\al(v^n)u_{12}^n)
\\
\sim&\,
(\id\oti\De_{\bhGG})(u_{12}^{n*}\al(v^{n*}))
\\
&\quad\cdot
\ph\circ\Ph_\rho^\be\circ\Ph_\ze^\be
\circ\Ph_\si^\al
\circ(\Ph_\pi^\al\oti\id_\rho\oti\id_\si\oti\id_\ze)
\\
&\quad\cdot
(\id\oti\De_{\bhGG})(\al(v^n)u_{12}^n).
\end{align*}
By Lemma \ref{lem: composition}, we have 
\[
\Ph_\si^\al\circ\Ph_\pi^\al
=
\sum_{\xi\prec \pi\cdot\si}
\sum_{S\in \ONB(\xi,\pi\cdot\si)}
\frac{d_\xi}{d_\pi d_\si}
(1\oti S)\Ph_\xi^\al(1\oti S^*)
\]
on $M\oti B(H_\pi)\oti B(H_\si)$ and 
\[
\Ph_\rho^\be\circ\Ph_\ze^\be
=
\sum_{\eta\prec \ze\cdot\rho}
\sum_{T\in \ONB(\eta,\ze\cdot\rho)}
\frac{d_\eta}{d_\ze d_\rho}
(1\oti \Si_{\rho,\ze}T)\Ph_\eta^\be(1\oti (\Si_{\rho,\ze}T)^*)
\]
on $M\oti B(H_\rho)\oti B(H_\ze)$, 
where $\Si_{\rho,\ze}$ is the flip unitary of $H_\rho$ and $H_\ze$. 
Note that $\rho$ sits right from $\ze$ in that lemma, 
but it does left here. 
The flip arises for this reason. 
Also note that $\Si_{\rho,\ze}T$ may not be an intertwiner between 
$\eta$ and $\rho\cdot\ze$. 
Using these, we have 
\begin{align*}
&\Ph_\rho^\be\circ\Ph_\ze^\be
\circ\Ph_\si^\al
\circ(\Ph_\pi^\al\oti\id_\rho\oti\id_\si\oti\id_\ze)
\\
=&\,
\sum_{\xi,\eta}
\sum_{S,T}
\frac{d_\xi}{d_\pi d_\si}
\frac{d_\eta}{d_\ze d_\rho}
(1\oti \Si_{\rho,\ze}T)\Ph_\eta^\be(1\oti (\Si_{\rho,\ze}T)^*)
\circ
(1\oti S)\Ph_\xi^\al(1\oti S^*)
\\
=&\,
\sum_{\xi,\eta}
\sum_{S,T}
\frac{d_\xi}{d_\pi d_\si}
\frac{d_\eta}{d_\ze d_\rho}
\Si_{\rho,\ze}
(1\oti S\cdot T)
\big{(}
\Ph_\eta^\be\circ(\Ph_\xi^\al\oti\id_\eta)\big{)}
(1\oti (S\cdot T)^*)\Si_{\rho,\ze},
\end{align*}
where the indices $S$, $T$ runs $\ONB(\xi,\pi\cdot \si)$ and 
$\ONB(\eta,\zeta\cdot\rho)$, respectively and 
$S\cdot T\in (\xi\cdot\eta,\pi\cdot\ze\cdot\si\cdot\rho)$ 
is naturally defined via $S$, $T$. 
Let $x\in \lhGG$ and then we have 
\[
\De_{\bhGG}(x)\Si_{\rho,\ze}(S\cdot T)
=\Si_{\rho,\ze}(S\cdot T)x_{\xi,\eta}. 
\]
Hence we have 
\begin{align*}
\ps
\sim&\,
\sum_{\xi,\eta}
\sum_{S,T}
\frac{d_\xi}{d_\pi d_\si}
\frac{d_\eta}{d_\ze d_\rho}
\Si_{\rho,\ze}
(1\oti S\cdot T)
\cdot
u_\xi^{n*}\al_\xi(v_\eta^{n*})
\\
&\hspace{3cm}
\cdot
\big{(}
\ph\circ\Ph_\eta^\be\circ(\Ph_\xi^\al\oti\id_\eta)
\big{)}
\cdot\al_\xi(v_\eta^n)u_\xi^n
\cdot
(1\oti (S\cdot T)^*)
\Si_{\rho,\ze}
\\
=&\,
\sum_{\xi,\eta}
\sum_{S,T}
\frac{d_\xi}{d_\pi d_\si}
\frac{d_\eta}{d_\ze d_\rho}
\Si_{\rho,\ze}
(1\oti S\cdot T)
\cdot 
u_\xi^{n*}
\\
&\hspace{3cm}
\cdot
\big{(}
(v_\eta^{n*}(\ph\circ\Ph_\eta^\be)\cdot v_\eta^n)
\circ(\Ph_\xi^\al\oti\id_\eta)
\big{)}
\cdot u_\xi^n
\cdot
(1\oti (S\cdot T)^*)
\Si_{\rho,\ze}
\\
\sim&\,
\sum_{\xi,\eta}
\sum_{S,T}
\frac{d_\xi}{d_\pi d_\si}
\frac{d_\eta}{d_\ze d_\rho}
\Si_{\rho,\ze}
(1\oti S\cdot T)
\cdot 
u_\xi^{n*}
\cdot
\big{(}
(\ph\oti\ta_\eta)
\circ(\Ph_\xi^\al\oti\id_\eta)
\big{)}
\cdot u_\xi^n
\\
&\hspace{3cm}
\cdot
(1\oti (S\cdot T)^*)
\Si_{\rho,\ze}
\\
\sim&\,
\sum_{\xi,\eta}
\sum_{S,T}
\frac{d_\xi}{d_\pi d_\si}
\frac{d_\eta}{d_\ze d_\rho}
\Si_{\rho,\ze}
(1\oti S\cdot T)
\cdot 
(\ph\oti\ta_\xi\oti\ta_\eta)
\cdot
(1\oti (S\cdot T)^*)
\Si_{\rho,\ze}
\\
=&\,
\ph\oti\ta_\pi\oti\ta_\rho\oti\ta_\si\oti\ta_\ze.
\end{align*}
\end{proof}

\textit{Proof of Lemma \ref{lem: coc-ultra}.}
Since the discrete Kac algebra $\bhGG$ is amenable, 
we can apply Lemma \ref{lem: vanish-ultra} to the 
cocycle action $(\ga,w)$. 
Then there exists a unitary 
$c\in M_\om\oti \lhGG$ such that 
\begin{equation}\label{eq: vw}
c_{123}\ga(c)w(\id\oti\De_{\bhGG})(c^*)=1. 
\end{equation}
Set the unitaries $c^\el=c_{\cdot\oti\btr}$, 
$c^r=c_{\btr\oti\cdot}$ in $M_\om\oti\lhG$. 
By definition we have
\[
\ga_{\cdot\oti\btr}=\Ad U^*\circ \al^\om,
\quad
\ga_{\btr\oti\cdot}=\Ad V^*(\cdot\oti1).
\]
Hence applying $1\oti1\oti e_\btr\oti1\oti e_\btr$, 
$1\oti e_\btr\oti1\oti e_\btr\oti1$ to (\ref{eq: vw}), 
we have 
\begin{align*}
&c_{12}^\el U_{12}^*\al^\om(c^\el)U_{12}
w_{\cdot\oti\btr\oti\cdot\oti\btr}
(\id\oti\De)(c^{\el*})=1,
\\
&c_{12}^r V_{12}^*c_{13}^r V_{12}w_{\btr\oti\cdot\oti\btr\oti\cdot}
(\id\oti\De^\opp)(c^{r*})=1.
\end{align*}
The equalities 
\[
w_{\cdot\oti\btr\oti\cdot\oti\btr}
=
U_{12}^*\al^\om(U^*)(\id\oti\De)(U),
\quad
w_{\btr\oti\cdot\oti\btr\oti\cdot}
=
V_{12}^* V_{13}^* (\id\oti\De^\opp)(V)
\]
yield 
\[
c_{12}^\el U_{12}^*\al^\om(c^\el U^*)
(\id\oti\De)(U c^{\el*})=1,
\quad
c_{12}^r V_{12}^*c_{13}^r V_{13}^*
(\id\oti\De^\opp)(Vc^{r*})=1.
\]
Hence $c^\el U^*$ is an $\al^\om$-cocycle and $V c^{r*}$ is 
a unitary representation of $\bhG$. 
Set $v=c^\el U^*$ and $\ovl{v}=V c^{r *}$. 
Set the perturbed action $\al^v=\Ad v\circ \al^\om$. 
We claim that the unitary representation 
$\ovl{v}$ is fixed by $\al^v$, and 
then it follows that 
$W=\ovl{v}v$ is an $\al^\om$-cocycle. 
If we prove this claim, the unitary $W$ is a desired one. 
Indeed, for $x\in M$ we have 
\begin{align*}
\Ad W \circ \al(x)
=&\,
\Ad \ovl{v} \circ\al^v(x)
=
\ovl{v}v\al^\om(x) v^* \ovl{v}^*
=
\ovl{v}v U(x\oti1)U^* v^* \ovl{v}^*
\\
=&\,
\ovl{v}c^\el (x\oti1) c^{\el*} \ovl{v}^*
=
\ovl{v}(x\oti1)\ovl{v}^*
=
V(x\oti1)V^*
=
\be(x). 
\end{align*}
We prove the claim as follows. 
Applying $1\oti1\oti e_\btr\oti e_\btr\oti1$ to (\ref{eq: vw}), 
we have 
\[
c_{12}^\el U_{12}^* \al^\om(c^r) U_{12}
w_{\cdot\oti\cdot\oti \btr\oti\btr\oti\cdot}
c_{123}^*=1. 
\]
The equality 
$
w_{\cdot\oti\cdot\oti \btr\oti\btr\oti\cdot}
=
1$ 
yields
\begin{equation}\label{eq: cU}
c_{12}^\el U_{12}^* \al^\om(c^r) U_{12}c_{123}^*=1. 
\end{equation}
Again applying $1\oti e_\btr\oti1\oti1\oti e_\btr$ to 
(\ref{eq: vw}), we have 
\[
c_{12}^r V_{12}^* c_{13}^\el V_{12}
w_{\cdot\oti \btr\oti\cdot\oti\cdot\oti \btr}
c_{132}^*=1. 
\]
The equality 
\[
w_{\cdot\oti \btr\oti\cdot\oti\cdot\oti \btr}=
V_{12}^* U_{13}^* \al^\om(V)_{132} U_{13}
\]
and (\ref{eq: cU}) implies
\begin{align*}
1=&\,
c_{12}^r V_{12}^* c_{13}^\el V_{12}
V_{12}^* U_{13}^* \al^\om(V)_{132} U_{13}c_{132}^* 
\\
=&\,
c_{12}^r V_{12}^* c_{13}^\el 
U_{13}^* \al^\om(V)_{132} U_{13}c_{132}^*
\\
=&\,
c_{12}^r V_{12}^* c_{13}^\el 
U_{13}^* \al^\om(V)_{132} U_{13}
\cdot
U_{13}^*\al^\om(c^{r*})_{132}U_{13}c_{13}^{\el*}
\\
=&\,
c_{12}^r V_{12}^* c_{13}^\el 
U_{13}^* \al^\om(V c^{r*})_{132}
U_{13}c_{13}^{\el*}
\\
=&\,
\ovl{v}_{12}^* \al^v(\ovl{v})_{132}. 
\end{align*}
Therefore $\al^v(\ovl{v})=\ovl{v}_{13}$ and 
we have proved Lemma \ref{lem: coc-ultra}. 
$\hspace{4.46cm}\Box$

\section{Rohlin type theorem}
The Rohlin theorem in \cite[Theorem 6.1]{Oc1} has been a main ingredient 
to show vanishing results on 1 and 2-cohomology for strongly free 
cocycle actions of discrete amenable groups. 
Even for amenable discrete Kac algebras we can prove the 
Rohlin type theorem which is, however, not a generalization of 
the Rohlin theorem in \cite{Oc1}. 
As one difference, 
which comes from difficulty of reducing a cocycle action 
to an action, 
we give a Rohlin tower which has a good estimate only for cocycle actions 
whose 2-cocycles are very small. 
Another difference is that we treat 
not paving families but one sufficiently large projection. 
Since our classification result is deduced from 
the Evans-Kishimoto type intertwining argument, 
we do not need a model action splitting method. 
Hence it is unnecessary to utilize a paving family. 
Also even in proving vanishing results on 2-cohomology, 
we do not need such a family. 
We use the Rohlin type theorem only to find a unitary 
perturbing a 1-cocycle to a smaller 1-cocycle by using 
the Shapiro lemma. 
Although it may seem to be an incomplete form, 
it is in fact a sufficiently powerful tool for our strategy. 

\subsection{Local quantization principle}
We begin with the local quantization principle 
proved by Popa \cite[Lemma A.1.1]{Po-amen}, 
\cite[Theorem A.1.2]{Po1}. 

\begin{thm}[Popa]
Let $A\subs B$ be an inclusion of finite von Neumann algebras. 
Let $\ta$ be a faithful normal trace of $B$. 
Assume that elements $\{x_i\}_{i=1}^n\subs B$ are orthogonal to 
$A\vee (A'\cap B)$ with respect to $\ta$. 
Then for any $\vep>0$, there exists a finite index set $J$ and 
a partition of unity 
$\{q_r\}_{r\in J} \subs A$ satisfying
\[
\Big{\|}
\sum_{r\in J} q_r x_i q_r
\Big{\|}_\ta<\vep
\quad
\mbox{for all}\ 1\leq i\leq n. 
\]
\end{thm}

Let $M$ be a von Neumann algebra. 
Let $\ga\col M_\om\ra M_\om\oti \lhG$ be a strongly free semiliftable 
action with the left inverses $\{\Ph_\pi\}_{\pi\in\IG}$. 
We assume that $\ta^\om\circ\Ph_\pi=\ta^\om\oti\ta_\pi$ 
on $M_\om\oti B(H_\pi)$ for all $\pi\in\IG$. 
Take a faithful state $\ph\in M_*$ and set $\ps=\ph\circ\ta^\om$. 
Consider the inclusion $M_\om \subs M_\om\rti_\ga\bhG$. 
Then the dual state $\hps$ is tracial 
by Proposition \ref{prop: dualtrace}. 
Let $S$ be a countably generated von Neumann subalgebra of $M^\om$, 
$A=S'\cap M_\om$ and $B=M_\om\rti_\ga\bhG$. 
Since $\ga$ is strongly free, 
we have the inclusion $A\vee (A'\cap B)\subs M_\om$ 
by Lemma \ref{lem: free-relcom}. 
Let 
$E_\hga \col M_\om\rti_\ga\bhG\ra M_\om$ 
be the conditional expectation defined by averaging the dual action $\hga$, 
and then it preserves the trace $\hps$. 
Let $E_{A\vee (A'\cap B)}$ be the trace preserving conditional expectation 
from $B$ onto $A\vee (A'\cap B)$. 
Then it factors through $M_\om$ with $E_\hga$. 
Hence by definition of $E_\hga$, 
we have $E_{A\vee (A'\cap B)}(\la_{\pi_{i,j}})=0$ 
for all $\pi\neq\btr$ and $i,j\in I_\pi$. 
Let $\mF$ be a finite subset of $\IG\setm\{\btr\}$. 
We apply the local quantization principle to the above 
$A$, $B$ 
and 
$\{\la_{\pi_{i,j}}\}_{\pi\in\mF,i,j\in I_\pi}$. 
Then for any $\vep>0$, 
we get a finite partition of unity 
$\{q_r\}_{r\in J} \subs S'\cap M_\om$ satisfying
\[
\Big{\|}\sum_{r\in J} \ga(q_r) \la_{\pi_{i,j}}\ga(q_r)
\Big{\|}_{\hps}
<\vep 
\]
for all $\pi\in\mF$ and $i,j\in I_\pi$. 
Since 
$\la_{\pi_{i,j}}\ga(q_r)
=
\sum_{k\in I_\pi}\ga
(\ga_{\pi_{i,k}}(q_r))\la_{\pi_{k,j}}$
and
\begin{align*}
\hps
\big{(}
\big{(}
\ga(q_r\ga_{\pi_{i,\el}}(q_r))\la_{\pi_{\el,j}}\big{)}^*
\ga(q_r\ga_{\pi_{i,k}}(q_r))\la_{\pi_{k,j}}
\big{)} 
=&\,
d_\pi^{-1} 
\ps\big{(}
\Ph_\pi\big{(}
|q_r\ga_{\pi_{i,k}}(q_r)|^2\oti e_{\pi_{\el,k}}
\big{)}
\big{)}
\\
=&\,
d_\pi^{-1} 
(\ps\oti\ta_\pi)\big{(}
|q_r\ga_{\pi_{i,k}}(q_r)|^2\oti e_{\pi_{\el,k}}
\big{)}
\\
=&\,
\de_{\el,k} d_\pi^{-1} 
\ps\big{(}
|q_r\ga_{\pi_{i,k}}(q_r)|^2\big{)}, 
\end{align*}
we have 
\begin{align*}
&\sum_{i\in I_\pi}\Big{\|}
\sum_{r\in J} 
\ga(q_r) \la_{\pi_{i,j}}\ga(q_r)
\Big{\|}_{\hps}^2
\\
=&
\sum_{i\in I_\pi}\sum_{r\in J} 
\Big{\|}
\ga(q_r) \la_{\pi_{i,j}} \ga(q_r)
\Big{\|}_{\hps}^2
\\
=&
\sum_{i\in I_\pi}
\sum_{r\in J} 
\Big{\|}
\sum_{k\in I_\pi}
\ga(q_r\ga_{\pi_{i,k}}(q_r))\la_{\pi_{k,j}}
\Big{\|}_\hps^2
\\
=&
\sum_{r\in J} 
\sum_{i,k,\el\in I_\pi}
\hps
\big{(}
\big{(}
\ga(q_r\ga_{\pi_{i,\el}}(q_r))\la_{\pi_{\el,j}}\big{)}^*
\ga(q_r\ga_{\pi_{i,k}}(q_r))\la_{\pi_{k,j}}
\big{)}
\\
=&
\sum_{r\in J} 
\sum_{i,k,\el\in I_\pi}
\de_{\el,k} d_\pi^{-1} 
\ps\big{(}
|q_r\ga_{\pi_{i,k}}(q_r)|^2\big{)}
\\
=&
\sum_{r\in J} 
\big{\|}
(q_r\oti1_\pi)\ga_\pi(q_r)
\big{\|}_{\ps\oti\ta_\pi}^2.
\end{align*}
Thus for $\pi\in\mF$, 
\[\sum_{r\in J} 
\big{\|}
(q_r\oti1_\pi)\ga_{\pi}(q_r)
\big{\|}_{\ps\oti\ta_\pi}^2
<
d_\pi\vep^2. 
\]
Summing up the above inequality with $\pi\in \mF$, 
we obtain 
\[
\sum_{r\in J}
\sum_{\pi\in\mF}
\big{\|}
(q_r\oti1_\pi)\ga_\pi(q_r)\big
{\|}_{\ps\oti\ta_\pi}^2
<
\Big{(}\sum_{\pi\in\mF}d_\pi\Big{)}\vep^2 .
\]
We use a Chebyshev inequality as follows. 
Define an index subset 
\[
J_0=\Big{\{}r\in J 
\mid 
\sum_{\pi\in\mF}
\big{\|}
(q_r\oti1_\pi)\ga_\pi(q_r)
\big{\|}_{\ps\oti\ta_\pi}^2
<
(\sum_{\pi\in\mF}d_\pi)\vep \|q_r\|_\ps^2
\Big{\}}.
\]
For $r\in J_0$, we have
\begin{align*}
\Big{(}
\sum_{\pi\in\mF}
\big{|}
(q_r\oti1_\pi)\ga_\pi(q_r)
\big{|}_{\ps\oti\ta_\pi}
\Big)^2
\leq&
|\mF|
\sum_{\pi\in\mF}
\big{|}(q_r\oti1_\pi)\ga_\pi(q_r)
\big{|}_{\ps\oti\ta_\pi}^2
\\
\leq&
|\mF|
\sum_{\pi\in\mF}
\|\ga_\pi(q_r)\|_{\ps\oti\ta_\pi}^2
\big{\|}(q_r\oti1_\pi)\ga_\pi(q_r)\big{\|}_{\ps\oti\ta_\pi}^2
\\
=&
|\mF|\|q_r\|_\ps^2
\sum_{\pi\in\mF}
\big{\|}(q_r\oti1_\pi)\ga_\pi(q_r)\big{\|}_{{\ps\oti\ta_\pi}}^2
\\
<&
|\mF|\|q_r\|_\ps^2
\Big{(}\sum_{\pi\in\mF}
d_\pi\Big{)}
\vep \|q_r\|_{\ps}^2
\\
=&
|\mF|
\Big{(}\sum_{\pi\in\mF}d_\pi\Big{)}
\vep |q_r|_{\ps}^2. 
\end{align*}
Hence the following inequality holds. 
\[
\sum_{\pi\in\mF}
\big{|}
(q_r\oti1_\pi)\ga_\pi(q_r)
\big{|}_{\ps\oti\ta_\pi}
\leq
\Big{(}
|\mF|
\sum_{\pi\in\mF}
d_\pi\Big{)}^{1/2} 
\vep^{1/2}|q_r|_\ps.
\]
On the size of $\sum_{r\in J\setm J_0}q_r$, we have
\begin{align*}
\sum_{r\in J\setm J_0}
|q_r|_\ps
<&
\Big{(}
\sum_{\pi\in\mF}d_\pi\Big{)}^{-1}\vep^{-1}
\sum_{r\in J\setm J_0}
\sum_{\pi\in\mF}
\big{\|}
(q_r\oti1_\pi)\ga_\pi(q_r)\big{\|}_{\ps\oti\ta_\pi}^2
\\
\leq&
\Big{(}
\sum_{\pi\in\mF}d_\pi\Big{)}^{-1}\vep^{-1}
\sum_{r\in J}
\sum_{\pi\in\mF}
\big{\|}
(q_r\oti1_\pi)\ga_\pi(q_r)\big{\|}_{\ps\oti\ta_\pi}^2
\\
<&
\Big{(}\sum_{\pi\in\mF}d_\pi\Big{)}^{-1}\vep^{-1}
\Big{(}\sum_{\pi\in\mF}d_\pi\Big{)}\vep^2 
\\
=&
\vep.
\end{align*}
We summarize these arguments. 
\begin{lem}
Let $M$ be a von Neumann algebra. 
Let $\ga$ be a strongly free semiliftable action of $\bhG$ 
on $M_\om$ whose left inverse $\Ph_\pi$ satisfies 
$\ta^\om\circ\Ph_\pi=\ta^\om\oti\ta_\pi$ for all $\pi\in\IG$. 
Let $\ph$ be a faithful normal state on $M$ and 
set $\ps=\ph\circ\ta^\om$. 
Then for any countably generated von Neumann algebra 
$S\subs M^\om$, 
any finite subset $\mF\subs\IG$, 
$\btr\nin\mF$ and any $0<\vep<1$, 
there exists $n\in \N$ and 
a partition of unity $\{q_r\}_{r=0}^n \subs S'\cap M_\om$ 
with the following properties. 
\begin{enumerate}

\item $|q_0|_\ps<\vep$. 

\item 
$
\displaystyle
\sum_{\pi\in\mF}
\big{|}
(q_r\oti1_\pi)\ga_\pi(q_r)\big{|}_{\ps\oti\ta_\pi}
<
\vep |q_r|_\ps$ 
for all $1\leq r\leq n$ and $\pi\in\mF$. 
\end{enumerate}
\end{lem}

We can strengthen this result as follows. 

\begin{lem}\label{lem: partition}
Let $M$ be a von Neumann algebra. 
Let $\ga$ be a strongly free semiliftable action of $\bhG$ 
on $M_\om$ whose left inverse $\Ph_\pi$ satisfies 
$\ta^\om\circ\Ph_\pi=\ta^\om\oti\ta_\pi$ for all $\pi\in\IG$. 
Let $\ph$ be a faithful normal state on $M$ and 
set $\ps=\ph\circ\ta^\om$. 
Then for any countably generated von Neumann algebra 
$S\subs M^\om$, 
any finite subset $\mF\subs\IG$, 
$\btr\nin\mF$ and any $0<\de<1$, 
there exists $n\in \N$ and 
a partition of unity $\{e_r\}_{r=0}^n \subs S'\cap M_\om$ 
with the following properties. 
\begin{enumerate}

\item $|e_0|_{\ps}\leq \de$. 

\item 
$(e_r\oti1_\pi)\ga_\pi(e_r)=0$ 
for all $1\leq r\leq n$ and $\pi\in \mF$.
\end{enumerate} 
\end{lem}

\begin{proof}
This is proved by a similar argument to one in \cite{Oc1} as follows. 
We may assume $S\subs M^\om$ is $\ga$-invariant by considering 
a von Neumann algebra generated by 
$S$ and $\ga_{\rho_{k,\el}}(S)$ for 
all $\rho\in\IG$ and $k,\el\in I_\rho$. 

\textbf{Step A}. 
Let $\mu>0$ and $f\in\Proj(S'\cap M_\om)$, $f\neq0$. 
We show that there exists $f'\in \Proj(S'\cap M_\om)$, $0\neq f'\leq f$, 
such that 
\[
\sum_{\pi\in\mF}
\big{|}
(f'\oti1_\pi)\ga_\pi(f')
\big{|}_{\ps\oti\ta_\pi} < 2\mu |f'|_{\ps}
,\]
for all $\pi\in \mF$. 

Let $\ovl{S}$ be a von Neumann subalgebra in $M^\om$ 
which is generated by $S$ and 
$\ga_{\rho_{i,j}}(f)$ for all $\rho\in\IG$ and $i,j\in I_\rho$. 
By the previous lemma, 
there exists a partition of unity $f_0,f_1,\dots,f_m$ 
in $\ovl{S}'\cap M_\om$ such that 
\begin{enumerate}

\item $|f_0|_{\ps}\leq 2^{-1} |f|_{\ps}$,

\item 
$\displaystyle
\sum_{\pi\in\mF}
\big{|}(f_i\oti1_\pi)\ga_\pi(f_i)\big{|}_{\ps\oti\ta_\pi}
<\mu |f|_{\ps} |f_i|_{\ps}$, $1\leq i\leq m$.

\end{enumerate}
Let $\bar{f_i}=ff_i\in\Proj(S'\cap M_\om)$. 
Since $f_i\oti1_\rho$ commutes $\ga_\rho(f)$ for all $\rho\in\IG$, 
we have 
$(\bar{f_i}\oti1_\rho)\ga_\rho(\bar{f_i})
=
(f\oti1_\rho)\ga_\rho(f)(f_i\oti1_\rho)\ga_\rho(f_i). 
$
Hence we have 
$
\big{|}
(\bar{f_i}\oti1_\rho)\ga_\rho(\bar{f_i})
\big{|}_{\ps\oti\ta_\pi}
\leq
\big{|}(f_i\oti1_\rho)\ga_\rho(f_i)
\big{|}_{\ps\oti\ta_\pi}.
$
Suppose that 
for each $i=1,\dots,m$ 
\[
\sum_{\pi\in\mF}
\big{|}(\bar{f_i}\oti1_\pi)\ga_\pi(\bar{f_i})
\big{|}_{\ps\oti\ta_\pi}
\geq 
2\mu |\bar{f_i}|_{\ps}.
\] 
Then 
\begin{align*}
\sum_{i=1}^m \sum_{\pi\in\mF}
\big{|}(f_i\oti1_\pi)\ga_\pi(f_i)\big{|}_{\ps\oti\ta_\pi}
\geq&
\sum_{i=1}^m \sum_{\pi\in\mF}
\big{|}(\bar{f_i}\oti1_\rho)
\ga_\rho(\bar{f_i})\big{|}_{\ps\oti\ta_\pi}
\\
\geq&
2\mu \sum_{i=1}^m |\bar{f_i}|_{\ps}
\\
=&
2\mu|(1-f_0)f|_{\ps}
\\
\geq&
2\mu(|f|_{\ps}-|f_0|_{\ps})
\\
\geq&
\mu |f|_{\ps}.
\end{align*}
On the other hand, we have 
\begin{align*}
\sum_{i=1}^m \sum_{\pi\in\mF}
\big{|}
(f_i\oti1_\pi)\ga_\pi(f_i)\big{|}_{\ps\oti\ta_\pi}
<&\mu
|f|_{\ps}\sum_{i=1}|f_i|_{\ps}
\\
\leq&
\mu|f|_{\ps}.
\end{align*}
This derives a contradiction and hence for some $i\in\{1,\dots,m\}$, 
the equality 
\[
\sum_{\pi\in\mF}
\big{|}
(\bar{f_i}\oti1_\pi)\ga_\pi(\bar{f_i})
\big{|}_{\ps\oti\ta_\pi}
<
2\mu |\bar{f_i}|_{\ps}
\]
holds and we take $f'=\bar{f_i}$. 

\textbf{Step B}. 
We show that for any $f\in\Proj(S'\cap M_\om)$ and 
any $\mu>0$, there exists $e\in\Proj(S'\cap M_\om)$ with
\begin{enumerate}

\item $e\leq f$,

\item 
$\big{|}(e\oti1_\pi)\ga_\pi(e)\big{|}_{\ps\oti\ta_\pi}
\leq \mu |e|_{\ps}$ for all $\pi\in\mF$,

\item 
$|e|_{\ps}\geq(1+\sum_{\pi\in\mF\cup\ovl{\mF}}2d_{\pi}^2)^{-1}
|f|_{\ps}$. 

\end{enumerate}
Set $T_\pi(x)=\Ph_{\ovl{\pi}}(x\oti1_{\ovl{\pi}})$ for $x\in M^\om$. 
Note $T_\pi(S'\cap M_\om)\subs S'\cap M_\om$ for all $\pi\in\IG$ 
as is seen below. 
Let $x\in S'\cap M_\om$. 
Since $T_\pi$ preserves $M_\om$, $T_\pi(x)\in M_\om$. 
Take an element $y\in S$. 
The $\ga$-invariance of $S$ implies 
\[
yT_\pi(x)
=
\Phi_\pi(\ga_\pi(y)(x\oti1_\pi))
=
\Phi_\pi((x\oti1_\pi)\ga_\pi(y))
=
T_\pi(x)y. 
\]
This shows $T_\pi(x)\in S'\cap M_\om$. 
Now the family of projections $e\in S'\cap M_\om$ satisfying 
(1) and (2) 
is not empty and inductively ordered, so let $e$ be maximal with these 
properties. 
We show that $e$ also satisfies 

\begin{equation}\label{cover}
e\vee \bigvee_{\pi\in\mF}\,s(T_\pi(e))\vee
\bigvee_{\pi\in\mF}\,s(T_{\bar{\pi}}(e))\vee(1-f)=1.
\end{equation}

Otherwise let $e'$ be a nonzero projection in $S'\cap M_\om$ orthogonal to 
the left member of (\ref{cover}). 
We claim that $(e\oti1_\pi)\ga_\pi(e')=0=(e'\oti1_\pi)\ga_\pi(e)$ 
for $\pi\in\mF\cup\ovl{\mF}$. 
Since 
\[
\Ph_\pi\big{(}\ga_\pi^\om(e')(e\oti1_\pi)\ga_\pi^\om(e')\big{)}
=
e'\Ph_\pi(e\oti1_\pi)e'
=
e' T_\pi(e)e'
=
0,
\]
we have $(e\oti1_\pi)\ga_\pi(e')=0$ 
for $\pi\in\mF\cup\ovl{\mF}$ by faithfulness of $\Ph_\pi$. 
Since 
\begin{align*}
\Ph_\pi\big{(}\ga_\pi(e)(e'\oti1_\pi)\ga_\pi(e)\big{)}
=&\,
e\Ph_\pi\big{(}e'\oti1_\pi\big{)}e
\\
=&\,
e(1\oti T_{\opi,\pi}^*)
(\ga_\opi(e')\oti1_\pi )(1\oti T_{\opi,\pi})e
\\
=&\,
(1\oti T_{\opi,\pi}^*)
\big{(}
(e\oti1_\opi)\ga_\opi(e')(e\oti1_\pi)\oti1_\pi
\big{)} 
(1\oti T_{\opi,\pi})
\\
=&\,
0,
\end{align*}
we also have $(e'\oti1_\pi)\ga_\pi(e)=0$ 
for $\pi\in\mF\cup\ovl{\mF}$. 
By Step A, under $e'$ we can find a nonzero projection $e''\in S'\cap M_\om$ 
with 
$
\big{|}
(e''\oti1_\pi)\ga_\pi(e'')\big{|}_{\ps\oti\ta_\pi}
\leq
\mu|e''|_\ps
$. 
Then $e''$ satisfies $e''\per e$, $e''\leq f$ and 
$(e\oti1_\pi)\ga_\pi(e'')=0=(e''\oti1_\pi)\ga_\pi(e)$ 
for all $\pi\in\mF\cup\ovl{\mF}$. 
Then we have 
\begin{align*}
\big{|}
(\bar{e}\oti1_\pi)\ga_\pi(\bar{e})
\big{|}_{\ps\oti\ta_\pi}
=&\,
\big{|}
(e\oti1_\pi)\ga_\pi(e)
+
(e''\oti1_\pi)\ga_\pi(e'')
\big{|}_{\ps\oti\ta_\pi}
\\
=&\,
\big{|}
(e\oti1_\pi)\ga_\pi(e)
\big{|}_{\ps\oti\ta_\pi}
+
\big{|}
(e''\oti1_\pi)\ga_\pi(e'')
\big{|}_{\ps\oti\ta_\pi}
\\
\leq&
\mu |e|_{\ps}+\mu |e''|_{\ps}
\\
=&
\mu |\bar{e}|_{\ps}.
\end{align*}
This is a contradiction to the maximality of $e$, 
and hence (\ref{cover}) holds. 
Now we estimate the size of $e$. 
Here we make use of tensor products to 
treat translations of projections. 
We index $\mF\cup \ovl{\mF}$ as $\{\pi_1,\cdots,\pi_N\}$.
Consider the tensor product algebra 
$\wdt{M}^\om =M^\om\oti 
B(\oti_{k=1}^N H_{\pi_k})$. 
Set a product state 
$\wdt{\ps}=\ps\oti \oti_{k=1}^N \ta_{\pi_k}$, 
which is tracial on 
$M_\om\oti B(\oti_{k=1}^N H_{\pi_k})$. 
We regard 
$M^\om$ and 
$M^\om \oti B(H_{\pi_k})$ as subalgebras 
of $\wdt{M}^\om$ via the natural embedding. 
For each $\pi_k$, take a finite group 
$\meU_k\cong \Z_{d_{\pi_k}}\!\rti \Z_{d_{\pi_k}}$ 
consisting of unitaries in $B(H_{\pi_k})$ 
which acts on $H_{\pi_k}$ irreducibly. 
Our claim is the following one. 
\begin{equation}\label{cover2}
e\vee
\bigvee_{k=1}^N \bigvee_{v\in\meU_k}
(1\oti v)
\ga_{\pi_k}(e)
(1\oti v^*)
\vee(1-f)
=
1.
\end{equation}
If not so, we have a nonzero projection $p\in \wdt{M}^\om$ which is 
orthogonal to the above left projection. 
Hence $p\leq f-e$ and satisfies 
$
p\cdot(1\oti v)\ga_{\pi_k}(e)
(1\oti v^*)=0
$ 
for all $k=1,\dots,N$ and $v\in\meU_k$. 
Irreducibility of the action of $\meU_k$ implies 
the equality 
$
\ta_{\pi_k}=d_{\pi_k}^{-2}\sum_{v\in \meU_k}\Ad v$ on 
$B(H_{\pi_k})$. 
Hence 
$p T_{\opi_k}(e)=p\cdot(\id\oti\ta_{\pi_k})
(\ga_{\pi_k}(e))=0$ 
for all $k$. 
It shows that $p$ must be orthogonal to 
$e$, $1-f$ and $\vee_{\pi\in\mF\cup\ovl{\mF}} \,s(T_\pi(e))$, 
but this contradicts to (\ref{cover}). 
Therefore the above claim holds. 
Applying the product trace $\wdt{\ps}$ to the both sides of (\ref{cover2}), 
we obtain 
\begin{align*}
1
\leq&\,
|e|_{\ps}+|1-f|_{\ps}\\
&\quad
+
\sum_{k=1}^N \sum_{v\in\meU_k}
(\ps\oti\ta_{\pi_k})\big{(}
(1\oti v)
\ga_{\pi_k}(e)
(1\oti v^*)
\big{)}
\\
\leq&\,
|e|_{\ps}+1-|f|_{\ps}
+
\sum_{k=1}^N \sum_{v\in\meU_k}
(\ps\oti\ta_{\pi_k})\big{(}\ga_{\pi_k}(e)\big{)}
\\
=&\,
1-|f|_{\ps}
+
(1+\sum_{k=1}^N 2d_{\pi_k}^2)|e|_{\ps}. 
\end{align*}
Hence  
\[
|e|_\ps\geq (1+\sum_{k=1}^N 2d_{\pi_k}^2)^{-1}|f|_\ps. 
\]

\textbf{Step C}. 
Let $q\in\N$ be such that 
$(1-(1+\sum_{\pi\in\mF\cup\ovl{\mF}} 2d_{\pi}^2)^{-1})^q<\de$. 
We show that for any $\mu>0$ there exists a partition of unity 
$\{e_k\}_{k=0}^{q} \subs S'\cap M_\om$ such that 
\begin{enumerate}

\item $|e_0|_\ps\leq \de$. 

\item $\big{|}(e_k\oti1_{\pi})\ga_{\pi}(e_k)
\big{|}_{\ps\oti\ta_\pi}
<\mu |e_k|_\ps$, 
for all $k=1,\dots,q$ and $\pi\in\mF$. 

\end{enumerate}

Set $f_1=1$. According to Step B, 
we construct projections $e_k$ and $f_{k+1}$ successively 
for $k=1,\dots,q$ such that 
\begin{enumerate}

\item $e_k\leq f_k$,

\item $f_{k+1}=f_k-e_k$,

\item 
$\big{|}(e_k\oti1_\pi)\ga_\pi(e_k)\big{|}_{\ps\oti\ta_\pi}
\leq \mu |e_k|_\ps$ 
for all $\pi\in \mF$. 

\item 
$|e_k|_\ps\geq (1+\sum_{\pi\in\mF\cup\ovl{\mF}} 2d_{\pi}^2)^{-1}
|f_k|_\ps$. 

\end{enumerate}
Then we have 
$|f_{q+1}|_\ps\leq 
(1-(1+\sum_{\pi\in\mF\cup\ovl{\mF}} 2d_{\pi}^2)^{-1})^q<\de$ 
and letting $e_0=f_{q+1}$, Step C is proved. 

\textbf{Step D}. 
We finish the proof by using the Index Selection Trick. 

Note that the partition number $q$ depends not on $\mu$ but on 
$\de$ and $\mF$. 
Letting $\mu=1/n$ and take a partition of unity 
$\{e_k^n\}_{k=0}^q\subs S'\cap M_\om$ for each $n\in\N$ such that 
\begin{enumerate}

\item $|e_0^{n}|_\ps \leq \de$. 

\item $\big{|}(e_k^n\oti1_{\pi})\ga_{\pi}(e_k^n)
\big{|}_{\ps\oti\ta_\pi}
<(1/n) |e_k^n|_\ps$, 
for all $k=1,\dots,q$ and $\pi\in\mF$. 
\end{enumerate}

Set the elements 
$\ovl{e}_i=(e_i^n)_n$ in $\el^\infty(\N,S'\cap M_\om)$ 
for $0\leq i\leq q$. 
Then apply the Index Selection Trick for 
$\meC=C^*(\{\ovl{e}_i\}_{i=0}^q)$ and 
$\mB=\{\ga_\pi\}_{\pi\in\IG}$. 
Let $\Ps$ be the index selection map with respect to them. 
Set $e_i=\Ps(\ovl{e}_i)$ which is in $S'\cap M_\om$. 
Then we have $|e_0|_\ps\leq \de$ and 
\begin{align*}
\big{|}(e_i\oti1_{\pi})\ga_{\pi}(e_i)
\big{|}_{\ps\oti\ta_\pi}
=&\,
(\ps\oti\ta_\pi)
\big{(}
\big{|}(e_i\oti1_{\pi})\ga_{\pi}(e_i)
\big{|}
\big{)}
\\
=&\,
(\ph\circ\ta^\om\circ\Ps\oti\ta_\pi)
\big{(}
\big{|}(\ovl{e}_i\oti1_{\pi})\ga_{\pi}(\ovl{e}_i)
\big{|}
\big{)}
\\
=&\,
\lim_{n\to\om}
(\ph\circ\ta^\om\oti\ta_\pi)
\big{(}
\big{|}(e_i^n\oti1_{\pi})\ga_{\pi}(e_i^n)
\big{|}
\big{)}
\\
\leq&\,
\lim_{n\to\om}1/n
\\
=&\,0.
\end{align*}
for $\pi\in \mF$.
\end{proof}

\subsection{Tower bases and diagonal elements}

In the previous subsection, 
we have obtained a projection in $M_\om$ which behaves 
like a tower base with respect to an action $\ga$ 
as in the case of group actions. 
Here one shall note that it gives projections in not 
$M_\om$ but von Neumann algebras $M_\om\oti B(H_\pi)$, $\pi\in\IG$. 
In the following lemma, 
we clarify the special properties of such a projection 
in a general situation. 

\begin{lem}\label{lem: towerbase}
Let $M$ be a von Neumann algebra and 
$(\al,u)$ a cocycle action of $\bhG$ on $M$. 
Let $\mK$ be a finite subset of $\bhG$. 
If a projection $e$ in $M$ satisfies 
$[e\oti1\oti1,u]=0$ and 
$(e\oti1_\rho)\al_\rho(e)=0$ for 
all $\rho\in \ovl{\mK}\cdot\mK\setm{\btr}$, 
then 

\begin{enumerate}

\item 
$(\al_\pi(e)\oti1_\si)(\al_\si(e))_{1,3}=0$ 
for all $\pi\neq\si\in\mK$. 

\item 
$(e\oti1_{\ovl{\pi}}\oti1_\pi)
(\al_{\ovl{\pi}}\oti\id)(\al_\pi(e))
=u_{\ovl{\pi},\pi}
(e\oti {}_\opi\De_\pi(e_\btr))
u_{\ovl{\pi},\pi}^*
$ 
for all $\pi\in\mK$.

\item 
The element $d_\pi^2\Ph_{\ovl{\pi}}(e\oti1_{\ovl{\pi}})$ 
is a projection which is equal to 
\[
q_\pi=\inf\{q\in \Proj(M)\mid 
(q\oti1_\pi)\al_\pi(e)=\al_\pi(e), 
\ 
[q\oti1\oti1,u]=0\}. 
\]

\item 
For all $\pi\in\mK$, 
$d_\pi^2 (\al_\pi(e)\oti1_{\ovl{\pi}}) 
u_{\pi, \ovl{\pi}}(1\oti T_{\pi,\opi}T_{\pi,\opi}^*)
u_{\pi, \ovl{\pi}}^*
(\al_\pi(e)\oti1_{\ovl{\pi}})
=\al_\pi(e)\oti1_{\ovl{\pi}}.
$ 
\end{enumerate} 

\end{lem}

\begin{proof}
(1). 
It suffices to show 
$(\Ph_\pi\oti\id)(\al_\pi(e)_{12}\al_\si(e)_{13}\al_\pi(e)_{12})=0$. 
Since $\pi\neq\si$, $\btr\nprec\opi\cdot\si$ 
and then
\begin{align*}
(e\oti1_\opi\oti1_\si)\al_\opi(\al_\si(e))
=&\,
(e\oti1_\opi\oti1_\si)u_{\opi,\si}(\id\oti{}_\opi\De_\si)(\al(e))
u_{\opi,\si}^*
\\
=&\,
u_{\opi,\si}
(e\oti1_\opi\oti1_\si)(\id\oti{}_\opi\De_\si)(\al(e))
u_{\opi,\si}^*
\\
=&\,
\sum_{\rho\prec\opi\cdot\si}
u_{\opi,\si}
(e\oti1_\opi\oti1_\si)(\id\oti{}_\opi\De_\si)(\al_\rho(e))
u_{\opi,\si}^*
\\
=&\,
0. 
\end{align*}
Using this equality, we have 
\begin{align*}
&(\Ph_\pi\oti\id)(\al_\pi(e)_{12}\al_\si(e)_{13}\al_\pi(e)_{12})
\\
=&\,
(e\oti1_\si)
(\Ph_\pi\oti\id)(\al_\si(e)_{13})
(e\oti1_\si)
\\
=&\,
(e\oti1_\si)
(1\oti T_{\opi,\pi}^*\oti1_\si)(u_{\opi,\pi}^*\oti1_\si)
\al_\opi(\al_\si(e))_{124}
(u_{\opi,\pi}\oti1_\si)
\\
&\hspace{3cm}\cdot
(1\oti T_{\opi,\pi}\oti1_\si)
(e\oti1_\si)
\\
=&\,
(1\oti T_{\opi,\pi}^*\oti1_\si)(u_{\opi,\pi}^*\oti1_\si)
((e\oti1_\opi\oti1_\si)\al_\opi(\al_\si(e)))_{124}
(u_{\opi,\pi}\oti1_\si)
\\
&\hspace{3cm}\cdot
(1\oti T_{\opi,\pi}\oti1_\si)
(e\oti1_\si)
\\
=&\,
0. 
\end{align*}
Since $\Ph_\pi$ is faithful, we have $\al_\pi(e)_{12}\al_\si(e)_{13}=0$. 

(2). It is verified as follows,
\begin{align*}
(e\oti1_\opi\oti1_\pi)\al_\opi(\al_\pi(e))
=&\,
(e\oti1_\opi\oti1_\pi)
u_{\opi,\pi}(\id\oti{}_\opi\De_\pi)(\al(e))u_{\opi,\pi}^*
\\
=&\,
\sum_{\rho\prec\opi\cdot\pi}
u_{\opi,\pi}(\id\oti{}_\opi\De_\pi)((e\oti1_\rho)\al_\rho(e))
u_{\opi,\pi}^*
\\
=&\,
u_{\opi,\pi}(\id\oti{}_\opi\De_\pi)((e\oti e_\btr)\al_\btr(e))
u_{\opi,\pi}^*
\\
=&\,
u_{\opi,\pi}(e\oti{}_\opi\De_\pi(e_\btr))
u_{\opi,\pi}^*.
\end{align*}

(3). 
Applying the map $\Ph_{\ovl{\pi}}\oti\id$ to 
the equality in (2), by Lemma \ref{lem: left inverse} 
we have 
\[
(d_\pi^2\Ph_{\opi}(e\oti1_{\opi})\oti1_\pi)\al_\pi(e)
=\al_\pi(e).
\]
Set $p_\pi=d_\pi^2\Ph_{\opi}(e\oti1_{\opi})$. 
Since $e$ commutes with $u$, so does with $\al(u)$. 
Hence $p_\pi$ commutes with $u$. 
Then 
\begin{align*}
p_\pi^2
=&\,
p_\pi d_\pi^2 \Ph_{\opi}(e\oti1_{\opi})
\\
=&\,
p_\pi d_\pi^2 (1\oti T_{\pi,\opi}^*)
u_{\pi,\opi}^* 
(\al_\pi(e)\oti1_\opi)
u_{\pi,\opi}
(1\oti T_{\pi,\opi})
\\
=&\,
d_\pi^2 (1\oti T_{\pi,\opi}^*)
u_{\pi,\opi}^* 
((p_\pi\oti1_\pi)\al_\pi(e)\oti1_\opi)
u_{\pi,\opi}
(1\oti T_{\pi,\opi})
\\
=&\,
d_\pi^2 (1\oti T_{\pi,\opi}^*)
u_{\pi,\opi}^* 
(\al_\pi(e)\oti1_\opi)
u_{\pi,\opi}
(1\oti T_{\pi,\opi})
\\
=&\,
p_\pi.
\end{align*}
Hence $p_\pi$ is a projection and 
$p_\pi\geq q_\pi$. 
The inequality $p_\pi\leq q_\pi$ is verified as follows, 
\begin{align*}
q_\pi p_\pi 
=&\,
q_\pi d_\pi^2 \Ph_{\opi}(e\oti1_{\opi})
\\
=&\,
d_\pi^2 q_\pi
(1\oti T_{\pi,\opi}^*)
u_{\pi,\opi}^* 
(\al_\pi(e)\oti1_\opi)
u_{\pi,\opi}
(1\oti T_{\pi,\opi})
\\
=&\,
d_\pi^2 
(1\oti T_{\pi,\opi}^*)
u_{\pi,\opi}^* 
((q_\pi\oti1)\al_\pi(e)\oti1_\opi)
u_{\pi,\opi}
(1\oti T_{\pi,\opi})
\\
=&\,
d_\pi^2 
(1\oti T_{\pi,\opi}^*)
u_{\pi,\opi}^* 
(\al_\pi(e)\oti1_\opi)
u_{\pi,\opi}
(1\oti T_{\pi,\opi})
\\
=&\,
p_\pi. 
\end{align*}

(4). 
Since 
$d_\pi^2\Ph_{\ovl{\pi}}(e\oti1_{\ovl{\pi}})
=
d_\pi^2 (1\oti T_{\pi,\opi}^*)u_{\pi,\opi}^*
(\al_\pi(e)\oti1_\opi)u_{\pi,\opi}(1\oti T_{\pi,\opi})
$ is a projection, 
the operator 
$v:=d_\pi(\al_\pi(e)\oti1_\opi)u_{\pi,\opi}
(1\oti T_{\pi,\opi})$ is a partial isometry. 
Hence $f:=vv^*$ is a projection. 
Clearly we have $f\leq \al_\pi(e)\oti 1_\opi$. 
In fact they are equal as is shown below, 
\begin{align*}
(\Ph_\pi\oti\id)(f)
=&\,
d_\pi^2 
(\Ph_\pi\oti\id)
((\al_\pi(e)\oti1_\opi)
u_{\pi,\opi}(1\oti{}_\pi\De_\opi(e_\btr))
u_{\pi,\opi}^*(\al_\pi(e)\oti1_\opi))
\\
=&\,
d_\pi^2
(e\oti1_\opi)
(\Ph_\pi\oti\id)
(u_{\pi,\opi}(1\oti{}_\pi\De_\opi(e_\btr))
u_{\pi,\opi}^*)
(e\oti1_\opi)
\\
=&\,
e\oti1_\opi
\\
=&\,
(\Ph_\pi\oti\id)(\al_\pi(e)\oti1_\opi), 
\end{align*}
where we have used Lemma \ref{lem: left inverse}. 
The faithfulness of $\Ph_\pi$ yields $f=\al_\pi(e)\oti1_\opi$. 
\end{proof}

We call a projection $e$ a \textit{tower base} of the tower along with $K$ 
if $e$ satisfies the conditions from (1) to (4) in this lemma. 

\begin{defn}
Let $(\al,u)$ be a cocycle action of $\bhG$ on a von Neumann algebra $M$. 
The \textit{diagonal} of $u$ is 
the element $a$ in $M\oti\lhG$ defined by 
\[
(a\oti1)(1\oti\De(e_\btr))=u(1\oti\De(e_\btr)). 
\]
\end{defn}
The diagonal $a\in M\oti \lhG$ has the following 
explicit form. 
For all $\pi\in\IG$, 
\[
a_\pi=d_\pi \ep_\pi (1\oti1_\pi\oti T_{\opi,\pi}^*)
(u_{\pi,\opi}\oti1_\pi)
(1\oti T_{\pi,\opi}\oti1_\pi),
\]
where $\ep_\pi\in \{\pm1\}$ is defined in \S2.2. 

\begin{lem}\label{lem: diagonal}
Let $a_\pi$, 
$\pi\in\IG$, be as above. Then one has 
\begin{enumerate}

\item 
$(\al_\opi\oti\id)(a_\pi^*)(a_\opi\oti1_\pi)(1\oti T_{\opi,\pi})
=1\oti T_{\opi,\pi}$,

\item 
$(\id\oti \ta_\pi)(a_\pi^* a_\pi)=1$,

\item
$\Ph_\pi(a_\pi a_\pi^*)=1$, 

\item
$d_\pi^2 (\Ph_\opi\oti\id_\pi)(x\oti \De(e_\btr))
=a_\pi^* \al_\pi(x) a_\pi$ for all $x\in M$, 

\item 
$a_\pi^* a_\pi=d_\pi^2(\Ph_\opi\oti\id_\pi)(1\oti {}_\opi\De_{\pi}(e_\btr))$.

\end{enumerate}

\end{lem}

\begin{proof}
(1). Since $u$ is a 2-cocycle, we have 
\begin{align*}
&(1\oti1\oti\De(e_\btr))
(\al(a^*)\oti1)(a\oti1\oti1)(1\oti \De(e_\btr)\oti1)
\\
=&\,(1\oti1\oti\De(e_\btr))
\al(u^*)(u\oti1)(1\oti \De(e_\btr)\oti 1)
\\
=&\,(1\oti1\oti\De(e_\btr))
\al(u^*)(u\oti1)(\id\oti\De\oti\id)(u)(1\oti \De(e_\btr)\oti1)
\\
=&\,(1\oti1\oti\De(e_\btr))
(\id\oti\id\oti\De)(u)(1\oti \De(e_\btr)\oti1)
\\
=&\,
(1\oti1\oti\De(e_\btr))(1\oti \De(e_\btr)\oti1). 
\end{align*}
Applying $\id\oti\id\oti\id\oti\vph$ to the both sides, we obtain 
the desired equality by using $(\id\oti\vph)(\De(e_\btr))=1$. 

(2). 
It is verified as 
\begin{align*}
(\id\oti\ta_\pi)(a_\pi^*a_\pi)
=&\,
(1\oti T_{\pi,\opi}^*)(a_\pi^*a_\pi\oti1_\opi)(1\oti T_{\pi,\opi})
\\
=&\,
(1\oti T_{\pi,\opi}^*)(u_{\pi,\opi}^*u_{\pi,\opi})(1\oti T_{\pi,\opi})
\\
=&\,
(1\oti T_{\pi,\opi}^*)(1\oti T_{\pi,\opi})
\\
=&\,
1. 
\end{align*}

(3). 
It is verified as 
\begin{align*}
\Ph_\pi(a_\pi a_\pi^*)
=&\,
d_\pi^2
\Ph_\pi
\big{(}
(\id\oti\id\oti\ta_\opi)
((a\oti1)(1\oti\De(e_\btr))(a^*\oti1))
\big{)}
\\
=&\,d_\pi^2
(\Ph_\pi\oti\ta_\opi)(u_{\pi,\opi}(1\oti\De(e_\btr))u_{\pi,\opi}^*)
\\
=&\,
1, 
\end{align*}
where we have used Lemma \ref{lem: left inverse}. 

(4). 
It is verified as
\begin{align*}
&d_\pi^2 (\Ph_\opi\oti\id)(x\oti{}_\opi\De_\pi(e_\btr))
\\
=&\,
d_\pi^2
(1\oti T_{\pi,\opi}^*\oti1_\pi)
(u_{\pi,\opi}^*\oti1_\pi)
(\al_\pi(x)\oti{}_\opi\De_\pi(e_\btr))
(u_{\pi,\opi}\oti1_\pi)
(1\oti T_{\pi,\opi}\oti1_\pi)
\\
=&\,
d_\pi^2
(1\oti T_{\pi,\opi}^*\oti1_\pi)
(a_{\pi}^*\oti1_\opi\oti1_\pi)
(\al_\pi(x)\oti{}_\opi\De_\pi(e_\btr))
(a_{\pi}\oti1_\opi\oti1_\pi)
(1\oti T_{\pi,\opi}\oti1_\pi)
\\
=&\,
a_\pi^*\al_\pi(x)a_\pi. 
\end{align*}

(5). 
It is obtained by putting $x=1$ in (4). 
\end{proof}

With diagonals, we obtain the following result 
for a tower base. 

\begin{lem}\label{lem: matrixunit}
If a projection $e\in M$ satisfies the condition of 
Lemma \ref{lem: towerbase}, then 

\begin{enumerate}

\item
$d_\pi^2(\al_\pi(e)a_\pi\oti1_\opi)(1\oti{}_\pi\De_{\opi}(e_\btr))
(a_\pi^*\al_\pi(e)\oti1_\opi)=\al_\pi(e)\oti1_\opi$. 

\item
$\al_\pi(e)a_\pi a_\pi^* \al_\pi(e)=\al_\pi(e)$, 
in particular $\al_\pi(e)a_\pi$ is a partial isometry. 

\item 
$d_\pi^2(\Ph_\opi\oti\id_\pi)(e\oti_{\opi}\De_{\pi}(e_\btr))
=a_\pi^*\al_\pi(e)a_\pi$ is a projection. 

\item 
$\big{(}d_\pi \al_\pi(e)a_\pi \big{)}_{i,k}
\big{(} d_\pi a_\pi^*\al_\pi(e)\big{)}_{\el,j}
=\de_{k,\el}
d_\pi \big{(}\al_\pi(e)\big{)}_{i,j}$ 
for all $\pi\in\IG$ and $i,j, k, \el\in I_\pi$. 

\item 
Decompose $a^*\al(e)a_K$ as 
\[
a^*\al(e)a_K
=\sum_{\pi\in\meK}\sum_{i,j\in I_\pi}
d_\pi^{-1} f_{\opi_{i,j}}\oti e_{\pi_{i,j}}.
\]
Then one has 
\[
f_{\opi_{i,j}}^*=f_{\opi_{j,i}}
,\quad 
f_{\opi_{i,j}}f_{\orho_{k,\el}}
=\de_{\pi,\rho}\de_{j,k}f_{\opi_{i,\el}}
\]
for all $\pi,\rho\in\IG$, 
$i,j\in I_\pi$ and $k,\el\in I_\rho$. 
\end{enumerate}
\end{lem}

\begin{proof}
(1). It is a direct consequence of Lemma \ref{lem: towerbase} (4). 

(2). Apply $\id\oti\id\oti\vph$ to the both side of (1). 

(3). It is derived by using (2) and Lemma \ref{lem: diagonal} (4). 

(4). 
Recall a system of matrix units 
$\{e_{\opi_{\ovl{i},\ovl{j}}}\}_{i,j\in I_\pi}$ defined in \S2.2. 
Set functionals 
$\om_{\pi_{i,j}}=\Tr_\pi e_{\pi_{j,i}}$, 
$\om_{\opi_{\ovl{k},\ovl{\el}}}=\Tr_\opi e_{\opi_{\ovl{\el},\ovl{k}}}$. 
Apply 
$\id\oti\om_{\pi_{i,j}}\oti\om_{\opi_{\ovl{k},\ovl{\el}}}$ 
to the both sides of (1). 
Then the left hand side is equal to
\begin{align*}
&(\id\oti\om_{\pi_{i,j}}\oti\om_{\opi_{\ovl{k},\ovl{\el}}})
(
d_\pi^2(\al_\pi(e)a_\pi\oti1_\opi)(1\oti\De_{\pi,\opi}(e_\btr))
(a_\pi^*\al_\pi(e)\oti1_\opi)
)
\\
=&\,
d_\pi^2 (\id\oti\om_{\pi_{i,j}})
(
\al_\pi(e)a_\pi 
d_\pi^{-1}(1\oti e_{\pi_{k,\el}})
a_\pi^*\al_\pi(e)
)
\\
=&\,
d_\pi\big{(}\al_\pi(e)a_\pi\big{)}_{i,k}
\big{(}a_\pi^*\al_\pi(e)\big{)}_{\el,j}
\end{align*}
and the right hand side is equal to 
$\de_{k,\el}\big{(}\al_\pi(e)\big{)}_{i,j}$. 

(5). 
The self-adjointness yields 
$f_{\opi_{i,j}}^*=f_{\opi_{j,i}}$. 
By (4), 
$f_{\opi_{i,j}}
=d_\pi \big{(}a_\pi^*\al_\pi(e)a_\pi\big{)}_{i,j}$ 
satisfies the desired equality. 
\end{proof}

In paticular, (5) shows the similarity of projections 
$a^*\al(e)a_K$ and $\De_K(e_\btr)$. 
We presume that the projection $a^*\al(e)a_K$ 
makes a copy of finite dimensional algebra $\lhG \ovl{K}$ 
in $M$. 
This is a reason for giving indices not $f_{\pi_{i,j}}$ 
but $f_{\opi_{i,j}}$. 
We close this subsection with the next useful lemma. 
It shows that the support projection of each tower element is 
invariant by perturbation if $e$ commutes with $u$ and $\al(v)$. 

\begin{lem}\label{lem: invariance}
Let $(\al,u)$ be a cocycle action of $\bhG$ on a von Neumann algebra 
$M$. 
Let $v\in U(M\oti \lhG)$ and $(\wdt{\al},\wdt{u})$ 
the perturbed cocycle action of $(\al,u)$ by $v$. 
Let $a$ and $\wdt{a}$ be the diagonals of 
$u$ and $\wdt{u}$, respectively. 
Let $\mK$ be a finite subset of $\IG$. 
Assume that $e\in\Proj(M)$ satisfies 
$[e\oti1\oti1,u]=0=[e\oti1\oti 1,\al(v)]$ and 
$(e\oti1_\rho)\al_\rho(e)=0$ for all 
$\rho\in\ovl{\mK}\cdot\mK\setm\{\btr\}$. 
Then one has
\[
(\id\oti\vph)(a^*\al(e)a_\pi)
=(\id\oti\vph)(\wdt{a}^*\wdt{\al}(e)\wdt{a}_\pi). 
\]
\end{lem}

\begin{proof}
By Lemma \ref{lem: towerbase} and Lemma \ref{lem: matrixunit}, 
the elements 
$q_\pi=(\id\oti\vph)(a^*\al(e)a_\pi)$ and 
$\wdt{q}_\pi=(\id\oti\vph)(\wdt{a}^*\wdt{\al}(e)\wdt{a}_\pi)$ 
are the projections which are given by
\begin{align*}
&q_\pi
=\inf\{q\in\Proj(M)\mid (q\oti1_\pi)\al_\pi(e)=\al_\pi(e)
,\ [q\oti1\oti1,u]=0\},
\\
&\wdt{q}_\pi
=\inf\{q\in\Proj(M)\mid (q\oti1_\pi)\wdt{\al}_\pi(e)=\wdt{\al}_\pi(e)
,\ [q\oti1\oti1,\wdt{u}]=0\}. 
\end{align*}
Since $q_\pi=d_\pi^2\Ph_\opi(e\oti1_\opi)$, we have 
\begin{align*}
(q_\pi\oti1)v
=&\,
(d_\pi^2\Ph_\opi(e\oti1_\opi)\oti1)v
\\
=&\,
d_\pi^2(\Ph_\opi\oti\id)((e\oti1_\opi\oti1)\al_\opi(v))
\\
=&\,
d_\pi^2(\Ph_\opi\oti\id)(\al_\opi(v)(e\oti1_\opi\oti1))
\\
=&\,
v(d_\pi^2\Ph_\opi(e\oti1_\opi)\oti1)
\\
=&\,
v(q_\pi\oti1). 
\end{align*}
Hence $q_\pi$ commutes with $v$. 
Similarly we can show that $q_\pi$ commutes with $\al(v)$, 
and so does with $\wdt{u}$. 
Then we have $(q_\pi\oti1_\pi)\wdt{\al}_\pi(e)=\wdt{\al}_\pi(e)$. 
It yields $q_\pi\geq \wdt{q}_\pi$. 
We prove $q_\pi\leq \wdt{q}_\pi$. 
Since $e$ commutes with $\wdt{\al}(v)=(v\oti1)\al(v)(v^*\oti1)$, 
$\wdt{q}_\pi$ commutes with $v$ and $\al(v)$, 
and hence so does with $u$, and then 
we have $(\wdt{q}_\pi\oti1_\pi)\al_\pi(e)=\al_\pi(e)$. 
This equality yields $q_\pi\leq \wdt{q}_\pi$. 
\end{proof}

\subsection{Rohlin type theorem}
We present a Rohlin type theorem. 
In order to simplify our arguments, we treat only McDuff factors. 

\begin{thm}\label{thm: jRohlin}
Let $M$ be a McDuff factor and 
$(\al,u)$  an approximately inner and strongly free cocycle action 
of $\bhG$ on $M$. 
Take a unitary $v\in M^\om\oti\lhG$ such that 
\begin{enumerate}[(i)]

\item $\al=\Ad v $ on $M\subs M^\om$, 

\item $(v^*\oti1)(\al^\om\oti\id)(v^*)u(\id\oti\De)(v)=1$,  

\item $\ga=\Ad v^*\circ\al^\om$ is an action preserving $M_\om$,

\item $(\ta^\om\oti\ta_\pi)\circ\ga_\pi=\ta^\om$ for all $\pi\in\IG$. 

\end{enumerate}
Let $\ph$ be a faithful normal state on $M$ 
and set $\ps=\ph\circ\ta^\om$. 
Let $0<\de,\ka<1$ and $F\in \Projf(Z(\lhG))$. 
Take $K\in\Projf(Z(\lhG))$ which is 
$(F,\de)$-invariant and satisfies $K\geq e_\btr$. 
Set $\mF=\supp(F)$ and $\mK=\supp(K)$. 
Assume that the 2-cocycle $u$ is small in the following sense, 
\[
\|u_{\pi,\rho}-1\oti1_\pi\oti1_\rho\|_{\psi\oti\vph\oti\vph}<\ka
\]
for all $\pi\in\mF\cup \mK$ and $\rho\in\ovl{\mK}$. 
Then for any countable set $S\subs M^\om$, 
there exists a projection $E$ 
in $S'\cap M_\om\oti\lhG$ satisfying the following conditions. 
\begin{enumerate}

\item 
$E=E(1\oti K)$. 

\item
(approximate equivariance)
\[
\big{|}\ga_F(E)-(\id\oti_F\De_K)(E)
\big{|}_{\ps\oti\vph\oti\vph}
<
5\de^{1/2}|F|_\vph.
\]

\item
Decompose $E$ as 
\[
E
=
\sum_{\rho\in\mK}\sum_{i,j\in I_\rho}
d_\rho^{-1} f_{\orho_{i,j}}\oti e_{\rho_{i,j}}. 
\]
Then $\{f_{\orho_{i,j}}\}_{i,j\in I_\rho}$ 
is a system of matrix units. 
Moreover, 
they are orthogonal in the following sense. 
For all $\rho\neq\pi\in\mK$, 
$i,j\in I_\rho$ and $k,\el\in I_\pi$,  
\[
f_{\orho_{i,j}} f_{\opi_{k,\el}}=0.
\]

\item (joint property of $U$) 
Let $a$ be the diagonal of $u$. 
Set an operator $U=a^*vE$ and decompose as 
\[
a^*vEv^*a_K
=
\sum_{\rho\in\mK}\sum_{i,j\in I_\rho}
d_\rho^{-1} f_{\orho_{i,j}}^\al \oti e_{\rho_{i,j}}, 
\]
\[
U
=
\sum_{\rho\in\mK}\sum_{i,j\in I_\rho}
d_\rho^{-1} \mu_{\orho_{i,j}}\oti e_{\rho_{i,j}}. 
\]
Then we have 
\[
\mu_{\orho_{i,j}}^*\mu_{\opi_{k,\el}}
=\de_{\rho,\pi}\de_{i,k} f_{\orho_{j,\el}},
\quad
\mu_{\orho_{i,j}}\mu_{\opi_{k,\el}}^*
=\de_{\rho,\pi}\de_{j,\el} f_{\orho_{i,k}}^\al,
\]
for all $\rho,\pi\in \mK$, $i,j\in I_\rho$ and $k,\el\in I_\pi$. 
In particular $U^*U=E$. 

\item 
$Ev^* aa^*vE=E$. 

\item 
For each $\rho\in\mK$, 
the projection $(\id\oti\vph_\rho)(E)$ is 
in $S'\cap M_\om$ 
and satisfies 
\[
(\id\oti\vph_\rho)(E)=(\id\oti\vph_\rho)(a^*vEv^*a). 
\]

\item (partition of unity)
\[
(\id\oti\vph)(E)=1=(\id\oti\vph)(a^*vE v^*a). 
\]

\item
(Shapiro lemma) 
Set $\mu=(\id\oti\vph)(U)$ and then $\mu$ is a unitary. 
If the weights $\ps\oti\vph$ and $\ps\oti\vph\oti\vph$ are 
invariant for $\Ad v$ and $\Ad(\id\oti\De)(v)$, respectively, 
then the following inequality holds. 
\[
|v_F \ga_F^\om(\mu)-\mu\oti F|_{\psi\oti\vph}
<
(9\de^{1/4}+3\ka^{1/2})|F|_\vph. 
\] 
\end{enumerate}
\end{thm}

In this situation, we call a projection $E$ and a unitary 
$\mu$ a \textit{Rohlin projection} 
and a \textit{Shapiro unitary}, respectively. 

\begin{lem}\label{lem: a}
Let $\pi\in\mK$. 
If $\|u_{\pi,\opi}-1\oti1_\pi\oti1_\opi\|_{\ps\oti\vph\oti\vph}<\ka$, 
then we have 
$\|a_\pi-1\oti1_\pi\|_{\ps\oti\vph}<\ka$. 
\end{lem}
\begin{proof} 
Set a functional $\th_\pi=T_{\pi,\opi}^*\cdot T_{\pi,\opi}$ on 
$B(H_\pi\oti H_\opi)$. 
Then we have 
\begin{align*}
(\ps\oti\vph_\pi)(|a_\pi-1\oti1_\pi|^2)
=&\,
d_\pi^2 (\ps\oti\th_\pi)(|a_\pi-1\oti1_\pi|^2\oti1_\opi)
\\
=&\,
d_\pi^2 (\ps\oti\th_\pi)(|u_{\pi,\opi}-1\oti1_\pi\oti1_\opi|^2)
\\
\leq&\,
(\ps\oti\vph\oti\vph)(|u_{\pi,\opi}-1\oti1_\pi\oti1_\opi|^2)
\\
<&\,
\ka^2, 
\end{align*}
where we have used $d_\pi^2\th_\pi\leq \vph_\pi\oti\vph_\opi$. 
\end{proof}

Define a set $\meJ$ which consists of a projection 
$E$ in $M_\om\oti\lhG$ 
satisfying the conditions (1), (3), (4), (5), (6) 
in Theorem \ref{thm: jRohlin}, and in addition, the following ones, 
\begin{enumerate}[(a)]

\item 
$(\ta^\om\oti\id)(E)=(\ta^\om\oti\id)(vEv^*)\in \C K$,

\item
$(\id\oti\vph_\rho)(E)=(\id\oti\vph_\rho)(a^*vEva)$ for all $\rho\in\mK$.

\end{enumerate}
Define functions $a$ and $b$ from $\meJ$ to $\R_+$ by 
\begin{align*}
a_E=&\,
|F|_\vph^{-1}
\big{|}\ga_F(E)-(\id\oti_F\De_K)(E)\big{|}_{\ps\oti\vph\oti\vph},
\\
b_E=&\,
|E|_{\ps\oti\vph}.
\end{align*}

\begin{lem}\label{lem: E'}
Let $E$ be an element of $\meJ$. 
Assume $b_E<1-\de^{1/2}$. 
Then there exists $E'\in\meJ$ satisfying the following inequalities. 
\begin{enumerate}

\item 
$a_{E'}-a_{E}\leq 3 \de^{1/2}(b_{E'}-b_E)$, 

\item
$0<(\de^{1/2}/2)|E'-E|_{\ps\oti\vph}\leq b_{E'}-b_E$.

\end{enumerate}
\end{lem}

\begin{proof}
We may assume that the entries of $E$, $v$ and $u$ are in $S$ 
and $S$ is $\al^\om$-invariant. 
Take a projection $e$ from $S'\cap M_\om$ such that 
$(e\oti1_\rho)\ga_\rho(e)=0$ for all 
$\rho\in\ovl{\mK}\cdot\mK\setm\{\btr\}$ 
by Lemma \ref{lem: partition}. 
Since $e$ commutes with $v$, $(e\oti1_\rho)\al_\rho^\om(e)=0$ also holds. 
By Lemma \ref{lem: invariance}, we have 
\[
(\id\oti\vph)(\ga_\rho(e))=(\id\oti\vph)(a_\rho^*\al_\rho^\om(e)a_\rho)
\in S'\cap M_\om. 
\]
Let $N$ be a von Neumann subalgebra in $M^\om$ which is generated 
by $M$ and the entries of $\{\al_\rho^\om(e)\}_{\rho\in\IG}$. 
Applying the Fast Reindexation Trick for $N$ and $S$, 
we have a map 
$\Ps\in\Mor(\wdt{N}, M^\om\oti \lhG)$ 
as in Lemma \ref{lem: fast}. 
Set $f=\Ps(e)$ and then $f\in S'\cap M_\om$. 
Since 
$(f\oti1_\rho)\al_\rho^\om(f)
=(\Ps\oti\id)((e\oti1_\rho)\al_\rho^\om(e))
=0$ for $\rho\in\ovl{\mK}\cdot\mK\setm\{\btr\}$, 
the equality $(f\oti1_\rho)\ga_\rho(f)=0$ holds. 
Then by Lemma \ref{lem: invariance}, we have 
\[
(\id\oti\vph)(\ga_\rho(f))
=
(\id\oti\vph)(a_\rho^* \al_\rho^\om(f)a_\rho)
\\
=
\Ps\big{(}
(\id\oti\vph)(a_\rho^* \al_\rho^\om(e)a_\rho)\big{)}. 
\]
This shows the following splitting property of $\ta^\om$
for $x\in S\oti \lhG$, 
\begin{align*}
&(\ta^\om\oti\id)
\big{(}x((\id\oti\vph)(\ga_\rho(f))\oti1)\big{)}
\\
=&\,
(\ta^\om\oti\id)
\big{(}x\big{(}\Ps\big{(}
(\id\oti\vph)(a_\rho^* \al_\rho^\om(e)a_\rho)\big{)}\oti1\big{)}\big{)}
\\
=&\,
(\ta^\om\oti\id)(x)
\cdot
(\ta^\om\oti\id)
\big{(}
\Ps\big{(}(\id\oti\vph)(a_\rho^* \al_\rho^\om(e)a_\rho)\big{)}\oti1
\big{)}
\\
=&\,
(\ta^\om\oti\id)(x)
\cdot
\big{(}
\ta^\om((\id\oti\vph)(\ga_\rho(f)))\oti1\big{)}.
\end{align*}
Set a projection $f'=(\id\oti\vph)(\ga_K(f))$ in $S'\cap M_\om$. 
The equality $(\ta^\om\oti\ta_\rho)\circ\ga_\rho=\ta^\om$ 
yields 
$\ta^\om((\id\oti\vph)(\ga_\rho(f)))=\ta_\om(f)d_\rho^2$ 
for $\rho\in\mK$, 
in particular $|f'|_\ps=|f'|_{\ta_\om}=\ta_\om(f)|K|_\vph$. 
Then set a projection $E'\in M_\om\oti \lhG$ by 
\[
E'=E(f'^\per\oti1)+\ga_K(f). 
\]
We verify that $E'$ is a desired projection. 
At first, we will show $E'\in \meJ$. 
The condition (1) in Theorem \ref{thm: jRohlin} is trivial. 
The conditions (3), (4) and (5) hold by Lemma \ref{lem: matrixunit}. 
The condition (6) follows from Lemma \ref{lem: invariance} 
on $\ga_K(f)$. 
We verify the remained conditions (a) and (b). 
We claim $(\ta^\om\oti\id)(\ga_K(f))=\ta_\om(f)K$. 
Indeed, by applying Proposition \ref{prop: free-left} 
to the action $\ga$ on $M_\om$, 
we have $(\ta_\om\oti\id)\circ\ga_\rho=\ta_\om$. 
Then the condition (a) is verified as 
\begin{align*}
(\ta^\om\oti\id)(E')
=&\,
(\ta^\om\oti\id)(E(f'^\per\oti1)+\ga_K(f))
\\
=&\,
(\ta_\om\oti\id)(E)(\ta_\om(f'^\per)\oti1)
+\ta_\om(f)K
\end{align*}
and 
\begin{align*}
(\ta^\om\oti\id)(vE'v^*)
=&\,
(\ta^\om\oti\id)(vEv^*(f'^\per\oti1)+\al_K^\om(f))
\\
=&\,
(\ta^\om\oti\id)(vEv^*)
(\ta_\om(f'^\per)\oti1)
+
\ta_\om(f)K
\\
=&\,
(\ta_\om\oti\id)(E)(\ta_\om(f'^\per)\oti1)
+\ta_\om(f)K. 
\end{align*}

Next we verify the condition (b). 
Since $(\id\oti\vph)(\ga_\rho(f))=(\id\oti\vph)(a^*\al^\om(e)a_\rho)$, 
we have 
\begin{align*}
(\id\oti\vph)(E_\rho')
=&\,
(\id\oti\vph)(E_\rho)f'^\per+(\id\oti\vph)(\ga_\rho(f))
\\
=&\,
(\id\oti\vph)(a^*vEv^*a_\rho)f'^\per
+(\id\oti\vph)(a^*\al^\om(e)a_\rho)
\\
=&\,
(\id\oti\vph)(a^*vE'v^*a_\rho). 
\end{align*}
Hence $E'\in \meJ$. Now we estimate $a_{E'}$ and $b_{E'}$. 
First we have 
\begin{align*}
|E'-E|_{\ps\oti\vph}
=&\,
|-E(f'\oti1)+\ga_K(f)|_{\ps\oti\vph}
\\
\leq&\,
|E(f'\oti1)|_{\ps\oti\vph}
+
|\ga_K(f)|_{\ps\oti\vph}
\\
=&\,
|(\id\oti\vph)(E)f'|_\ps
+
|f'|_\ps
\\
\leq&\,
2|f'|_\ps. 
\end{align*}
Since 
\begin{align*}
b_{E'}
=&\,
|E(f'^\per\oti1)+\ga_K(f)|_{\ps\oti\vph}
\\
=&\,
|E(f'^\per\oti1)|_{\ps\oti\vph}+|\ga_K(f)|_{\ps\oti\vph}
\\
=&\,
\ps(f'^\per)|E|_{\ps\oti\vph}+|f'|_\ps
\\
=&\,
\ps(f'^\per)b_E+|f'|_\ps,
\end{align*}
we have the condition (2) in Lemma \ref{lem: E'} as follows. 
\begin{align*}
b_{E'}-b_{E}
=&\,
\ps(f')(1-b_E)
\\
>&\,
\de^{1/2}|f'|_\ps
\\
\geq&\,
\frac{\de^{1/2}}{2}|E'-E|_{\ps\oti\vph}. 
\end{align*}
Secondly we verify the condition (1) in Lemma \ref{lem: E'}.  
By direct calculation, we have 
\begin{align}
&(\ga_F\oti\id)(E')-(\id\oti_F\De_K)(E')
\notag\\
=&\,
(\ga_F\oti\id)(E(f'^\per\oti1))-(\id\oti_F\De_K)(E(f'^\per\oti1))
\notag\\
&\quad+
(\ga_F\oti\id)(\ga_K(f))-(\id\oti_F\De_K)(\ga_K(f))
\notag\\
=&\,
(\ga_F\oti\id)(E)(\ga_F(f'^\per)\oti1-f'^\per\oti1\oti1)
\label{eq: a1}\\
&\quad
+(\ga_F\oti\id)(E)(f'^\per\oti1\oti1)-(\id\oti_F\De_K)(E(f'^\per\oti1))
\label{eq: a2}\\
&\quad
+(\id\oti_F\De_K)(\ga_{K^\per}(f))
\label{eq: a3}. 
\end{align}
We estimate the trace norms of the above three terms. 

On (\ref{eq: a1}), 
we know $[(\ga_F\oti\id)(E),\ga_F(f'^\per)\oti1-f'^\per\oti F\oti1]=0$, 
and then
\begin{align*}
|(\ref{eq: a1})|_{\ps\oti\vph\oti\vph}
=
&\,\big{|}
(\ga_F\oti\id)(E)(\ga_F(f'^\per)\oti1-f'^\per\oti F\oti1)
\big{|}_{\ps\oti\vph\oti\vph}
\\
=&\,
(\ps\oti\vph\oti\vph)
\big{(}
(\ga_F\oti\id)(E)
\big{(}\big{|}\ga_F(f'^\per)-f'^\per\oti F \big{|}\oti1\big{)}
\big{)}
\\
=&\,
(\ps\oti\vph)
\big{(}
\ga_F((\id\oti\vph)(E))
\big{|}\ga_F(f'^\per)-f'^\per\oti F\big{|}
\big{)}
\\
\leq&\,
(\ps\oti\vph)
\big{(}
\big{|}\ga_F(f'^\per)-f'^\per\oti F\big{|}
\big{)}
\\
=&\,
(\ps\oti\vph)
\big{(}
\big{|}\ga_F(f')-f'\oti F\big{|}
\big{)}
\\
=&\,
\big{|}
(\id\oti\id\oti\vph)
\big{(}
(\ga_F\oti\id)(\ga_K(f))-(\id\oti_F\De)(\ga_K(f))
\big{)}
\big{|}_{\ps\oti\vph}
\\
=&\,
\big{|}
(\id\oti\id\oti\vph)
\big{(}
-(\id\oti_F\De_{K^\per})(\ga_K(f))+(\id\oti_F\De_K)(\ga_{K^\per}(f))
\big{)}
\big{|}_{\ps\oti\vph}
\\
\leq&\,
\big{|}
(\id\oti_F\De_{K^\per})(\ga_K(f))
\big{|}_{\ps\oti\vph\oti\vph}
+
\big{|}
(\id\oti_F\De_K)(\ga_{K^\per}(f))
\big{|}_{\ps\oti\vph\oti\vph}
\end{align*}
holds. 
On the last terms, we have 
\begin{align*}
\big{|}
(\id\oti_F\De_{K^\per})(\ga_K(f))
\big{|}_{\ps\oti\vph\oti\vph}
=&\,
(\ps\oti\vph\oti\vph)\big{(}
(\id\oti_F\De_{K^\per})(\ga_K(f))
\big{)}
\\
=&\,
(\vph\oti\vph)\big{(}
(\ta^\om\oti_F\De_{K^\per})(\ga_K(f))
\big{)}
\\
=&\,
\ta_\om(f)(\vph\oti\vph)((F\oti K^\per)\De(K))
\\
<&\,
|f|_\ps \de|F|_\vph|K|_\vph
\\
=&\,
\de|F|_\vph|f'|_\ps,
\end{align*}
where we have used the $(F,\de)$-invariance of $K$. 
Similarly we get 
\[
\big{|}
(\id\oti_F\De_K)(\ga_{K^\per}(f))
\big{|}_{\ps\oti\vph\oti\vph}
<\de|F|_\vph|f'|_\ps. 
\]
Hence we obtain 
\[
|(\ref{eq: a1})|_{\ps\oti\vph\oti\vph}<2\de|F|_\vph|f'|_\ps.
\]

On (\ref{eq: a2}), we have 
\begin{align*}
|(\ref{eq: a2})|_{\ps\oti\vph\oti\vph}
=&\,\big{|}
\big{(}(\ga_F\oti\id)(E)-(\id\oti_F\De_K)(E)\big{)}
(f'^\per\oti1\oti1)
\big{|}_{\ps\oti\vph\oti\vph}
\\
\leq&\,
\big{|}
(\ga_F\oti\id)(E)-(\id\oti_F\De_K)(E)
\big{|}_{\ps\oti\vph\oti\vph}
\\
=&\,
|F|_\vph a_E. 
\end{align*}

On (\ref{eq: a3}), we have 
\begin{align*}
|(\ref{eq: a3})|_{\ps\oti\vph\oti\vph}
=&\,
\big{|}
(\id\oti_F\De_K)(\ga_{K^\per}(f))
\big{|}_{\ps\oti\vph\oti\vph}
\\
=&\,
\ps(f)\big{|}
{}_F\De_K(K^\per)
\big{|}_{\vph\oti\vph}
\\
<&\,
\de|F|_\vph|f'|_\ps. 
\end{align*}
Summarizing these calculations, we have 
\begin{align*}
|F|_\vph a_{E'}
\leq&\,
|(\ref{eq: a1})|_{\ps\oti\vph\oti\vph}
+
|(\ref{eq: a2})|_{\ps\oti\vph\oti\vph}
+
|(\ref{eq: a3})|_{\ps\oti\vph\oti\vph}
\\
<&\,
2\de|F|_\vph|f'|_\ps+|F|_\vph a_E+ \de|F|_\vph|f'|_\ps
\\
=&\, 
|F|_\vph a_E+ 3\de|F|_\vph|f'|_\ps
\\
\leq&\,
|F|_\vph a_E+3\de^{1/2} |F|_\vph(b_{E'}-b_E). 
\end{align*}
Hence the condition (1) in Lemma \ref{lem: E'} holds. 
\end{proof}

\textit{Proof of Theorem \ref{thm: jRohlin}}. 

Consider a subset $\meS\subs \meJ$ whose element $E$ satisfies 
$a_E\leq 3 \de^{1/2}b_E$. 
We order $\meS$ by $E\prec E'$ if $E=E'$ or the inequalities (1) and (2) 
in Lemma \ref{lem: E'} hold. 
Since $\meS$ contains 0, $\meS$ is nonempty. 
The order of $\meS$ is inductive as is shown below. 
By Lemma \ref{lem: E'} (2), 
the map $b$ is an order preserving isomorphism 
on a totally ordered subset $\meL\subs \meS$ onto 
a subset in $[0,1]$. 
Hence $\meL$ is cofinal. 
Then again with (2), the cofinal subsequence of $\meL$ strongly converges 
to a projection. 
We can easily observe that $\meS$ is strongly closed. 
Hence the supremum of $\meL$ exists in $\meS$. 
By Zorn's lemma, there exists a maximal element $\ovl{E}$ in $\meS$. 
Assume $b_{\bar E}<1-\de^{1/2}$ and then by Lemma \ref{lem: E'} 
we can take an element $E'\in\meJ$ which satisfies the conditions in 
the lemma for $\ovl{E}$. 
It is easy to see that $E'\in \meS$, $\bar{E}\prec E'$ and $\bar{E}\neq E'$, 
but this is a contradiction. 
Hence we have $b_{\bar E}\geq1-\de^{1/2}$. 
Set a projection 
$p=1-(\id\oti\vph)(\ovl{E})=1-(\id\oti\vph)(a^*v\ovl{E}v^*a)$ in 
$S'\cap M_\om$ and then we have $\ta_\om(p)\leq\de^{1/2}$. 
Then set a projection $E=\ovl{E}+p\oti e_\btr$. 
Since $K\geq e_\btr$, $E=E(1\oti K)$. 
We verify all the conditions in Theorem \ref{thm: jRohlin}. 
The conditions (3), (4), (5), (6), (7) are 
immediately verified. 

On the condition (2), we estimate as follows, 
\begin{align*}
&\big{|}\ga_F(E)-(\id\oti_{F}\De_K)(E)
\big{|}_{\ps\oti\vph\oti\vph}
\\
\leq&\,
\big{|}\ga_F(\ovl{E})-(\id\oti_{F}\De_K)(\ovl{E})
\big{|}_{\ps\oti\vph\oti\vph}
+
\big{|}\ga_F(p)\oti e_\btr- p\oti_{F}\De_K(e_\btr)
\big{|}_{\ps\oti\vph\oti\vph}
\\
\leq&\,
3\de^{1/2}|F|_\vph b_{\ovl{E}}
+
2|p|_\ps|F|_\vph
\\
\leq&\,
5\de^{1/2}|F|_\vph. 
\end{align*}

Finally we show the condition (8). 
By the conditions (4) and (7), the element 
$\mu=(\id\oti\vph)(U)$ is a unitary. 
We claim the following inequalities. 
\begin{align}
&\|(a^*-1\oti K)vE\|_{\ps\oti\vph}<\ka,
\label{eq: jroh10-1}\\
&\|\ga_F\big{(}(a^*-1\oti K)vE\big{)}\|_{\ps\oti\vph\oti\vph}
<\ka\|F\|_\vph,
\label{eq: jroh10-2}\\
&\|(u_{F,K}-1\oti F\oti K)(\id\oti\De)(vE)\|_{\ps\oti\vph\oti\vph}
<(4\de^{1/2}+\ka^2)^{1/2}\|F\|_\vph. 
\label{eq: jroh10-3}
\end{align}
The inequality (\ref{eq: jroh10-2}) 
is an immediate consequence from (\ref{eq: jroh10-1}) 
since $\ps\oti\ta_\pi\circ\ga=\ps$ for all $\pi\in\IG$. 
The inequality (\ref{eq: jroh10-1}) is proved as follows, 
\begin{align*}
&\|(a^*-1\oti K)vE\|_{\ps\oti\vph}^2
\\
=&\,
(\ps\oti\vph)(Ev^*|a^*-1\oti K|^2vE)
\\
=&\,
(\ps\oti\vph)(Ev^*aa^*vE+E)-2\Re (\ps\oti\vph)(Ev^*avE)
\\
=&\,
2|E|_{\ps\oti\vph}-2\Re (\ps\oti\vph)(v^*avE)
\\
=&\,
2|E|_{\ps\oti\vph}-2\Re (\ps\oti\vph)(avEv^*)
\\
=&\,
2|E|_{\ps\oti\vph}-2\Re (\ph\circ\ta^\om\oti\vph)(avEv^*)
\\
=&\,
2|E|_{\ps\oti\vph}
-2\Re (\ph\oti\vph)\big{(}
a(\ta^\om\oti\id)(v\ovl{E}v^*+p\oti e_\btr)
\big{)}
\\
=&\,
2|E|_{\ps\oti\vph}
-2\Re (\ph\oti\vph)\big{(}
a(b_{\ovl{E}}|K|_\vph^{-1}(1\oti K)+\ta_\om(p)(1\oti e_\btr)
\big{)}
\\
=&\,
2|E|_{\ps\oti\vph}
-2\Re b_{\ovl{E}}|K|_\vph^{-1}(\ph\oti\vph)(a_K)
-2\ta_\om(p)
\\
=&\,
b_{\ovl{E}}|K|_\vph^{-1}(2|K|_\vph-2\Re (\ph\oti\vph)(a_K))
\\
\leq&\,
|K|_\vph^{-1} \|a_K-1\oti K\|_{\ph\oti\vph}^2
\\
<&\,
\ka^2,
\end{align*}
where we have used Lemma \ref{lem: a}. 
The inequality (\ref{eq: jroh10-3}) is obtained as follows, 
\begin{align*}
&\|(u_{F,K}-1\oti F\oti K)(\id\oti\De)(vE)\|_{\ps\oti\vph\oti\vph}^2
\\
=&\,
(\ps\oti\vph\oti\vph)
((\id\oti\De)(Ev^*)|u_{F,K}-1\oti F\oti K|^2(\id\oti\De)(vE))
\\
=&\,
(\ps\oti\vph\oti\vph)
(|u_{F,K}-1\oti F\oti K|^2(\id\oti\De)(vEv^*))
\\
=&\,
(\ph\circ\ta^\om\oti\vph\oti\vph)
(|u_{F,K}-1\oti F\oti K|^2(\id\oti\De)(vEv^*))
\\
=&\,
(\ph\oti\vph\oti\vph)
(|u_{F,K}-1\oti F\oti K|^2(\ta^\om\oti\De)(vEv^*))
\\
=&\,
(\ph\oti\vph\oti\vph)
(|u_{F,K}-1\oti F\oti K|^2
\\
&\hspace{3cm}
\cdot(b_{\ovl{E}}|K|_\vph^{-1}(1\oti \De(K))+\ta_\om(p)(1\oti \De(e_\btr))))
\\
\leq&\,
b_{\ovl{E}}|K|_\vph^{-1}
\|u_{F,K}-1\oti F\oti K\|_{\ps\oti\vph\oti\vph}^2
\\
&\hspace{3cm}+
\ta_\om(p)
\|(u_{F,K}-1\oti F\oti K)(1\oti\De(e_\btr))\|_{\ps\oti\vph\oti\vph}^2
\\
\leq&\,
b_{\ovl{E}}|K|_\vph^{-1} \ka^2|\mF||\mK|
+
4\ta_\om(p)\|1\oti_F\De_K(e_\btr)\|_{\ps\oti\vph\oti\vph}^2
\\
\leq&\,
\ka^2b_{\ovl{E}}|K|_\vph^{-1} |F|_\vph|K|_\vph
+4\ta_\om(p)|F|_\vph
\\
<&\,
\ka^2|F|_\vph+4\de^{1/2}|F|_\vph. 
\end{align*}

Let 
\[
v_F\ga_F(\mu)-\mu\oti F=w|v_F\ga_F(\mu)-\mu\oti F|
\] be the polar decomposition with the partial isometry 
$w\in M^\om\oti \lhG$. 
Then we have 
\begin{align}
&|v_F\ga_F(\mu)-\mu\oti F|_{\ps\oti\vph}
\notag\\
=&\,
(\ps\oti\vph)(w^* (v_F\ga_F(\mu)-\mu\oti F))
\notag\\
=&\,
(\ps\oti\vph)
\big{(}w^*
(\id\oti\id\oti\vph)((v_F\oti K)\ga_F(a^* v E)-(\id\oti_F\De)(a^*vE))
\big{)}
\notag\\
=&\,
(\ps\oti\vph\oti\vph)
\big{(}
(w^*\oti K)(v_F\oti K)\ga_F((a^*-1\oti K) v E)
\label{eq: jroh-shap1}\\
&\quad
+
(w^*\oti K)\big{(}(v_F\oti K)\ga_F(v)(\ga_F(E)-(\id\oti_F\De_K)(E))\big{)}
\label{eq: jroh-shap2}\\
&\quad
+
(w^*\oti K)\big{(}\big{(}(v_F\oti K)\ga_F(v)
-(\id\oti_F\De_K)(v)\big{)}
(\id\oti_F\De_K)(E)\big{)}
\label{eq: jroh-shap3}\\
&\quad
+
(w^*\oti 1)\big{(}(\id\oti_F\De_K)(vE)-(\id\oti_F\De)(vE)\big{)}
\label{eq: jroh-shap4}\\
&\quad
+
(w^*\oti 1)(\id\oti_F\De)((1\oti K-a^*)vE)
\big{)}. 
\label{eq: jroh-shap5}
\end{align}
On (\ref{eq: jroh-shap1}), using (\ref{eq: jroh10-2}), we have 
\begin{align*}
|(\ref{eq: jroh-shap1})|
=&\,
\big{|}(\ps\oti\vph\oti\vph)
\big{(}
(w^*\oti K)(v_F\oti K)\ga_F((a^*-1\oti K) v E)
\big{)}
\big{|}
\\
=&\,
\big{|}
(\ps\oti\vph\oti\vph)
\big{(}
\ga_F(E)(w^*\oti K)(v_F\oti K)\ga_F((a^*-1\oti K) vE )
\big{)}
\big{|}
\\
\leq&\,
\big{\|}(v_F^*\oti K)(w\oti K)\ga_F(E)\big{\|}_{\ps\oti\vph\oti\vph}
\big{\|}\ga_F((a^*-1\oti K) vE )\big{\|}_{\ps\oti\vph\oti\vph}
\\
\leq &\,
\|F\|_\vph \big{\|}\ga_F((a^*-1\oti K) v E)\big{\|}_{\ps\oti\vph\oti\vph}
\\
<&\,\ka |F|_\vph. 
\end{align*}
On (\ref{eq: jroh-shap2}), 
the term $\ga_F(E)-(\id\oti_F\De_K)(E)$ is in the centralizer of the 
weight $\ps\oti\vph\oti\vph$. 
Using approximate equivalence of $E$, we have 
\begin{align*}
|(\ref{eq: jroh-shap2})|
=&\,
\big{|}(\ps\oti\vph\oti\vph)
\big{(}
(w^*\oti K)\big{(}(v_F\oti K)\ga_F(v)(\ga_F(E)-(\id\oti_F\De_K)(E))\big{)}
\big{)}
\big{|}
\\
\leq&\,
\|(w^*\oti K)(v_F\oti K)\ga_F(v)\|
\big{|}\ga_F(E)-(\id\oti_F\De_K)(E)\big{|}_{\ps\oti\vph\oti\vph}
\\
<&\,
5\de^{1/2}|F|_\vph. 
\end{align*}
On (\ref{eq: jroh-shap3}), 
using $u=(v\oti1)\ga(v)(\id\oti\De)(v^*)$ and 
(\ref{eq: jroh10-3}), we have 
\begin{align*}
&|(\ref{eq: jroh-shap3})|
\\
=&\,
\big{|}(\ps\oti\vph\oti\vph)
\big{(}
(w^*\oti K)\big{(}(v_F\oti K)\ga_F(v)
-(\id\oti_F\De_K)(v)\big{)}
(\id\oti_F\De_K)(E)
\big{)}
\big{|}
\\
=&\,
\big{|}(\ps\oti\vph\oti\vph)
\big{(}(\id\oti_F\De_K)(E)
(w^*\oti K)(u_{F,K}-1\oti F\oti K)
(\id\oti_F\De_K)(vE)
\big{)}
\big{|}
\\
\leq&\,
\|(w\oti K)(\id\oti_F\De_K)(E)\|_{\ps\oti\vph\oti\vph}
\|(u_{F,K}-1\oti F\oti K)(\id\oti_F\De_K)(vE)\|_{\ps\oti\vph\oti\vph}
\\
<&\,
\|F\|_\vph(4\de^{1/2}+\ka)^{1/2}\|F\|_\vph
\\
=&\,
(4\de^{1/2}+\ka)^{1/2}|F|_\vph. 
\end{align*}
On (\ref{eq: jroh-shap4}), we have 
\begin{align*}
|(\ref{eq: jroh-shap4})|
=
&\,
\big{|}(\ps\oti\vph\oti\vph)
\big{(}
(w^*\oti 1)\big{(}(\id\oti_F\De_K)(vE)-(\id\oti_F\De)(vE)\big{)}
\big{)}
\big{|}
\\
=&\,
\big{|}(\ps\oti\vph\oti\vph)
\big{(}
-(w^*\oti 1)(\id\oti_F\De_{K^\per})(vE)
\big{)}
\big{|}
\\
=&\,
\big{|}
(\ps\oti\vph\oti\vph)
\big{(}
-(\id\oti_F\De_{K^\per})(E)(w^*\oti 1)(\id\oti_F\De_{K^\per})(vE)
\big{)}
\big{|}
\\
\leq&\,
\|(w\oti 1)(\id\oti_F\De_{K^\per})(E)\|_{\ps\oti\vph\oti\vph}
\|(\id\oti_F\De_{K^\per})(vE)\|_{\ps\oti\vph\oti\vph}
\\
\leq&\,
\|F\|_\vph\|(\id\oti_F\De_{K^\per})(E)\|_{\ps\oti\vph\oti\vph}
\\
=&\,
\|F\|_\vph
(\ps\oti\vph\oti\vph)
\big{(}(\id\oti_F\De_{K^\per})(\ovl{E}+p\oti e_\btr)\big{)}^{1/2}
\\
=&\,
\|F\|_\vph
(\ph\circ\ta^\om\oti\vph\oti\vph)
\big{(}(\id\oti_F\De_{K^\per})(\ovl{E}+p\oti e_\btr)\big{)}^{1/2}
\\
=&\,
\|F\|_\vph
(\ph\oti\vph\oti\vph)
\big{(}(\id\oti_F\De_{K^\per})(b_{\ovl{E}}|K|_\vph^{-1}1\oti K
+\ta_\om(p)\oti e_\btr)\big{)}^{1/2}
\\
\leq&\,
\|F\|_\vph 
\big{(}b_{\ovl{E}}|K|_\vph^{-1}
(\vph\oti\vph)((F\oti K^\per)\De(K))
+\ta_\om(p)|F|_\vph
\big{)}^{1/2}
\\
<&\,
\|F\|_\vph 
\big{(}b_{\ovl{E}}|K|_\vph^{-1}\de|F|_\vph|K|_\vph
+\ta_\om(p)|F|_\vph
\big{)}^{1/2}
\\
\leq&\,
|F|_\vph(\de+\de^{1/2})^{1/2}
\\
<&\,
2^{1/2}\de^{1/4}|F|_\vph.
\end{align*}
Finally on (\ref{eq: jroh-shap5}), 
using (\ref{eq: jroh10-1}), we have 
\begin{align*}
|(\ref{eq: jroh-shap5})|
=&\,
\big{|}(\ps\oti\vph\oti\vph)
\big{(}
(w^*\oti 1)(\id\oti_F\De)((1\oti K-a_K^*)vE)
\big{)}
\big{|}
\\
=&\,
\big{|}
(\ps\oti\vph\oti\vph)
\big{(}
(\id\oti_F\De)(E)
(w^*\oti 1)(\id\oti_F\De)((1\oti K-a_K^*)vE)
\big{)}
\big{|}
\\
\leq&\,
\|(w\oti 1)(\id\oti_F\De)(E)\|_{\ps\oti\vph\oti\vph}
\|(\id\oti_F\De)((1\oti K-a_K^*)vE)\|_{\ps\oti\vph\oti\vph}
\\
\leq&\,
\|F\|_\vph 
\|F\|_\vph\|(1\oti K-a_K^*)vE\|_{\ps\oti\vph}
\\
<&\,
\ka|F|_\vph. 
\end{align*}
Therefore we obtain 
\begin{align*}
|v_F\ga_F(\mu)-\mu\oti F|_{\ps\oti\vph}
\leq&\,
|(\ref{eq: jroh-shap1})|
+
|(\ref{eq: jroh-shap2})|
+
|(\ref{eq: jroh-shap3})|
+
|(\ref{eq: jroh-shap4})|
+
|(\ref{eq: jroh-shap5})|
\\
<&\,
\big{(}
\ka+5\de^{1/2}+(4\de^{1/2}+\ka)^{1/2}+2^{1/2}\de^{1/4}+\ka
\big{)}
|F|_\vph
\\
\leq&\,
(9\de^{1/4}+3\ka^{1/2})|F|_\vph. 
\end{align*}
\begin{flushright}
$\Box$
\end{flushright}

\section{Cohomology vanishing II}

\subsection{2-cohomology vanishing in McDuff factors of type \II}
In Lemma \ref{lem: vanish-ultra}, 
we have proved the 2-cohomology vanishing result 
in an ultraproduct von Neumann algebra. 
This result ensures the existence of a unitary $v$
which is able to perturb a 2-cocycle $u$ to 
a much smaller 2-cocycle $\wdt{u}$, 
but the problem is how small $v$ is. 
In the construction of $v$ in the proof of Lemma \ref{lem: perturb}, 
we see that $v$ is not small even when $u$ is small. 
It is, however, an approximate 1-cocycle and then 
the Rohlin type theorem enables us to perturb $v$ 
to a small new unitary $v'$ by a 1-coboundary constructed from 
a Shapiro unitary. 
Since the perturbation of $v$ does not change 2-cocycle 
$\wdt{u}$ essentially, 
we can perturb $u$ to make it much smaller by a small unitary $v'$. 
Successive perturbations yield a vanishing result of 2-cocycles 
in the original von Neumann algebra. 
We mention that this strategy has been seen in the context of \cite{Oc2}. 

From now on, we assume that $M$ is a McDuff factor of type \II\ 
with the tracial state $\ta$. 
Then the technical assumption in 
Theorem \ref{thm: jRohlin} (8) automatically stands up for $\ph=\ta$. 
The trace $\ta\circ\ta^\om$ on $M^\om$ is also denoted by $\ta$. 
We choose a suitable net of F{\o}lner sets as follows. 
If $\IG$ is finite, we set $F_0=K_0=1$ and $\de_0=0$. 
When $\IG$ is infinite, 
for each $n\geq0$ we will take 
finitely supported central projections $F_n$, $K_n$ in $\lhG$ 
and $\de_n>0$ inductively such that 
\begin{enumerate}[(1)]

\item $\{F_n\}_{n=0}^\infty$ and 
$\{K_n\}_{n=0}^\infty$ are increasing and strongly converge to $1$, 

\item $\{\de_n\}_{n=0}^\infty$ is decreasing and 
$\displaystyle \lim_{n\to\infty}\de_n=0$,

\item $K_n$ is $(F_n,\de_n)$-invariant and $K_n\geq e_\btr$,

\item $(9\de_n^{1/4}+3\de_{n+1}^{1/2})|F_n|_\vph<(1/2)\de_{n-1}$ 
$(n\geq 1)$. 

\end{enumerate} 

Fix a sequence of ascending finitely supported central projections 
$\{S_n\}_{n=0}^\infty$ in $\lhG$ with $S_0=e_\btr$ 
and
$\cup_{n=0}^\infty \supp (S_n)=\IG$. 
Let $F_0=e_\btr=K_0$ and $\de_0=1$. 
First take $F_1$ and $\de_1>0$ such that 
$F_1\geq F_0\vee K_0\vee \ovl{K_0}\vee S_1$ and 
$9\de_1^{1/4}|F_1|_\vph<(1/2)\de_0$. 
Then take an $(F_1,\de_1)$-invariant 
finitely supported central projection $K_1$ with $K_1\geq e_\btr$. 
Second take $F_2$ and $\de_2>0$ such that 
$F_2\geq F_1\vee K_1\vee \ovl{K_1}\vee S_2$ and 
$(9\de_1^{1/4}+3 \de_2^{1/2})|F_1|_\vph<(1/2)\de_0$ and 
$9\de_2^{1/4}|F_2|_\vph <(1/2)\de_1$. 
Suppose we have chosen finitely supported central projections 
$F_n$, $K_n$ in $\lhG$ and $\de_n>0$ with 
$9\de_{n}^{1/4}|F_{n}|_\vph<(1/2)\de_{n-1}$. 
Then take $F_{n+1}$ and $\de_{n+1}$ 
such that $F_{n+1}\geq F_n\vee K_n\vee \ovl{K_n}\vee S_{n+1}$, 
$(9\de_{n}^{1/4}+3\de_{n+1}^{1/2})|F_{n}|_\vph<(1/2)\de_{n-1} $ 
and 
$9\de_{n+1}^{1/4}|F_{n+1}|_\vph<(1/2)\de_{n}$. 
Then take $K_{n+1}$ which is  
$(F_{n+1},\de_{n+1})$-invariant 
and $K_{n+1}\geq e_\btr$. 
Set $\mF_n=\supp(F_n)$ and $\mK_n=\supp(K_n)$. 
Since $F_n\geq S_n$, 
$\displaystyle\cup_{n=0}^\infty \mF_n=\IG$. 

Let $(\al,u)$ be an approximately inner cocycle action of $\bhG$ 
on $M$. 
By Lemma \ref{lem: reduction-to-action}, 
we can take a unitary $v$ in $M^\om\oti\lhG$ such that 
$\al=\Ad v$ on $M\subs M^\om$, 
\[(v^*\oti1)\al^\om(v^*)u(\id\oti\De)(v)=1,
\]
\[
\lim_{n\to\infty}
\|(\ph\oti\ta_\pi)\circ\Ad v_\pi^*\circ \al_\pi-\ph\|=0
\quad \mbox{for all}\ \ph\in M_*,\ \pi\in\IG.
\]
Then $\ga=\Ad v^*\circ\al^\om$ is an action on $M^\om$ 
fixing $M$ and preserving $M_\om$. 
Assume that 
\[
\|u_{\pi,\rho}-1\oti1_\pi\oti1_\rho\|_{\ta\oti\vph\oti\vph}<\de_{n+1}
\] 
for all $\pi, \rho\in\mF_{n+1}$. 
Then $v$ is an approximate coboundary, 
and by Theorem \ref{thm: jRohlin} 
there exists a unitary $\mu\in M^\om$ such that 
\[
|v_{F_{n}}\ga_{F_{n}}(\mu)-\mu\oti F_{n}|_{\ta\oti\vph}
<(9\de_{n}^{1/4}+3\de_{n+1}^{1/2})|F_{n}|_\vph
<\de_{n-1}. 
\]
Then we have 
\begin{align*}
|(\mu\oti F_{n})v_{F_{n}}^*\al_{F_{n}}^\om(\mu^*)
-1\oti F_{n}|_{\ta\oti\vph}
=&\,
|(\mu\oti F_{n})\ga_{F_{n}}(\mu^*)v_{F_{n}}^*
-1\oti F_{n}|_{\ta\oti\vph}
\\
=&\,
|\mu\oti F_{n}-v_{F_{n}}\ga_{F_{n}}(\mu)|_{\ta\oti\vph}
\\
<&\,
\de_{n-1}.
\end{align*}
Set a perturbed unitary 
$\wdt{v}=(\mu\oti1)v^*\al^\om(\mu^*)$ 
and then we have 
\begin{align*}
&(\wdt{v} \oti1)\al^\om(\wdt{v})u(\id\oti\De)(\wdt{v}^*)=1,
\\
&|\wdt{v}_{F_{n}}-1\oti F_{n}|_{\ta\oti\vph}<\de_{n-1}.
\end{align*}
Let $(\wdt{v}_m)_m$ be a representing sequence of $\wdt{v}$. 
Then there exists $m\in\N$ such that 
\begin{align*}
&\|((\wdt{v}_m)_\pi \oti1_\rho)\al_\pi((\wdt{v}_m)_\rho)
u_{\pi,\rho}(\id\oti{}_\pi\De_\rho)(\wdt{v}_m^*)
-1\oti1_\pi\oti1_\rho\|_{\ta\oti\vph\oti\vph}<\de_{n+2},
\\
&|\wdt{v}_m(1\oti F_{n})-1\oti F_{n}|_{\ta\oti\vph}<\de_{n-1}
\end{align*}
for all $\pi,\rho\in \mF_{n+2}$. 
Hence we obtain the following lemma. 

\begin{lem}
Let $M$ be a McDuff factor of type \II\ and 
$(\al,u)$ an approximately inner strongly free cocycle action 
of $\bhG$ on $M$. 
If the inequality 
\[
\|u_{\pi,\rho}-1\oti1_\pi\oti1_\rho\|_{\ta\oti\vph\oti\vph}<\de_{n+1}
\]
holds for all $\pi,\rho\in\mF_{n+1}$, 
then there exists a unitary $w\in M\oti \lhG$ such that 

\begin{enumerate}
\item

for all $\pi,\rho \in \mF_{n+2}$,
\[
\big{\|}
(w_\pi\oti1_\rho)\al_\pi(w_\rho)u_{\pi,\rho}(\id\oti{}_\pi\De_\rho)(w^*)
-1\oti1_\pi\oti1_\rho
\big{\|}_{\ta\oti\vph\oti\vph}<\de_{n+2},
\] 

\item
$|w_{F_{n}}-1\oti F_{n}|_{\ta\oti\vph}<\de_{n-1}.
$ 
\end{enumerate}

\end{lem}

Then as similar to the proof of \cite[Theorem 7.6]{Oc1}, 
we obtain a 2-cohomology vanishing result in $M$. 

\begin{thm}[2-cohomology vanishing theorem]\label{thm: 2vanish}
Let $M$ be a McDuff factor of type \II\  with the tracial state $\ta$. 
Let $(\al,u)$ be an approximately inner 
strongly free cocycle action of $\bhG$ on $M$. 
Then $u$ is a coboundary. 
Moreover, 
assume for fixed $n\geq2$, the inequality 
\[
\|u_{\pi,\rho}-1\oti1_\pi\oti1_\rho\|_{\ta\oti\vph\oti\vph}<\de_{n+1}
\] 
holds for all $\pi, \rho\in \mF_{n+1}$. 
Then there exists a unitary $w\in M\oti\lhG$ such that
\begin{enumerate}

\item 
$(w\oti1)\al(w)u(\id\oti\De)(w^*)=1$, 

\item
$|w_{F_n}-1\oti F_n|_{\ta\oti\vph}<\de_{n-2}$. 
\end{enumerate}

\end{thm}

\begin{proof}
It suffices to prove the theorem only in the case that 
$\big{\|}
u_{\pi,\rho}-1\oti1_\pi\oti1_\rho
\big{\|}_{\ta\oti\vph\oti\vph}<\de_{n+1}$ 
holds for all $\pi, \rho\in \mF_{n+1}$. 
By the previous lemma, 
there exists a unitary $w\in M\oti\lhG$ such that 
\begin{enumerate}[(i)]
\item 
for all $\pi,\rho\in \mF_{n+2}$, 
\[\big{\|}
(w_\pi\oti1_\rho)\al_\pi(w_\rho)u_{\pi,\rho}(\id\oti{}_\pi\De_\rho)(w^*)
-1\oti1_\pi\oti1_\rho
\big{\|}_{\ta\oti\vph\oti\vph}<\de_{n+2},
\] 
\item
$
|w_{F_{n}}-1\oti F_{n}|_{\ta\oti\vph}<\de_{n-1}.
$ 
\end{enumerate}
Set $u^n=u$, $w^n=w$, $\al^{n+1}=\Ad w^n\circ\al$ and 
$u^{n+1}=(w_\pi \oti1_\rho)\al_\pi(w_\rho)u(\id\oti_\pi\De_\rho)(w^*)$. 
Then $(\al^{n+1},u^{n+1})$ is a strongly free cocycle action on $M$ 
with 
$\|u_{\pi,\rho}^{n+1}-1\oti1_\pi\oti1_\rho\|_{\ta\oti\vph\oti\vph}
<\de_{n+2}$ 
for all $\pi, \rho\in \mF_{n+2}$. 
With a repetition of the above argument, 
we get a family of cocycle actions 
$\{(\al^m,u^m)\}_{m\geq n}$ and unitaries $\{w^m\}_{m\geq n}$ satisfying 
the following conditions
\begin{enumerate}[(1.$m$)]

\item
$\al^{m+1}=\Ad w^{m}\circ \al^m$, 

\item
$u^{m+1}=(w^m\oti1)\al^m(w^m)u^m (\id\oti\De)(w^{m*})$,

\item
$\|u_{\pi,\rho}^m-1\oti1_\pi\oti1_\rho\|_{\ta\oti\vph\oti\vph}<\de_{m+1}$ 
for all $\pi, \rho\in \mF_{m+1}$,

\item
$|w_{F_m}^m-1\oti F_m|_{\ta\oti\vph}<\de_{m-1}$. 

\end{enumerate}
Then set $\ovl{w}^m=w^m w^{m-1}\dots w^n$ 
and we have $\al^{m+1}=\Ad \ovl{w}^m \circ \al$ and 
$u^{m+1}=(\ovl{w}^m\oti1)\al(\ovl{w}^m)u(\id\oti\De)(\ovl{w}^{m*})$. 
By the condition (4.$m$), it is easy to see that the sequence of unitaries 
$\{\ovl{w}^m\}_{m\geq n}$ strongly converges to a unitary $\ovl{w}$. 
Then by (3.$m$), 
$(\ovl{w}\oti1)\al(\ovl{w})u(\id\oti\De)(\ovl{w}^*)=1$ holds. 
Moreover we have 
\begin{align*}
|\ovl{w}_{F_n}^m-1\oti F_n|_{\ta\oti\vph}
\leq&\,
|w_{F_n}^m-1\oti F_n|_{\ta\oti\vph}
+\dots+
|w_{F_n}^n-1\oti F_n|_{\ta\oti\vph}
\\
<&\,
\de_{m-1}+\cdots + \de_{n-1}
\\
\leq&\,
\de_{n-1}(1+1/2+1/{2^2}+\dots)
\\
<&\,
\de_{n-2}. 
\end{align*}
Hence we are done. 
\end{proof}

By virtue of 2-cohomology vanishing, 
we can show the following results. 

\begin{cor}
Let $M$ be a McDuff factor of type \II\  with the tracial state $\ta$ 
and $\al$ an approximately inner 
strongly free action of $\bhG$ on $M$. 
Let $v$ be a unitary in $M\oti\lhG$. 
Assume for fixed $n\geq2$, the inequality 
\[
\|(v_\pi\oti1_\rho)\al_\pi(v_\rho)(\id\oti_\pi\De_\rho)(v^*)
-1\oti1_\pi\oti1_\rho\|_{\ta\oti\vph\oti\vph}<\de_{n+1}
\]
holds for all $\pi, \rho\in \mF_{n+1}$. 
Then there exists a unitary $w\in M\oti\lhG$ such that
\begin{enumerate}

\item 
$(wv\oti1)\al(wv)(\id\oti\De)(((wv)^*)=1$, 

\item
$|w_{F_n}-1\oti F_n|_{\ta\oti\vph}<\de_{n-2}$. 
\end{enumerate}
\end{cor}

\begin{proof}
Let $\wdt{\al}=\Ad v\circ\al$ and 
$\wdt{u}=(v\oti1)\al(v)(\id\oti\De)(v^*)$. 
Apply the previous theorem to the cocycle action 
$(\wdt{\al},\wdt{u})$ and it is done. 
\end{proof}
 
With Lemma \ref{lem: coc-ultra} and the previous corollary, 
we can show the following result. 

\begin{cor}\label{cor: albe}
Let $M$ be a McDuff factor of type \II\ with the tracial state $\ta$, 
and  $\al$, $\be$  approximately inner strongly free actions. 
Then for any $\vep>0$, finite sets $\mF\Subs \IG$ 
and $T\Subs M$, there exists an $\al$-cocycle $v$ satisfying 
\[
\big{\|}\be_\pi(x)-\Ad v_\pi (\al_\pi(x))
\big{\|}_{\ta\oti\vph}<\vep 
\]
for all $\pi\in \mF$ and $x\in T$. 
\end{cor}

\begin{proof}
By Lemma \ref{lem: coc-ultra}, we can find an $\al^\om$-cocycle 
$W\in M^\om\oti \lhG$ with $\be=\Ad W\circ \al$ on $M$. 
Let $(w^\nu)_{\nu=0}^\infty$ be a representing 
sequence of $W$. 
Since $(W\oti1)\al^\om (W)(\id\oti\De)(W^*)=1$, for 
all $\pi,\rho\in\IG$ we have 
\[\lim_{\nu\to\om}
\big{\|}
(w_\pi^\nu \oti1_\rho)\al_\pi(w_\rho^\nu)(\id\oti_\pi\De_\rho)(w^{\nu*})
-1\oti1_\pi\oti1_\rho
\big{\|}_{\ta\oti\vph\oti\vph}=0.
\] 
Take a large $n$ so that $\mF\subs\mF_{n}$ 
and $\sup_{x\in T}(1+2\sqr{2}\|x\|)\de_{n-2}^{1/2}<\vep$. 
Also take a large $\nu$ so that 
$\big{\|}
\be_\pi(x)-\Ad w_\pi^\nu (\al_\pi(x))
\big{\|}_{\ta\oti\vph}<\de_{n-2}$ and 
\[
\big{\|}
(w_\pi^\nu \oti1_\rho)\al_\pi(w_\rho^\nu)(\id\oti_\pi\De_\rho)(w^{\nu*})
-1\oti1_\pi\oti1_\rho
\big{\|}_{\ta\oti\vph\oti\vph}<\de_{n+1} 
\] 
for all $\pi, \rho\in \mF_{n+1}$ and $x\in T$. 
Apply the previous corollary 
and then we have a unitary $\ovl{v}\in M\oti \lhG$ 
such that 
$|\ovl{v}_{F_n}-1\oti F_n|_{\ta\oti\vph}<\de_{n-2}$ 
and $\ovl{v}w^\nu$ is an $\al$-cocycle. 
Set $v=\ovl{v}w^\nu$ and then 
\begin{align*}
&
\big{\|}
\be_\pi(x)-\Ad v_\pi(\al_\pi(x))
\big{\|}_{\ta\oti\vph}
\\
\leq&\,
\big{\|}
\be_\pi(x)-\Ad w_\pi^\nu (\al_\pi(x))
\big{\|}_{\ta\oti\vph}
+
\big{\|}
\Ad w_\pi^\nu (\al_\pi(x))-\Ad v_\pi (\al_\pi(x))
\big{\|}_{\ta\oti\vph}
\\
\leq&\,
\big{\|}
\be_\pi(x)-\Ad w_\pi^\nu (\al_\pi(x))
\big{\|}_{\ta\oti\vph}
+
2\|w_\pi^\nu-v_\pi\|_{\ta\oti\vph}\|x\|
\\
=&\,
\big{\|}
\be_\pi(x)-\Ad w_\pi^\nu (\al_\pi(x))
\big{\|}_{\ta\oti\vph}
+
2\|1\oti1_\pi-\ovl{v}_\pi\|_{\ta\oti\vph}\|x\|
\\
\leq&\,
\big{\|}
\be_\pi(x)-\Ad w_\pi^\nu (\al_\pi(x))
\big{\|}_{\ta\oti\vph}
+
2\sqr{2}|1\oti1_\pi-\ovl{v}_\pi|_{\ta\oti\vph}^{1/2}\|x\|
\\
<&\,
\de_{n-2}+2\sqr{2}\de_{n-2}^{1/2}\|x\|
\\
<&\,\vep
\end{align*}
for all $\pi\in \mF_n$ and $x\in T$. 
\end{proof}

\subsection{Shapiro unitary}
We represent a Rohlin type theorem for two actions in order to study 
a commutation property of a Shapiro unitary, 
which is stated in Theorem \ref{thm: jRohlin-action} (9). 
We will explain a motivation for the study by considering 
a simple case. 
Let $\al$ be an action on a von Neumann algebra $M$ 
and $v$  an $\al$-cocycle. 
Put $\be=\Ad v\circ\al$. 
Let $K\in \Proj(Z(\lhG))$ as before. 
Assume that $\al^\om$ has a Rohlin projection $E\in M^\om\oti \lhG K$ 
in the sense of Theorem \ref{thm: jRohlin}. 
In addition, we assume that 
the projection $vEv^*$ is a Rohlin projection for $\be$. 
Now take a finite subset $T\subs M$. 
Suppose that 
$[E_\rho\oti1_\orho, \al_\rho(\al_\orho(x))]=0$ 
and 
$[v_\rho\oti1_\orho, \al_\rho(\al_\orho(x))]=0$ 
for all $x\in T$ and $\rho\in\mK$. 
Then the Shapiro unitary $\mu=(\id\oti\vph)(vE)$ 
commutes with all $x\in T$. 
We present an approximate version of the above argument. 

\begin{thm}\label{thm: jRohlin-action}
Let $M$ be a McDuff factor of type \II\ with the tracial state $\ta$ 
and $\al$ an approximately inner strongly free action of 
$\bhG$ on $M$. 
Let $0<\de<1$, $\vep>0$ and $F\in\Projf(Z(\lhG))$. 
Take an $(F,\de)$-invariant $K\in\Projf(Z(\lhG))$ 
with $K\geq e_\btr$. 
Set $\mF=\supp(F)$ and $\mK=\supp(K)$. 
Let $v$ be a unitary $\al$-cocycle in $M\oti \lhG$ and 
set a perturbed strongly free action $\be=\Ad v\circ \al$. 
Then for any countable set $S\subs M^\om$ and a finite subset 
$T\subs M_1$, 
there exist projections $E^\al$ and $E^\be$ in $M^\om\oti\lhG$ 
satisfying the following conditions. 

\begin{enumerate}

\item $E^\al=E^\al(1\oti K)$, $E^\be=E^\be(1\oti K)$. 

\item The following splitting properties of $\ta^\om$, 
\[(\ta^\om\oti\id)(xE^\al)=(\ta^\om\oti\id)(x)(\ta^\om\oti\id)(E^\al),\] 
\[(\ta^\om\oti\id)(xE^\be)=(\ta^\om\oti\id)(x)(\ta^\om\oti\id)(E^\be)\] 
hold for all $x\in S\oti\lhG K$. 

\item (approximate equivariance)
\[\big{|}\al_F^\om(E^\al)-(\id\oti_F\De_K)(E^\al)\big{|}_{\ta\oti\vph}
<5\de^{1/2}|F|_\vph, \]
\[\big{|}\be_F^\om(E^\be)-(\id\oti_F\De_K)(E^\be)\big{|}_{\ta\oti\vph}
<5\de^{1/2}|F|_\vph. \]

\item 
Decompose $E^\al$ and $E^\be$ as 
\[
E^\al=\sum_{\rho\in\mK}
\sum_{i,j\in I_\rho}d_\rho^{-1} f_{\orho_{i,j}}^\al\oti e_{\rho_{i,j}}, 
\]
\[
E^\be=\sum_{\rho\in\mK}
\sum_{i,j\in I_\rho}d_\rho^{-1} f_{\orho_{i,j}}^\be\oti e_{\rho_{i,j}}.
\]
Then $\{f_{\orho_{i,j}}^\al\}$ and $\{f_{\orho_{i,j}}^\be\}$ satisfy 
\[
f_{\orho_{i,j}}^\al f_{\opi_{k,\el}}^\al
=\de_{\rho,\pi}\de_{j,k}f_{\orho_{i,\el}}^\al,
\quad
f_{\orho_{i,j}}^\be f_{\opi_{k,\el}}^\be
=\de_{\rho,\pi}\de_{j,k}f_{\orho_{i,\el}}^\be
\]
for all $\rho,\pi\in\mK$, $i,j\in I_\rho$ and $k,\el\in I_\pi$. 

\item (joint property of $U$)
Set $U=vE^\al$ and decompose $U$ as 
\[
U=\sum_{\rho\in\mK}
\sum_{i,j\in I_\rho}d_\rho^{-1} \mu_{\orho_{i,j}}\oti e_{\rho_{i,j}}.
\]
Then we have 
\[
\mu_{\orho_{i,j}}^*\mu_{\opi_{k,\el}}
=\de_{\rho,\pi}\de_{i,k}f_{\orho_{j,\el}}^\al,
\]
\[
\mu_{\orho_{i,j}}\mu_{\opi_{k,\el}}^*
=\de_{\rho,\pi}\de_{j,\el}f_{\orho_{i,k}}^\be
\]
for all $\rho,\pi\in\mK$, $i,j\in I_\rho$ and $k,\el\in I_\pi$. 
In particular, $U^* U=E^\al$ and $UU^*=E^\be$ holds. 

\item 
For each $\rho\in\mK$, 
the projections $(\id\oti\vph)(E_\rho^\al)$ and $(\id\oti\vph)(E_\rho^\be)$ 
are equal. 
In addition, they are in $S'\cap M_\om$. 

\item (partition of unity)
\[(\id\oti\vph)(E^\al)=1=(\id\oti\vph)(E^\be). 
\]
\item (Shapiro lemma)
Set $\mu=(\id\oti\vph)(U)$ and then $\mu$ is a unitary satisfying 
\[|v_F \al^\om(\mu)-\mu\oti F|_{\ta\oti\vph}
<9\de^{1/4}|F|_\vph. 
\]
\item 
Further assume 
\[\|\be_\pi(\al_{\opi}(x))-\al_\pi(\al_{\opi}(x))
\|_{\ta\oti\vph}<\vep\] 
for all $x\in T$ and $\pi\in\mK$, 
then the unitary $\mu$ satisfies 
\[
|[\mu,x]|_\ta<\vep
\]
for all $x\in T$. 

\end{enumerate}

\end{thm}

\begin{proof}
As in Lemma \ref{lem: reduction-to-action}, 
we take a unitary $V\in M^\om\oti\lhG$ such that 
$\al=\Ad V$ on $M$ and $V^*$ is an $\al^\om$-cocycle. 
Set the strongly free action $\ga=\Ad V^*\circ\al^\om$. 
We use the same notations in Theorem \ref{thm: jRohlin}. 
We may assume that $S$ contains the entries of $\al_\pi(x)$ 
for all $\pi\in\IG$ and $x\in T$. 
Recall the set $\meJ$ defined in \S5.3. 
We denote by $\meI$ the subset of $\meJ$ 
whose elements satisfy the following conditions. 
For $\ovl{E}\in\meI$, 
\begin{enumerate}[(i)]
\item
$[\ovl{E},\ga(x)]=0$ for all $x\in S$. 

\item 
The projections 
$\ovl{E}^\al:=V\ovl{E}V^*$ and $\ovl{E}^\be:=vV\ovl{E}V^*v^*$ 
satisfy the conditions (1), (2), (4), (5) and (6) 
in Theorem \ref{thm: jRohlin-action}. 

\item The equality 
$(\id\oti\vph)(\ovl{E}_\rho^\al)
=(\id\oti\vph)(\ovl{E}_\rho)=(\id\oti\vph)(\ovl{E}_\rho^\be)$ 
holds for all $\rho\in\mK$. 

\item 
$(\ta^\om\oti\id)(\ovl{E}^\al)=b_{\ovl{E}}|K|_\vph^{-1}K
=(\ta^\om\oti\id)(\ovl{E}^\be)$. 

\end{enumerate}

Recall a subset $\meS\subs \meJ$ defined in the proof of 
Theorem \ref{thm: jRohlin}. 
Then we can easily see that $\meI\cap \meS$ 
is an inductive ordered set as similar to the proof of Theorem 
\ref{thm: jRohlin}. 
Let $\ovl{E}$ be a maximal element of $\meI$. 
Since the proof of Lemma \ref{lem: E'} is applicable for 
$\meI$ with the additional assumption (i), 
the projection $p=1-(\id\oti\vph)(\ovl{E})$ 
satisfies $\ta_\om(p)\leq\de^{1/2}$. 
Let $E^\al=VEV^*$ and $E^\be=vVEV^*v^*$ where 
$E=\ovl{E}+p\oti e_\btr$. 
Then the conditions (1), (2), (4), (5), (6), 
(7) of Theorem \ref{thm: jRohlin-action} are satisfied. 
Since $V^*$ is an $\al^\om$-cocycle, 
the condition (3) follows from the tracial property of $\ta$. 
On (8), 
a similar proof to that of 
Theorem \ref{thm: jRohlin} (8) 
is applicable to $(\al,E^\al)$ and an $\al$-cocycle $v$. 
Finally we verify (9). 
For $\rho\in\mK$ and $x\in T$, 
\begin{align*}
&\|[v_\rho\oti1_\orho,\al_\rho(\al_\orho(x))](\ovl{E}_\rho^\al\oti1_\orho)
\|_{\ta\oti\vph\oti\vph}^2
\\
=&\,
(\ta\oti\vph\oti\vph)
\big{(}(\ovl{E}_\rho^\al\oti1_\orho)
|[v_\rho\oti1_\orho,\al_\rho(\al_\orho(x))]|^2
(\ovl{E}_\rho^\al\oti1_\orho)\big{)}
\\
=&\,
(\ta\oti\vph\oti\vph)
\big{(}
|[v_\rho\oti1_\orho,\al_\rho(\al_\orho(x))]|^2
((\ta^\om\oti\id)(\ovl{E}_\rho^\al)\oti1_\orho)\big{)}
\\
=&\,
(\ta\oti\vph\oti\vph)
\big{(}
|[v_\rho\oti1_\orho,\al_\rho(\al_\orho(x))]|^2
(b_{\ovl{E}}|K|_\vph^{-1}1_\rho\oti1_\orho)\big{)}
\\
=&\,
b_{\ovl{E}}|K|_\vph^{-1}
\big{\|}[v_\rho\oti1_\orho,\al_\rho(\al_\orho(x))]
\big{\|}_{\ta\oti\vph\oti\vph}^2
\\
<&\,
\vep^2 |K|_\vph^{-1}. 
\end{align*}

Using it, for $\rho\in\mK$ and $x\in T$, 
\begin{align*}
|[U_\rho\oti1_\orho,\al_\rho(\al_\orho(x))]|_{\ta\oti\vph\oti\vph}
=&\,
|[vE^\al\oti1_\orho,\al_\rho(\al_\orho(x))]|_{\ta\oti\vph\oti\vph}
\\
=&\,
|[v\ovl{E}^\al\oti1_\orho,\al_\rho(\al_\orho(x))]|_{\ta\oti\vph\oti\vph}
\\
\leq&\,
|[v_\rho\oti1_\orho,\al_\rho(\al_\orho(x))] 
(\ovl{E}_\rho^\al\oti1_\orho)|_{\ta\oti\vph\oti\vph}
\\
&\quad
+
|(v_\rho\oti1_\orho) 
[\ovl{E}_\rho^\al\oti1_\orho,\al_\rho(\al_\orho(x))]
|_{\ta\oti\vph\oti\vph}
\\
=&\,
|[v_\rho\oti1_\orho,\al_\rho(\al_\orho(x))] 
(\ovl{E}_\rho^\al\oti1_\orho)|_{\ta\oti\vph\oti\vph}
\\
\leq&\,
\|[v_\rho\oti1_\orho,\al_\rho(\al_\orho(x))](\ovl{E}_\rho^\al\oti1_\orho)
\|_{\ta\oti\vph\oti\vph}
\\
&\quad\cdot
\|\ovl{E}_\rho^\al\oti1_\orho
\|_{\ta\oti\vph\oti\vph}
\\
<&\,
\vep\|K\|_\vph^{-1}\|\ovl{E}_\rho^\al\|_{\ta\oti\vph}d_\rho
\\
=&\,
\vep\|K\|_\vph^{-1}b_{\ovl{E}}^{1/2}\|K\|_\vph^{-1}d_\rho^2
\\
\leq&\,
\vep d_\rho^2|K|_\vph^{-1}. 
\end{align*}
Set a state 
$\th_\rho=T_{\rho,\orho}^*\cdot T_{\rho,\orho}$ on $B(H_\rho\oti H_\orho)$. 
Then we have $\vph_\rho(a)=d_\rho^2 \th_\rho(a\oti1_\orho)$ for 
all $a\in B(H_\rho)$ and 
$d_\rho^2 \th_\rho\leq \vph_\rho\oti\vph_\orho$ as positive functionals. 
We claim that 
\[
d_\rho^2 |(\id\oti\th_\rho)(x)|_\ta\leq 
|x|_{\ta\oti\vph_\rho\oti\vph_\orho}
\quad\mbox{for all}\ x\in M^\om \oti B(H_\rho\oti H_\orho). 
\]
Let $(\id\oti\th_\rho)(x)=w|(\id\oti\th_\rho)(x)|$ and 
$x=w'|x|$ be the polar decompositions. 
Then the claim is verified as follows
\begin{align*}
&d_\rho^2 |(\id\oti\th_\rho)(x)|_\ta
\\
=&\,
d_\rho^2 \ta(w^*(\id\oti\th_\rho)(x))
\\
=&\,
d_\rho^2(\ta\oti\th_\rho)((w^*\oti1_\rho\oti1_\orho)x)
\\
=&\,
d_\rho^2(\ta\oti\th_\rho)((w^*\oti1_\rho\oti1_\orho)w'|x|^{1/2} |x|^{1/2})
\\
\leq&\,
d_\rho^2 
(\ta\oti\th_\rho)\big{(}(w^*\oti1_\rho\oti1_\orho)w'|x|{w'}^*
(w\oti1_\rho\oti1_\orho)\big{)}^{1/2}
(\ta\oti\th_\rho)(|x|)^{1/2}
\\
\leq&\,
(\ta\oti\vph_\rho\oti\vph_\orho)
\big{(}(w^*\oti1_\rho\oti1_\orho)w'|x|{w'}^*
(w\oti1_\rho\oti1_\orho)\big{)}^{1/2}
(\ta\oti\vph_\rho\oti\vph_\orho)(|x|)^{1/2}
\\
=&\,
(\ta\oti\vph_\rho\oti\vph_\orho)
\big{(}|x|^{1/2}{w'}^*(ww^*\oti1_\rho\oti1_\orho)w'
|x|^{1/2}\big{)}^{1/2}
(\ta\oti\vph_\rho\oti\vph_\orho)(|x|)^{1/2}
\\
\leq&\,
(\ta\oti\vph_\rho\oti\vph_\orho)(|x|)^{1/2}
(\ta\oti\vph_\rho\oti\vph_\orho)(|x|)^{1/2}
\\
=&\,
|x|_{\ta\oti\vph_\rho\oti\vph_\orho}.
\end{align*}

Then we obtain 
\begin{align*}
|[\mu,x]|_\ta
=&\,
|[(\id\oti\vph)(U),x]|_\ta
\\
=&\,
\Big{|}\sum_{\rho\in\mK}
\big{[}
d_\rho^2(\id\oti\th_\rho)(U_\rho\oti1_\orho),x
\big{]}
\Big{|}_\ta
\\
=&\,
\Big{|}\sum_{\rho\in\mK}
d_\rho^2(\id\oti\th_\rho)([U_\rho\oti1_\orho,\al_\rho(\al_\orho(x))])
\Big{|}_\ta
\\
\leq&\,
\sum_{\rho\in\mK}
\big{|}
d_\rho^2
(\id\oti\th_\rho)
([U_\rho\oti1_\orho,\al_{\rho}(\al_\orho(x))])
\big{|}_\ta
\\
\leq&\,
\sum_{\rho\in\mK}
\big{|}
[U_\rho\oti1_\orho,\al_{\rho}(\al_\orho(x))]
\big{|}_{\ta\oti\vph\oti\vph}
\\
\leq&\,
\sum_{\rho\in\mK}
\vep d_\rho^2 |K|_\vph^{-1}
\\
=&\,\vep. 
\end{align*} 
\end{proof}

Although a 1-cohomology does not vanish in $M$ in general, 
it approximately vanishes. 
The following theorem shows the approximate vanishing 
with a commutation property of a Shapiro unitary. 
Since it is easily proved 
by considering a representing sequence of $\mu$ in the previous theorem, 
we omit the proof. 

\begin{thm}[Approximate vanishing of 1-cohomology]
\label{thm: appv1}
Let $M$ be a McDuff factor of type \II\ with the tracial state $\ta$ 
and $\al$ a strongly free action of $\bhG$ on $M$. 
Let $F\in\Projf(Z(\lhG))$ and $\de,\vep>0$. 
Take an $(F,\de)$-invariant $K\in\Projf(Z(\lhG))$ with $K\geq e_\btr$. 
Let $T$ be a finite subset in the unit ball of $M$. 
If an $\al$-cocycle $v$ satisfies 
\[
\|\Ad v_\rho\circ\al_\rho(\al_\orho(x))
-\al_\rho(\al_\orho(x))
\|_{\ta\oti\vph\oti\vph}<\vep
\] 
for all $x\in T$ and $\rho\in K$,
then there exists a unitary $w$ in $M$ satisfying 
\begin{enumerate}[(i)]

\item
$|v_F-(w\oti1)\al_F(w^*)|_{\ta\oti\vph}
<9\de^{1/4}|F|_\vph$, 

\item
$|[w,x]|_\ta<\vep,\,x\in T$. 

\end{enumerate}
\end{thm}

\section{Main theorem}

\subsection{Intertwining argument}

For a proof of the cocycle conjugacy of two actions, 
we make use of so-called an intertwining argument initiated by 
Evans-Kishimoto in \cite{EvKi}. 
The results Corollary \ref{cor: albe} and Theorem \ref{thm: appv1} 
are necessary for the argument. 
We briefly explain the outline. 
Let $\ga^0:=\al$ and $\ga^{-1}:=\be$ be approximately inner strongly free 
actions on $M$. 
First by Corollary \ref{cor: albe}, 
we perturb the action $\ga^{-1}$ to $\ga^1$ 
by a $\ga^{-1}$-cocycle $v^1$ so that 
$\ga^1$ is close to $\ga^0$. 
Second by Corollary \ref{cor: albe}, 
we perturb the action $\ga^0$ to $\ga^2$ 
by a $\ga^0$-cocycle $v^2$ so that $\ga^2$ is close to $\ga^1$. 
We construct families of actions and 1-cocycles inductively 
and achieve the equality at the limit. 
However in that process, 
we have to use Theorem \ref{thm: appv1} 
in order to treat successive multiplications of unitaries. 

\begin{thm}
Let $M$ be a McDuff factor of type \II. 
Let $\al$ and $\be$ be approximately inner strongly free actions 
of $\bhG$ on $M$. 
Then they are cocycle conjugate, that is, 
there exist an automorphism $\th$ in $\oInt(M)$ and 
an $\al$-cocycle $v$ with 
\[
\Ad v\circ \al=(\th^{-1}\oti\id)\circ\be\circ\th. 
\]

\end{thm}

\begin{proof}
Let $S=\{a_i\}_{i=1}^\infty$ be a strongly dense countable subset 
of the unit ball of $M$. 
Put $S_n=\{a_i\}_{i=1}^n$ and $\vep_n=2^{-n}$. 
Take the sequences of the finitely supported central projections 
$\{F_n\}_{n=1}^\infty$, $\{K_n\}_{n=1}^\infty$ 
and the positive numbers 
$\{\de_n\}_{n=1}^\infty$ as in \S6.1. 
Set $\ga^0=\al$, $\ga^{-1}=\be$, 
$\ovl{u}_{-1}=\ovl{u}_0=1$, 
$\th_{0}=\th_{-1}=\id\in \Int(M)$ and 
$T_0=\{1\}\Subs M$. 
For each $n\geq1$, 
we construct inductively the following members. 
\begin{enumerate}[(i)]

\item
an action $\ga^n$ of $\bhG$ on $M$, 

\item 
$w_n\in U(M)$, 

\item 
$\th_n\in \Int(M)$, 

\item 
an $\Ad(w_n\oti1)\circ\ga^{n-2}\circ\Ad w_n^*$-cocycle $u^n$, 

\item 
a $(\th_n\oti\id)\circ\ga^{\ovl{n}}\circ \th_n^{-1}$-cocycle $\ovl{u}^n$, 
where $\ovl{n}$ is equal to 0 or $-1$ 
according to that $n$ is even or odd respectively, 

\item 
a finite subset $T_n\Subs M$. 
\end{enumerate}
The induction conditions are

\begin{enumerate}[(1,$n$)]

\item 
$\|\ga_\rho^n(x)-\ga_\rho^{n-1}(x)\|_{\ta\oti\vph}<\vep_n$ 
for $x\in S_n$ and $\rho\in\mF_n$ ($n\geq1$), 

\item
$\|\ga_\rho^n(\ga_\orho^k(x))
-\ga_\rho^{n-1}(\ga_\orho^k(x))
\|_{\ta\oti\vph\oti\vph}<\vep_n$ 
for all $x\in T_{n-1}$, $\rho\in \mK_n\cup \mK_{n+1}$ 
and $1\leq k\leq n-1$ ($n\geq2$), 

\item
$|u_{F_n}^n-1\oti F_n|_{\ta\oti\vph}<9\de_n^{1/2}|F_n|_\vph$ 
($n\geq3$), 

\item
$|[w_n,x]|_\ta<\vep_{n-2}$, for $x\in T_{n-2}$ ($n\geq3$), 

\item
$\ovl{u}^n
=u^n (w_n\oti1)\ovl{u}^{n-2}(w_n^*\oti1)$ ($n\geq1$),

\item
$\th_n=\Ad w_n \circ \th_{n-2}$ ($n\geq1$),

\item
$\ga^n
=
\Ad u^n \circ \Ad(w_n\oti1)\circ \ga^{n-2}\circ \Ad w_n^*$ 
($n\geq1$),

\item
$T_n=T_{n-1}\cup S_n\cup \th_n(S_n)\cup 
\{\ovl{u}_{\pi_{i,j}}^n\mid \pi\in \mF_n, i,j\in I_\pi\}$ 
($n\geq1$). 

\end{enumerate}

\textbf{1st step.}

Since $\al$ and $\be$ are approximately inner and strongly free,  
there exists a $\ga^{-1}$-cocycle $v^1$ with 
$
\big{\|}
\Ad v_\rho^1(\ga_\rho^{-1}(x))-\ga_\rho^0(x)
\big{\|}_{\ta\oti\vph}
<\vep_1$ 
for all $x\in S_1$ and $\rho\in \mF_1$ 
by Corollary \ref{cor: albe}. 
Applying Theorem \ref{thm: appv1} to $v^1$, we can take 
a unitary $w_1\in M$ with 
$\big{|}
v_{F_1}^1-(w_1\oti 1)\ga_{F_1}^{-1}(w_1^*)
\big{|}_{\ta\oti\vph}
<
9\de_1^{1/4}|F_1|_\vph$. 
Then set $u^1=\ovl{u}^1=v^1 \ga^{-1}(w_1)(w_1^*\oti1)$ and 
it is an $\Ad (w_1\oti1) \circ \ga^{-1}\circ \Ad w_1^*$-cocycle. 
Set an action 
$\ga^1=\Ad v^1\circ \ga^{-1}=\Ad u^1 
\circ \Ad (w_1\oti1) \circ \ga^{-1}\circ \Ad w_1^*$, 
an automorphism $\th_1=\Ad w_1$ 
and $T_1$ as in $(8,1)$. 

\textbf{2nd step.}

Next take a $\ga^0$-cocycle $v^2$ with 
\[
\big{\|}
\Ad v_\rho^2(\ga_\rho^{0}(x))-\ga_\rho^1(x)
\big{\|}_{\ta\oti\vph}<\vep_2 
\quad\mbox{for all}\ 
x\in S_2,\  \rho\in \mF_2,
\]
\[
\big{\|}
\Ad v_\rho^2\circ\ga_\rho^{0}(\ga_\orho^k(x))
-\ga_\rho^1(\ga_\orho^k(x))
\big{\|}_{\ta\oti\vph}
<\vep_2
\quad\mbox{for all}\ 
x\in T_1,\  \rho\in \mK_2\cup\mK_3,\  k=0,1\] by Corollary \ref{cor: albe}. 
Then also by Theorem \ref{thm: appv1}, we can take 
a unitary $w_2$ with 
$
\big{|}
v_{F_2}^2-(w_2\oti 1)\ga_{F_2}^{0}(w_2^*)
\big{|}_{\ta\oti\vph}
<9\de_2^{1/4}|F_2|_\vph$. 
Set $u^2=v^2\ga^{0}(w_2)(w_2^*\oti1)$ and $\ga^2=\Ad v^2\circ\ga^0$. 
Set $\th_2=\Ad w_2$, $\ovl{u}^2=u^2$ and 
$T_2$ as in $(8,2)$. 

\mathversion{bold}
$(n+1)$ \textbf{-st step.}
\mathversion{normal}

Suppose that we have done up to $n$-th step. 
By Corollary \ref{cor: albe}, 
we can take a $\ga^{n-1}$-cocycle $v^{n+1}$ with 
\[
\big{\|}
\Ad v_\rho^{n+1}(\ga_\rho^{n-1}(x))-\ga_\rho^n(x)
\big{\|}_{\ta\oti\vph}
<\vep_{n+1} 
\quad\mbox{for all}\ 
x\in S_{n+1},\  \rho\in \mF_{n+1} 
\]
and for $0\leq k\leq n$, 
\[
\big{\|}
\Ad v_\rho^{n+1}\circ\ga_\rho^{n-1}(\ga_\orho^k(x))
-\ga_\rho^n(\ga_\orho^k(x))
\big{\|}_{\ta\oti\vph\oti\vph}
<\vep_{n+1}
\ \mbox{for all}\ 
x\in T_{n},\  \rho\in \mK_{n+1}\cup\mK_{n+2}.
\]  
Since we have the condition (2,$n$) (with $k=n-1$), i.e., 
\[
\big{\|}
\ga_\rho^{n}(\ga_\orho^{n-1}(x))
-\ga_\rho^{n-1}(\ga_\orho^{n-1}(x))
\big{\|}_{\ta\oti\vph\oti\vph}
<\vep_{n} 
\quad\mbox{for all}\ 
x\in T_{n-1},\  \rho\in \mK_{n}\cup\mK_{n+1}, 
\]
we obtain for all $x\in T_{n-1}$ and $\rho\in \mK_{n+1}$, 
\[
\big{\|}
\Ad v_\rho^{n+1}\circ\ga_\rho^{n-1}(\ga_\orho^{n-1}(x))
-\ga_\rho^{n-1}(\ga_\orho^{n-1}(x))
\big{\|}_{\ta\oti\vph\oti\vph}
<\vep_n+\vep_{n+1}.\]
Then by Theorem \ref{thm: appv1}, there exists a unitary 
$w_{n+1}$ in $M$ with 
\[
\big{|}
v_{F_{n+1}}^{n+1}-(w_{n+1}\oti1)\ga^{n-1}_{F_{n+1}}(w_{n+1}^*)
\big{|}_{\ta\oti\vph}
<9\de_{n+1}^{1/4}|F_{n+1}|_\vph
\]
and 
\[
|[w_{n+1}, x]|_{\ta\oti\vph}<\vep_n+\vep_{n+1}<\vep_{n-1}
\quad\mbox{for all}\  x\in T_{n-1}. 
\]
Set $u^{n+1}=v^{n+1}\ga^{n-1}(w_{n+1})(w_{n+1}^*\oti1)$ 
and $\ga^{n+1}=\Ad v^{n+1}\circ\ga^{n-1}$. 
Set $\ovl{u}^{n+1}$, $\th_{n+1}$ and $T_{n+1}$ as in 
$(5,n+1)$, $(6,n+1)$ and $(8,n+1)$, respectively. 
Then all the conditions have been verified. 
Thus we have constructed the members in (i),\dots,(vi) inductively. 

We show the existence of $\lim\limits_{m\rightarrow \infty}\th_{2m}$ and 
$\lim\limits_{m\rightarrow\infty}\th_{2n+1}$. 
Let $x\in S_n$. 
For $m\geq n+2$, we have 
\begin{align*}
\|\th_{m}(x)-\th_{m-2}(x)\|_{\ta}
=&\,
\|w_m\th_{m-2}(x)w_m^*-\th_{m-2}(x)\|_{\ta}
\\
=&\,
\|[w_m,\th_{m-2}(x)]\|_\ta
\\
\leq&\,
\sqr{2}
|[w_m,\th_{m-2}(x)]|_\ta^{1/2}
\\
<&\,
\sqr{2}\vep_{m-2}^{1/2}. 
\end{align*}
Similarly we have 
$\|\th_m^{-1}(x)-\th_{m-2}^{-1}(x)\|_\ta
<\sqr{2}\vep_{m-2}^{1/2}$. 
Hence the strong limits 
$\displaystyle\lim_{m\to\infty}\theta_{2m}(x)$ and 
$\displaystyle\lim_{m\to\infty}\theta_{2m+1}(x)$ exist for all $x\in S$. 
It clearly derives the existence of the limits for all $x\in M$. 
Let $\displaystyle\ovl{\th}_0=\lim_{n\to\infty}\th_{2n}$ 
and $\displaystyle\ovl{\th}_1=\lim_{n\to\infty}\th_{2n+1}$. 
Then they are approximately inner automorphisms. 

Next we show the existence of  
$\lim\limits_{m\rightarrow \infty}\ovl{u}^{2m}$
and $\lim\limits_{m\rightarrow \infty}\ovl{u}^{2m+1}$.
For $m\geq n$ we have 
\begin{align*}
\big{|}
\ovl{u}_{F_n}^{m+2}-\ovl{u}_{F_n}^m
\big{|}_{\ta\oti\vph}
=&\,
\big{|}
u_{F_n}^{m+2}(w_{m+2}\oti 1) \ovl{u}_{F_n}^m (w_{m+2}^*\oti1)
- \ovl{u}_{F_n}^m
\big{|}_{\ta\oti\vph}
\\
\leq&\,
\big{|}
u_{F_n}^{m+2}-1\oti F_n
\big{|}_{\ta\oti\vph}
+
\big{|}
[w_{m+2}\oti F_n, \ovl{u}_{F_n}^m]
\big{|}_{\ta\oti\vph}
\\
\leq&\,
\big{|}
u_{F_{m+2}}^{m+2}-1\oti F_{m+2}
\big{|}_{\ta\oti\vph}
+
\sum_{\pi\in \mF_{n}}\sum_{i,j\in I_\pi}
|[w_{m+2}, \ovl{u}_{\pi_{i,j}}^m]\oti e_{\pi_{i,j}}|_{\ta\oti\vph}
\\
<&\,
9\de_{m+2}^{1/4}|F_{m+2}|_\vph
+
\sum_{\pi\in \mF_{n}}\sum_{i,j\in I_\pi}
d_\pi |[w_{m+2}, \ovl{u}_{\pi_{i,j}}^m]|_\ta 
\\
<&\,
9\de_{m+2}^{1/4}|F_{m+2}|_\vph
+
\sum_{\pi\in \mF_{n}}\sum_{i,j\in I_\pi}
d_\pi \vep_{m+2}
\\
\leq&\,
\de_{m+1}+\vep_{m+2}|F_n|_\vph^2. 
\end{align*}
Hence $\{\ovl{u}_\pi^{2n}\}_{n=1}^\infty$ and
$\{\ovl{u}_\pi^{2n-1}\}_{n=1}^\infty$ are Cauchy sequences 
for all $\pi\in\IG$, 
and the strong limits 
$\displaystyle\hat{u}^0=\lim_{n\to\infty}\ovl{u}^{2n}$ and 
$\displaystyle\hat{u}^1=\lim_{n\to\infty}\ovl{u}^{2n+1}$ exist. 
It is easy to see that 
$\ga^{2n}=\Ad \ovl{u}^{2n}
\circ (\th_{2n}\oti\id)\circ \al\circ \th_{2n}^{-1}$ and 
$
\ga^{2n+1}=\Ad \ovl{u}^{2n+1}
\circ (\th_{2n+1}\oti\id)\circ \be\circ \th_{2n+1}^{-1}
$. 
Since for fixed $n$, 
$\displaystyle\lim_{m\to\infty}
\|\ga_\pi^{2m+1}(x)-\ga_\pi^{2m}(x)\|_{\ta\oti\vph}=0$ 
for all $x\in S_n$ and $\pi\in\IG$, 
the next equality holds on 
$\displaystyle\cup_{n\geq1}S_n$, and so does on $M$, 
\[
\Ad \hat{u}^0\circ (\ovl{\th}_0\oti\id)\circ\al\circ\ovl{\th}_0^{-1}
=
\Ad \hat{u}^1\circ (\ovl{\th}_1\oti\id)\circ\be\circ\ovl{\th}_1^{-1}.
\]
Moreover, $\hat{u}^0$ and $\hat{u}^1$ are 1-cocycles for 
$(\ovl{\th}_0\oti\id)\circ\al\circ\ovl{\th}_0^{-1}$ 
and 
$(\ovl{\th}_1\oti\id)\circ\be\circ\ovl{\th}_1^{-1}$, respectively. 
Therefore $\al$ and $\be$ are cocycle conjugate. 
\end{proof}

Since strong freeness and freeness are equivalent notions for 
the AFD factor of type \II\  (see Appendix), 
we obtain the following result. 

\begin{cor}\label{cor: uniqueness}
Any two free actions of an amenable discrete Kac algebra 
on the AFD factor of type \II\  are cocycle conjugate. 
\end{cor}

\subsection{Classification of minimal actions}

We show the uniqueness of minimal actions of 
a compact Kac algebra $\bG=(\lG,\de,h)$ 
with  amenable dual on the AFD factor of type \II.

\begin{lem}\label{lem: stability}
Let $M$ be a finite von Neumann algebra, $\alpha$ an action of 
a compact Kac algebra $\bG$ on M. 
If $M \rtimes_\alpha \bG$ is a factor, 
then any $\al$-cocycle is a coboundary.
\end{lem}
\begin{proof}
The proof is similar to that of \cite[Theorem 12]{Wa2}. 
Let $w$ be an $\al$-cocycle. 
Set $N:=M_2(\C)\oti M$, $\wdt{\al}:=\id\otimes \al$ 
and $\wdt{w}:=e_{11}\oti w + e_{22}\oti 1\oti 1$. 
Then $\wdt{w}$ is an $\wdt{\al}$-cocycle, 
and $\beta:=\Ad \wdt{w}\circ\wdt{\al}$ is an action of $\bG$ on $N$. 
Since $M_2(\C)\oti (M\rtimes_\alpha \bG)
\cong N\rtimes_{\wdt{\al}}\bG
\cong N\rtimes_\beta \bG$, $N^\beta$ is a factor by 
\cite[Corollary 5]{Se} or
\cite[Corollary 3.9]{Ya1}, 
and the restriction of any trace on $M_{2}(\C)\oti M$ is 
the unique trace on $N^\beta$. 
Since $e_{11}\oti1, e_{22}\oti1\in N^\beta$ and their values of trace 
are $1/2$, they are equivalent. 
Let $v\in N^\beta$ be such that $v^*v=e_{11}\oti1$ 
and $vv^*=e_{22}\oti1$. 
Then $v$ is of the form $v=e_{12}\oti u$ for some $u\in U(M)$, 
and $v\in N^\beta$ implies that $w=(u\otimes 1)\alpha(u^*)$.
\end{proof}

Recall that an action $\alpha$ of a compact Kac algebra $\bG$ on $M$ 
is said to be \textit{minimal} if 
$(M^{\alpha})'\cap M=\C$ and the linear span of 
$\{(\phi\otimes \id)(\alpha(M))\mid \phi\in M_*\}$ 
is weakly dense in $\lG$ \cite{ILP}. 
Since this definition is equivalent to 
$(M^{\alpha})'\cap M=\C$ and the factoriality  of
$M\rtimes_\alpha \bG$ 
by \cite[Corollary 7]{Se} and \cite[Corollary 3.10]{Ya1}, 
any 1-cocycle for a minimal action is a coboundary. 
Readers are referred to \cite{HY} or \cite{V} for 
constructions of minimal actions. 

\begin{cor}
Let $\bG=(\lG,\de,h)$ be a compact Kac algebra with  amenable dual. 
Let $M$ be the AFD factor of type \II\ and 
$\alpha$ a minimal action of $\bG$ on $M$. 
Then $\alpha$ is dual.
\end{cor}
\begin{proof}
On $M\oti B(H_\pi)$, consider a minimal action $\ovl{\al}^\pi$ defined by 
$\ovl{\al}^\pi(x)=(\al\oti\id)(x)_{132}$. 
Recall the multiplicative unitaries 
$V$ and $\tV$ defined in \S 2.2 and \S 2.5, respectively. 
Since $\tV_{23}{}V_{21}\tV_{23}^*=V_{21}V_{31}$, 
$1\otimes (V_\pi)_{21}\in M\otimes B(H_\pi)\otimes \lG$ 
is an $\ovl{\al}^\pi$-cocycle. 
By Lemma \ref{lem: stability}, 
there exists $v_\pi\in U(M\otimes B(H_\pi))$ with 
$(v_\pi^*\oti1)\ovl{\al}^\pi(v_\pi)=1\otimes(V_\pi)_{21}$, that is, 
$\al(v_\pi)=(v_\pi)_{13} (1\oti V_\pi)$. 
Define $v\in U(M\otimes \lhG)$ by 
$v(1\oti 1_\pi)=v_\pi$ for all $\pi\in\IG$. 
Set 
$\beta(x):=\Ad v(x\otimes 1)$ for $x\in M^\alpha$.
We claim that $\beta$ is a cocycle action of $\bhG$ on 
$M^\alpha$ with a 2-cocycle $u:=v_{12}v_{13}(\id_M\otimes \Delta)(v^*)$. 
Since for $x\in M^\al$, we have 
\begin{align*}
(\al\oti\id)(\be(x))
=&\,
\al(v)(\al(x)\oti1)\al(v^*)
\\
=&\,
v_{13}(1\oti V)(x\oti1\oti1)(1\oti V^*)v_{13}^*
\\
=&\,
\be(x)_{13}. 
\end{align*}
Hence $\be$ preserves $M^\al$. 
We verify $u\in M^\al\oti\lhG\oti\lhG$ as follows. 
\begin{align*}
\al(u)
=&\,
\al(v)_{123} \al(v)_{124}(\id\oti\id\oti \Delta)(\al(v^*))_{1234}
\\
=&\,
v_{13}V_{23} v_{14}V_{24}
(\id \oti\id\oti\Delta)(V_{23}^* v_{13}^*)_{1234} 
\\
=&\,
v_{13}V_{23} v_{14}V_{24} (V_{23}V_{24})^* (\id\oti\De)(v^*)_{134}
\\
=&\,
v_{13}v_{14}(\id\oti\De)(v^*)_{134}
\\
=&\,
u_{134}.
\end{align*}
Hence $(\be,u)$ is a cocycle action on $M^\al$. 
Next we show the freeness of the cocycle action $(\be,u)$. 
Assume that for an element $\pi\in\IG\setm\{\btr\}$, 
there exists a nonzero $a\in M^\al\oti B(H_\pi)$ such that 
$\be_\pi(x)a=a(x\oti1_\pi)$ for all $x\in M^\al$. 
Then since $v_\pi^*a\in ((M^\al)'\cap M)\oti B(H_\pi)=\C\oti B(H_\pi)$, 
there exists $b\in B(H_\pi)$ such that $a=v_\pi(1\oti b)$. 
Applying $\al\oti\id$ to the both sides, 
we have $a_{13}=(v_\pi)_{13}(1\oti V_\pi)(1\oti1\oti b)$ and hence 
$1\oti1\oti b=(V_\pi)_{23}(1\oti1\oti b)\in \C\oti \lG_\pi\oti B(H_\pi)$, 
but this is a contradiction. 
Therefore $(\be,u)$ is a free cocycle action 
on the AFD factor $M^\al$ of type \II. 
By Theorem \ref{thm: 2vanish}, 
there exists $w\in U(M^\alpha\otimes \lhG)$ 
with $(w\otimes 1)
(\beta\otimes \id)(w)u(\id\otimes \Delta)(w^*)
=(w_{12}v_{12})(w_{13}v_{13})(\id\otimes \Delta)(v^*w^*)=1$. 
Then $wu$ is a 
unitary representation of $\lhG$. 
Since $\al(v)=v_{13}V_{23}$ 
and $\al(w)=w_{13}$, it follows 
$\al(wv)=(wv)_{13}V_{23}$. 
Hence $\al$ is a dual action 
for a free action $\Ad w\circ\be$ on $M^\al$.
\end{proof}

In the end, we prove the following main theorem. 

\begin{thm}
Let $M$ be the AFD factor of type \II,  and 
$\bG=(\lG,\de,h)$  a compact Kac algebra with  amenable dual. 
Let $\al$ and $\be$ be minimal actions of $\bG$ on $M$. 
Then they are conjugate. 
\end{thm}

\begin{proof}
By the previous corollary, 
a minimal action of $\bG$ on the AFD factor of type \II\  is dual, 
and $\al$ and $\be$ are of the form $\hat{\ga_0}$ and $\hat{\ga_1}$ 
where $\ga_0$ and $\ga_1$ are free actions of $\bhG$ on $M^\al$ and 
$M^\be$ respectively. 
Since $M^\al$ and $M^\be$ are injective factors of type \II, they are 
isomorphic by Connes's result \cite{Co}. 
By Corollary \ref{cor: uniqueness},
$\ga_0$ and $\ga_1$ are cocycle conjugate. 
Hence their dual actions $\al$ and $\be$ are conjugate. 
\end{proof}

\section{Appendix}

Let $M$ be a von Neumann algebra, and $K$ a finite dimensional 
Hilbert space. 
For $\be\in \Mor_0(M, M\oti B(K))$, 
we prepare several properties. 

\begin{defn}
We say that $\be$ is 
\begin{enumerate}

\item
\textit{properly outer} if there exists no nonzero $a\in M\oti B(K)$ 
such that 
$\be(x)a=a(x\oti1)$ for all $x\in M$, 

\item
\textit{centrally trivial} if $\be^\om(x)=x\oti1$ for all $x\in M_\om$,

\item
\textit{centrally nontrivial} if $\be$ is not centrally trivial,

\item
\textit{properly centrally nontrivial} if 
there exists no nonzero element $a\in M\oti B(K)$ such that 
$\be^\om(x)a=(x\oti1)a$ for all $x\in M_\om$. 

\end{enumerate}
\end{defn}

\begin{lem}
A map $\be\in \Mor_0(M,M\oti B(K))$ is properly centrally nontrivial 
if and only if 
it is strongly outer. 
\end{lem}

\begin{proof}
The ``if'' part is trivial. 
We show the ``only if'' part. 
Let $\be$ be a properly centrally nontrivial homomorphism. 
Assume that $\be$ is not strongly outer. 
Then there exists a nonzero 
$a\in M^\om \oti B(K)$ and a countably 
generated von Neumann algebra $S\subs M^\om$ 
such that $\be^\om(x)a=a(x\oti1)$ holds for all $x\in S'\cap M_\om$. 
We claim that the proper central nontriviality implies 
\[
1=\bigvee_{z\in M_\om} 
s\big{(}
(\ta^\om\oti\id)(|\be^\om(z)-z\oti 1|^2)
\big{)}.
\]
Indeed, if $b\in M\oti B(K)$ satisfies 
$(\ta^\om\oti\id)(|\be^\om(z)-z\oti 1|^2)b=0$ for all $z\in M_\om$, 
then 
$(\be^\om(z)-z\oti 1)b=0$, 
but this is a contradiction. 
Hence there exists $z\in M_\om$ 
such that it satisfies 
$(\ta^\om\oti\id)(|a^*|^2)(\ta^\om\oti\id)(|\be^\om(z)-z\oti 1|^2)
\neq 0$. 
Let $\{e_{i,j}\}_{i,j=1}^n$ be a system of matrix units 
for $B(K)$. 
Decompose $a$ as $a=\sum_{i,j=1}^n a_{i,j}\oti e_{i,j}$. 
Set 
$N=W^*(z)$, $\wdt{S}=W^*(S,\{a_{i,j}\}_{1\leq i,j\leq n})$ 
and $\mB=\{\be^\om\}$. 
Then by Lemma \ref{lem: fast}, 
there exists $\Ps\in \Mor(\wdt{N}, M^\om)$ such that 
\begin{enumerate}

\item 
$\Ps(z)\in \wdt{S}'\cap M_\om$, 

\item 
$\ta^\om(b \Ps(z))=\ta^\om(b)\ta^\om(z)$ for all 
$b\in \wdt{S}$, 

\item $\be^\om(\Ps(z))=(\Ps\oti\id)(\be^\om(z))$, 

\end{enumerate}
Set $y=\Ps(z)\in \wdt{S}'\cap M_\om$ and 
then $\be^\om(y)a=a(y\oti1)=(y\oti1)a$. 
However the following equality 
\begin{align*}
(\ta^\om\oti\id)(|a^*|^2|\be^\om(y)-y\oti1|^2)
=&\,
(\ta^\om\oti\id)\big{(}|a^*|^2(\Ps\oti\id)(|\be^\om(z)-z\oti1|^2)\big{)}
\\
=&\,
(\ta^\om\oti\id)(|a^*|^2)(\ta^\om\oti\id)(|\be^\om(z)-z\oti1|^2)
\\
\neq&\,0
\end{align*}
implies $(\be^\om(y)-(y\oti1))a\neq0$. 
This is a contradiction. 
\end{proof}
Hence $\be$ is 
\begin{equation}\label{eq: proper}
\mbox{strongly outer}
\Leftrightarrow
\mbox{properly centrally nontrivial}
\Rightarrow
\mbox{centrally nontrivial}.
\end{equation}

A map $\be\in \Mor(M, M\oti B(K))$ is 
said to be \textit{irreducible} if 
$\be(M)'\cap (M\oti B(K))=\C$. 

\begin{lem}\label{lem: cent-proper}
Let $\be\in \Mor_0(M,M\oti B(K))$ be irreducible. 
Then 
$\be$ is centrally nontrivial 
if and only if 
$\be$ is properly centrally nontrivial. 
\end{lem}

\begin{proof}
The ``if'' part is trivial. 
We show the ``only if'' part. 
Let $\be$ be a centrally nontrivial homomorphism. 
Assume that there exists $a\in M\oti B(K)$ 
satisfying 
$\be^\om(x)a=(x\oti1)a$ for all $x\in M_\om$. 
Let $a=v|a|$ be the polar decomposition. 
Then it is easy to see that 
$\be^\om(x)vv^*=(x\oti1)vv^*$ for all $x\in M_\om$. 
Put $p=vv^*\in M\oti B(K)$. 
If $u\in U(M)$, then 
$\be^\om(x)\be(u)p\be(u^*)=(x\oti1)\be(u)p\be(u^*)$. 
Therefore the projection 
$z=\vee_{u\in U(M)}\be(u)p\be(u^*)$ satisfies 
$\be^\om(x)z=(x\oti1)z$ for all $x\in M_\om$. 
Since $z\in \be(M)'\cap M\oti B(K)=\C$, we have $z=1$, 
and $\be$ is properly centrally nontrivial. 
\end{proof}
Hence 
under the assumption on irreducibility, 
all the properties of (\ref{eq: proper}) are equivalent. 
In addition, if $M$ is the AFD factor of type \II, 
they are equivalent to proper outerness. 

\begin{lem}\label{lem: cent-proper-outer}
Let $\meR_0$ be the AFD factor of type \II. 
Let $\be\in\Mor_0(\meR_0,\meR_0\oti B(K))$ be irreducible. 
Then the following properties on $\be$ are equivalent: 
\begin{enumerate}

\item central nontriviality,

\item proper central nontriviality,

\item strong outerness, 

\item proper outerness.

\end{enumerate}
\end{lem}

\begin{proof}
We know the equivalence of $(1)$, $(2)$ and $(3)$. 
It is trivial that $(3)$ implies $(4)$. 
We show that $(4)$ implies $(1)$. 
Let $\ta$ be the trace on $\meR_0$. 
We assume the following lemma for a moment, and 
we prove this implication. 
Let $\meK_1\subs\meK_2\subs\dots$ be an increasing net 
of finite dimensional subfactors in $\meR_0$ with 
$(\cup_{n\geq1}\meK_n)''=\meR_0$. 
For each $n$, there exists a unitary $u_n\in \meK_n'\cap \meR_0$ 
with $\|\be(u_n)-u_n\oti1\|_{\ta\oti\ta_K}>1/2$. 
Then the sequence $(u_n)_n$ is central which defines $u$ in $(\meR_0)_\om$. 
It satisfies $\|\be^\om(u)-u\oti1\|_{\ta\oti\ta_K}\geq1/2$. 
Hence $\be^\om$ is not trivial on $(\meR_0)_\om$. 
\end{proof}

We adapt \cite[Lemma 3.4]{Co-peri} 
to the case of a homomorphism as follows. 

\begin{lem}
Let $M$ be a factor of type II$_1$ and $K$ a finite dimensional 
Hilbert space. 
Let $\be\in\Mor(M, M\oti B(K))$ be irreducible. 
If there exists a finite dimensional subfactor $\meK\subset M$ with 
\[
\sup\{\|\be(u)-u\oti1\|_{\ta\oti\ta_K}
\mid u\in U(\meK'\cap M)\}<1,
\]
then $\be$ is not properly outer. 
\end{lem}

\begin{proof}
Consider a weakly closed convex set in $M\oti B(K)$, 
$C=\ovl{\co}^w\{(u\oti1)\be(u^*)\mid u\in U(\meK'\cap M)\}$. 
In $C$ 
take a unique point $y_0$ attaining the minimal distance from 0 with
 repect to $\|\cdot\|_\tau$. 
By assumption, $\|y_0-1\|_\ta<1$, in particular, $y_0\neq0$. 
Unicity yields $(u\oti1)y_0=y_0\be(u)$ for any 
$u\in U(\meK'\cap M)$. 
Hence we have $(x\oti1)y_0=y_0\be(x)$ for any 
$x\in \meK'\cap M$. 
Let $\{e_{i,j}\}_{i,j=1}^n$ be a system of matrix units for $\meK$. 
Since $\be(e_{1,1})$ and $e_{1,1}\oti1$ are equivalent 
in $M\oti B(K)$, 
there exists a partial isometry $v\in M\oti B(K)$ with 
$\be(e_{1,1})=vv^*$ and $e_{1,1}=v^*v$. 
Set a unitary 
$
u=\sum_{i=1}^n \be(e_{i,1})v e_{1,i}
$ in $M\oti B(K)$. 
Then we have $u(e_{i,j}\oti1)=\be(e_{i,j})u$. 
Hence $u(x\oti1)=\be(x)u$ for all $x\in \meK$. 
Let $\ga=\Ad(u^*)\circ \be$ and $y_1=y_0 u$. 
Then $\ga$ is trivial on $\meK$ and 
$y_1 \ga(x)=(x\oti1)y_1$ holds for $x\in \meK'\cap M$. 
Let $E\col M\oti B(K)\ra (\meK'\cap M)\oti B(K)$ be a 
faithful conditional expectation. 
Then there exists an element $a$ in $\meK$ with 
$z=E((a\oti1)y_1)\neq0$. 
Since $\ga(\meK'\cap M)\subs (\meK'\cap M)\oti B(K)$, 
$z$ satisfies 
$z \ga(x)=(x\oti1)z$ for $x\in \meK'\cap M$. 
In fact, this equality is valid for any $x\in M$ 
because $\ga$ is trivial on $\meK$. 
This shows $\ga$ is not outer. 
Hence $\be=\Ad(u)\circ \ga$ is not outer. 
\end{proof}

Let $(\al,u)$ be a cocycle action of $\bhG$ 
on a von Neumann algebra $M$. 
We call $(\al,u)$ 
\textit{centrally free} 
if $\al_\pi$ is properly centrally nontrivial 
for each $\pi\in\IG\setm\{\btr\}$. 

\begin{cor}\label{cor: cent-strong-free}
Let $(\al,u)$ be a cocycle action of $\bhG$ 
on a von Neumann algebra $M$. 
Then the following properties of $(\al,u)$ are equivalent: 
\begin{enumerate}

\item
central freeness,

\item 
strong freeness.
\end{enumerate}
In addition, if $M$ is the AFD factor of type \II, 
they are also equivalent to 
\\
\hspace{15pt}(3) freeness. 
\end{cor}

\begin{proof}
We know that (1) and (2) are equivalent and 
(2) implies (3). 
We show that (3) implies (1). 
By Lemma \ref{lem: relcom}, 
each map $\al_\pi$, $\pi\in\IG$, is irreducible. 
Then by Lemma \ref{lem: cent-proper-outer}, 
$\al_\pi$ is properly centrally nontrivial. 
\end{proof}

\end{document}